\documentclass[a4paper,12pt]{article}
\usepackage[cp1251]{inputenc}
\usepackage{amssymb,amsmath,amsthm}
\usepackage[russian]{babel}
\usepackage{indentfirst}
\usepackage{amssymb,amsmath,amsthm}
\usepackage{euler}
\usepackage{hyperref,hypbmsec}
\newtheorem{Th}{\hskip\parindent Theorem}
\newtheorem{Le}{\hskip\parindent Lemma}
\newtheorem{Utv}{\hskip\parindent Statement}

\newtheorem{Sl}{\hskip\parindent Corollary}
\newtheorem{Zam}{\hskip\parindent Remark}

\newcommand{\Q}{Q’}

\newcounter{propet}

\renewcommand{\le}{\leqslant}\renewcommand{\ge}{\geqslant}

\begin{document}
\author{D.\,A.\,Frolenkov\footnote{The research was supported by the grant RFBR № 11-01-00759-а} }
\title{
The mean value of Frobenius numbers with three arguments.
}
\date{}
\maketitle

\setcounter{Zam}0
This paper is written in Russian but we translate the statement of our main result in English. \par

The Frobenius number $g(a_1,…,a_n)$ of relatively prime positive integers $a_1,…,a_n$ is defined as the largest number $k$ that is not representable as a non-negative integer combination of  $a_1,…,a_n$
\begin{equation}\label{defg}
x_1a_1+…+x_na_n=k,\qquad\quad x_1\ge 0,…,x_n\ge 0.
\end{equation}
Sometimes it is easier to use the following function
\begin{equation}\label{deff}
f(a_1,…,a_n)=g(a_1,…,a_n)+a_1+…+a_n
\end{equation}
This function may be defined as the largest number $k,$ that is not representable as a positive integer combination of  $a_1,…,a_n$
Some general results on the Frobenius problem are available in the book ~\cite{Ramirez}. A.V. Ustinov proved (see ~\cite{Ustinov-frob}) that the average value of $f(a,b,c)$ is equal to $\frac{8}{\pi}\sqrt{abc}$.
\begin{Th}
\label{Ustinov-frob-1}
Let $a$  be a positive integer,and let $x_1>0,x_2>0,\varepsilon>0$  be real numbers. Then the following asymptotic formula is valid
\begin{equation}\label{sred2}
\frac{1}{a^{3/2}|M_a(x_1,x_2)|}\sum_{(b,c)\in M_a(x_1,x_2)}\left(f(a,b,c)-\frac{8}{\pi}\sqrt{abc}\right)=
O_\varepsilon(R_\varepsilon(a;x_1,x_2)),
\end{equation}
where
\begin{gather}
R_\varepsilon(a;x_1,x_2)=\left(a^{-1/6}(x_1+x_2)+a^{-1/4}(x_1^{3/2}+x_2^{3/2})(x_1x_2)^{-1/4}+a^{-1/2}\right)a^\varepsilon\notag\\
\ll_{x_1,x_2}a^{-1/6+\varepsilon}\label{sred2ost}
\end{gather}
and
$$M_a(x_1,x_2)=\left\{(b,c): 1\le b\le x_1a, 1\le c\le x_2a, (a,b,c)=1\right\}$$
\end{Th}

\begin{Zam}
With the help of Theorem \ref{Ustinov-frob-1} it can be shown that
\begin{gather*}
\frac{1}{x_1x_2x_3N^{9/2}}\sum_{a\le x_1N}\sum_{b\le x_2N}\sum\limits_{c\le x_3N \atop (a,b,c)=1}
\left(f(a,b,c)-\frac{8}{\pi}\sqrt{abc}\right)=O_{\varepsilon,x_1,x_2,x_3}(N^{-1/6+\varepsilon}).
\end{gather*}
\end{Zam}
In this paper we prove the following result.
\begin{Th}\label{Th1}
The following asymptotic formula is valid
\begin{gather*}
\frac{1}{x_1x_2x_3N^{9/2}}\sum_{a\le x_1N}\sum_{b\le x_2N}\sum\limits_{c\le x_3N \atop (a,b,c)=1}
\left(f(a,b,c)-\frac{8}{\pi}\sqrt{abc}\right)=
O\left(N^{-1/2+\varepsilon}\left(
x_1^{1+\varepsilon}\frac{x_3^{\varepsilon}}{x_2^{\varepsilon}}+
x_2^{1+\varepsilon}\frac{x_3^{\varepsilon}}{x_1^{\varepsilon}}
+x_3^{1+\varepsilon}
\right)\right),
\end{gather*}
where $x_1>0,x_2>0, x_3>0, \varepsilon>0$ are real numbers.
\end{Th}
This theorem was conjectured by A.V. Ustinov in ~\cite{Ustinov-frob}.
To prove Theorem \ref{Th1} we use the ideas from the Ustinov's papers  ~\cite{Ustinov-frob} and ~\cite{Ustinov2007} and from the E.N.Zhabitskaya's paper~\cite{Jabizkay}. We also use classical bounds on exponential sums.

I am very grateful to my advisor, N.G.Moshchevitin, for helpful discussions. I would also like to thank I.D.Kan for meticulously reading and commenting this paper.
\newpage
\section{Введение}
Числом Фробениуса $g(a_1,…,a_n)$ натуральных чисел $a_1,…,a_n$, взаимно простых в совокупности, называется наибольшее целое $k,$ не представимое в виде суммы
\begin{equation}\label{defg}
x_1a_1+…+x_na_n=k,\quad\mbox{где}\quad x_1\ge 0,…,x_n\ge 0.
\end{equation}
Во многих задачах оказывается удобнее рассматривать функцию
\begin{equation}\label{deff}
f(a_1,…,a_n)=g(a_1,…,a_n)+a_1+…+a_n
\end{equation}
равную наибольшему целому $k,$ не представимому в виде суммы  \eqref{defg}, но уже с натуральными коэффициентами $x_1,…,x_n.$
Наиболее обширный обзор результатов и задач, связанных с числом Фробениуса, приведен в книге ~\cite{Ramirez}. А.В. Устиновым в работе ~\cite{Ustinov-frob} было доказано, что функция $f(a,b,c)$ в среднем ведет себя как $\frac{8}{\pi}\sqrt{abc}$.
\begin{Th}
\label{Ustinov-frob-1}
Пусть a~---натуральное число, $x_1>0,x_2>0,\varepsilon>0$~---действительные числа. Тогда
\begin{equation}\label{sred2}
\frac{1}{a^{3/2}|M_a(x_1,x_2)|}\sum_{(b,c)\in M_a(x_1,x_2)}\left(f(a,b,c)-\frac{8}{\pi}\sqrt{abc}\right)=
O_\varepsilon(R_\varepsilon(a;x_1,x_2)),
\end{equation}
где
\begin{gather}
R_\varepsilon(a;x_1,x_2)=\left(a^{-1/6}(x_1+x_2)+a^{-1/4}(x_1^{3/2}+x_2^{3/2})(x_1x_2)^{-1/4}+a^{-1/2}\right)a^\varepsilon\notag\\
\ll_{x_1,x_2}a^{-1/6+\varepsilon}\label{sred2ost}
\end{gather}
и
$$M_a(x_1,x_2)=\left\{(b,c): 1\le b\le x_1a, 1\le c\le x_2a, (a,b,c)=1\right\}$$
\end{Th}

\begin{Zam}
Используя теорему \ref{Ustinov-frob-1}, легко получить оценку
\begin{gather*}
\frac{1}{x_1x_2x_3N^{9/2}}\sum_{a\le x_1N}\sum_{b\le x_2N}\sum\limits_{c\le x_3N \atop (a,b,c)=1}
\left(f(a,b,c)-\frac{8}{\pi}\sqrt{abc}\right)=O_{\varepsilon,x_1,x_2,x_3}(N^{-1/6+\varepsilon}).
\end{gather*}
\end{Zam}
Сформулируем основной результат данной работы.
\begin{Th}\label{Th1}
Справедлива оценка
\begin{gather*}
\frac{1}{x_1x_2x_3N^{9/2}}\sum_{a\le x_1N}\sum_{b\le x_2N}\sum\limits_{c\le x_3N \atop (a,b,c)=1}
\left(f(a,b,c)-\frac{8}{\pi}\sqrt{abc}\right)=
O\left(N^{-1/2+\varepsilon}\left(
x_1^{1+\varepsilon}\frac{x_3^{\varepsilon}}{x_2^{\varepsilon}}+
x_2^{1+\varepsilon}\frac{x_3^{\varepsilon}}{x_1^{\varepsilon}}
+x_3^{1+\varepsilon}
\right)\right),
\end{gather*}
где $x_1>0,x_2>0, x_3>0, \varepsilon>0$~---действительные числа.
\end{Th}
Это утверждение было сформулировано А.В. Устиновым в работе ~\cite{Ustinov-frob} в виде гипотезы.
Доказательство теоремы \ref{Th1} использует идеи из работ А.В. Устинова ~\cite{Ustinov-frob} и ~\cite{Ustinov2007} и работы Е.Н. Жабицкой ~\cite{Jabizkay}. Также используются классические оценки тригонометрических сумм.

Автор благодарен И.Д.Кану за указания на недочеты, имевшиеся в первоначальной версии статьи и глубоко признателен Н.Г. Мощевитину за неоднократные обсуждения полученных результатов.
\section{Вспомогательные утверждения и обозначения}
Разложим рациональное число $r$ в стандартную цепную дробь
\begin{equation}
r=[0;a_1\ldots,a_s]=\cfrac{1}{a_1+{\atop\ddots\,\displaystyle{+\cfrac{1}{a_s}}}}
\end{equation}
длины $s=s(r)$, в которой $a_1$, \ldots,$a_s$~--- натуральные
 и $a_s\ge 2$ при $s\ge1$.
Через $s_{1}(r)$ будем обозначать  сумму неполных частных
$$s_{1}(r)=\sum\limits_{1\le i \le s}a_i.$$
\begin{Le}
\label{s1}
Для любого натурального $b>1$ выполнено
$$\sum_{a\le b}s_1(a/b)\ll b\log^2 b.$$
\begin{proof} См в работе Д.Кнута ~\cite{Knuth}. \end{proof}
\end{Le}
\begin{Le}
\label{tau}
Пусть $\tau(n)=\sum\limits_{d|n}1$---число делителей натурального числа n, тогда
$$\tau(n)=o(n^{\varepsilon})$$
для любого $\varepsilon>0$.
\begin{proof} См книгу  К. Чандрасекхарана ~\cite[гл 6, \S3, теорема 5]{Chandra}. \end{proof}
\end{Le}
Следующая лемма общеизвестна (преобразование Абеля)
\begin{Le}
\label{abel}
Пусть $f(x)$---непрерывно-дифференцируема на $[a;b],$ $c_n$---произвольные числа,
$$C(x)=\sum_{a<n\le x}c_n.$$
Тогда
$$
\sum_{a<n\le b}c_nf(n)=C(b)f(b)-\int\limits_a^bC(x)f’(x)dx.
$$
\begin{proof} См книгу А.А.Карацубы ~\cite[гл. 2, \S5]{Karazuba}. \end{proof}
\end{Le}

\begin{Le}
\label{trigsum}
Пусть $\alpha$---произвольное действительное число, $Q$---целое и $P$---натуральное. Тогда
$$
\left|\sum_{x=Q+1}^{Q+P}\exp(2\pi i\alpha x)\right|\le\min\left(P,\frac{1}{2\|\alpha\|}\right),
$$
где $\|\alpha\|$ ---расстояние от $\alpha$ до ближайшего целого.
\begin{proof} См книгу Н.М. Коробова ~\cite[гл. 1, \S1]{Korobov}. \end{proof}
\end{Le}
\begin{Sl}
\label{trigsum12}
Пусть $\alpha$---произвольное действительное число.
Если $f(x)$ непрерывно дифференцируема на отрезке $[Q,Q+P]$ и монотонна, то
 $$
\left|\sum_{x=Q+1}^{Q+P}\exp(2\pi i\alpha x)f(x)\right|\ll
\left(\left|f(P+Q)\right|+\left|f(Q)\right|\right)\min\left(P,\frac{1}{\|\alpha\|}\right).
$$
\begin{proof}
После применения леммы \ref{abel} и леммы \ref{trigsum}, получаем
\begin{gather*}
\left|\sum_{x=Q+1}^{Q+P}\exp(2\pi i\alpha x)f(x)\right|\ll
\left(\left|f(P+Q)\right|+\int\limits_Q^{Q+P}|f’(x)|dx\right)\min\left(P,\frac{1}{\|\alpha\|}\right).
\end{gather*}
Используя монотонность функции $f(x)$, получаем нужную оценку.
\end{proof}
\end{Sl}

\begin{Le}
\label{dozelogo}
Пусть $\tau,\, p$ ---произвольные натуральные числа, a и b ---действительные. Тогда
\begin{enumerate}
  \item[(a)]
$$
\sum_{\substack{a<n\le b\\\tau\nmid n}}\frac{1}{\|\frac{n}{\tau}\|}\ll (b-a)\log\tau+\tau\log\tau,
$$
\item[(b)]
$$
\sum_{\substack{a<n\le b\\\tau\nmid n}}\frac{1}{n\|\frac{n}{\tau}\|}\ll
\log\frac{b}{a}\log\tau+\log\tau+\frac{\tau}{a}\log\tau,
$$
\item[(c)]
$$
\sum_{\substack{a<n\le b\\\tau\nmid n}}\frac{1}{n^2\|\frac{n}{\tau}\|}\ll
\frac{1}{a}\log\tau+\frac{\tau}{a^2}\log\tau,
$$
\item[(d)]
$$
\sum_{\substack{a<n\le b\\\tau\nmid n}}\frac{n^p}{\|\frac{n}{\tau}\|}\ll b^p(b-a)\log\tau+b^p\tau\log\tau.
$$
\end{enumerate}
\begin{proof}
Доказательство (a) можно найти в книге \mbox{Н.М. Коробова ~\cite[гл. 1, \S1]{Korobov}.} Для доказательства остальных пунктов достаточно воспользоваться леммой \ref{abel}.
\end{proof}
\end{Le}

Обозначим, следуя Н.М. Коробову, через $\delta_q(a)$ характеристическую функцию делимости на натуральное число q
\begin{gather}\label{deltasum}
\delta_q(a)=\frac{1}{q}\sum_{x=1}^q\exp\left(2\pi i\frac{ax}{q}\right)=
\left\{
              \begin{array}{ll}
                1, & \hbox{если $a\equiv 0 \pmod{q}$;} \\
                0, & \hbox{иначе.}
              \end{array}
\right.
\end{gather}
В дальнейшем нам понадобятся следующие обозначения. \par
Знак звездочки в суммах вида
$$
\mathop{{\sum}^*}_{x=1}^{a}, \quad
\sum_{a}\mathop{{\sum}^*}_{x}
$$
означает, что суммирование ведется по числам, удовлетворяющим условию $(a,x)=1$.
В суммах вида
$$\sum_{d|n}$$
суммирование ведется по делителям числа $n$.
В суммах вида
$$\sum_{n\le R\atop d|n},\qquad \sum_{n\le R\atop d\nmid n}$$
суммирование ведется по $n$, удовлетворяющим условию
$$n\equiv 0 \pmod{d},\qquad n\not\equiv 0 \pmod{d}$$
соответственно. Если $A$~--- некоторое утверждение, то $[A]$ означает $1$, если $A$ истинно, и $0$ в противном случае.\par

\section{О функции Редсета}
В работе ~\cite{Rodset} Редсет доказал следующий метод для подсчета функции $f(a,b,c)$.
Пусть $a,b,c$ ---натуральные числа и $(a,b)=1$, $(a,c)=1$, тогда существует натуральное \mbox{число $l$}
$$c\equiv bl \pmod{a}, \, 1\le l\le a, \,(l,a)=1.$$
Разложим число $\frac{a}{l}$ в приведенную регулярную  цепную дробь (см. например ~\cite[гл. 1]{Perron}.)
\begin{equation}
\frac{a}{l}=\langle a_1;a_2\ldots,a_m\rangle=a_1-\cfrac{1}{a_2-{\atop\ddots\,\displaystyle{-\cfrac{1}{a_m}}}}
\end{equation}
и определим последовательности $\{s_j\}\, \mbox{и}\,\{q_j\}, \mbox{при} -1\le j\le m$ следующим образом
$$
s_m=0,\quad  s_{m-1}=1, \quad q_{-1}=0,\quad  q_0=1
$$
$$
s_{j-1}=a_{j+1}s_{j}-s_{j+1}, \quad  q_{j+1}=a_{j+1}q_{j}-q_{j-1},\quad 0\le j\le m-1.
$$
Легко доказать (см ~\cite{Rodset} или ~\cite{Ustinov-frob}) что, для последовательностей $\{s_j\}$ $\{q_j\}$ выполнено
$$
0=\frac{s_m}{q_m}<\frac{s_{m-1}}{q_{m-1}}<\cdots<\frac{s_0}{q_0}<\frac{s_{-1}}{q_{-1}}=\infty.
$$
Функция Рёдсета $\rho_{l,a}(t_1,t_2)$ для $t_1\ge 0, t_2\ge 0$ таких, что
$$
\frac{s_n}{q_n}\le \frac{t_2}{t_1}<\frac{s_{n-1}}{q_{n-1}}
$$
задается формулой
\begin{gather*}
\rho_{l,a}(t_1,t_2)=t_1s_{n-1}+t_2q_{n}-\min\{t_1s_{n},t_2q_{n-1}\},
\end{gather*}
тогда функция $f(a,b,c)$, определенная в \eqref{deff}, находится по формуле
\begin{equation}
f(a,b,c)=\rho_{l,a}(b,c).
\end{equation}
Так же нам понадобится следующее утверждение о последовательностях \mbox{$\{s_j\}$, $\{q_j\}$ (см ~\cite{Ustinov-frob}).}
\begin{Utv}\label{sv5}
Четверки $(q_n,s_{n-1},q_{n-1},s_n)$ при $0\le n\le m$ находятся во взаимно однозначном соответствии с решениями $(u_1,u_2,v_1,v_2)$ уравнения
$$u_1u_2-v_1v_2=a,$$
для которых
$$
0\le v_1<u_1\le a, \quad (u_1,v_1)=1,\quad 0\le v_2<u_2\le a,\quad  (u_2,v_2)=1.
$$
\end{Utv}
Определим функцию $\rho_{l,a}(\alpha)$ следующим образом
$$
\rho_{l,a}(\alpha)=s_{n-1}+\alpha q_n-\min\{s_n,\alpha q_{n-1}\}, \quad \frac{s_n}{q_n}\le \alpha<\frac{s_{n-1}}{q_{n-1}}.
$$
Тогда в силу однородности функции Рёдсета получим
\begin{equation}\label{1}
\rho_{l,a}(t_1,t_2)=t_1\rho_{l,a}\left(\frac{t_2}{t_1}\right).
\end{equation}
В дальнейшем нам понадобятся еще две функции
\begin{equation}\label{2}
\rho_{a}^{*}(t_1,t_2)=\mathop{{\sum}^*}_{l=1}^{a}\rho_{l,a}(t_1,t_2), \quad
\rho_{a}^{*}(\alpha)=\mathop{{\sum}^*}_{l=1}^{a}\rho_{l,a}(\alpha).
\end{equation}
Очевидно, что
\begin{equation}\label{3}
\rho_{a}^{*}(t_1,t_2)=t_1\rho_{a}^{*}\left(\frac{t_2}{t_1}\right).
\end{equation}
Используя утверждение \ref{sv5}, можно получить следующую формулу (см.~\cite[\S 6]{Ustinov-frob})
\begin{equation}\label{4}
\rho_{a}^{*}(\alpha)=\lambda^{*}(a,\alpha)+\alpha\lambda^{*}\left(a,\frac{1}{\alpha}\right),
\end{equation}
где
\begin{equation}
\lambda^{*}(a,\alpha)=\sum_{x=1}^a\mathop{{\sum}^*}_{z=1}^x\sum_{y=1}^a\mathop{{\sum}^*}_{w=0}^{a-1}
\left[xy+wz=a, \frac{w}{x}\le \alpha <\frac{w}{x-z}\right](y+w+\alpha z).
\end{equation}
Освобождаясь от условий взаимной простоты, получаем
\begin{equation}\label{ljmbda*}
\lambda^{*}(a,\alpha)=\sum_{d_1d_2|a}\mu(d_1)\mu(d_2)d_2\lambda\left(\frac{a}{d_1d_2},\frac{d_1\alpha}{d_2}\right),
\end{equation}
где
\begin{equation}\label{ljmbda}
\lambda(a,\alpha)=\sum_{x\ge 1}\sum_{z=1}^x\sum_{y\ge 1}^a\sum_{w\ge 0}
\left[xy+wz=a, \frac{w}{x}\le \alpha <\frac{w}{x-z}\right](y+w+\alpha z).
\end{equation}

\section{Выделение плотности}
Этот параграф является переработкой соответствующего параграфа из работы Устинова (см.~\cite[\S 5]{Ustinov-frob}) с учетом того, что теперь усреднение ведется по трем переменным. Обозначим за F и G следующие суммы
\begin{equation}\label{FG}
F=\sum_{a\le x_3N}\sum_{b\le x_1N}\sum\limits_{c\le x_2N \atop (a,b,c)=1}f(a,b,c),\quad
G=\sum_{a\le x_3N}\sum_{b\le x_1N}\sum\limits_{c\le x_2N \atop (a,b,c)=1}\frac{8}{\pi}\sqrt{abc}.
\end{equation}
\begin{Le}
\label{F}
Пусть $x_1>0,x_2>0,\varepsilon>0$~---действительные числа. Тогда справедливо следующее равенство
\begin{gather}\nonumber
F=
\sum\limits_{d_1d_2\le x_3N \atop (d_1,d_2)=1}d_1d_2\sum_{\delta_1\le \frac{x_3N}{d_1}}\sum_{\delta_2\le \frac{x_3N}{d_2}}
\frac{\mu(\delta_1)}{\delta_1}\frac{\mu(\delta_2)}{\delta_2}(\delta_1,\delta_2)
\int\limits_0^{\frac{x_1N}{d_1}}\int\limits_0^{\frac{x_2N}{d_2}}\sum\limits_{a\le \frac{x_3N}{d_1d_2} \atop \delta|a}
\frac{t_1\lambda^*(a,\frac{t_2}{t_1})+t_2\lambda^*(a,\frac{t_1}{t_2})}{a} dt_1dt_2+\\\nonumber+
O\left(x_1x_2x_3^{2+\varepsilon}N^{4+\varepsilon}\right),
\end{gather}
где
$$\delta=HOK\left(\frac{\delta_1}{(\delta_1,d_2)},\frac{\delta_2}{(\delta_2,d_1)}\right).$$
\begin{proof}
Обозначим $$d_1=(a,b), d_2=(a,c), a_1=\frac{a}{(d_1d_2)}, b_1=\frac{b}{d_1}, c_1=\frac{c}{d_2}.$$ Так как $(a,b,c)=1$, то $(d_1,d_2)=1$. Получаем
\begin{gather*}
F=\sum_{a\le x_3N}\sum\limits_{d_1d_2|a \atop (d_1,d_2)=1}\sum\limits_{b\le x_1N \atop (b,a)=d_1}
\sum\limits_{c\le x_2N \atop (c,a)=d_2}f(d_1d_2a_1,d_1b_1,d_2c_1)=\\
\sum\limits_{d_1d_2\le x_3N \atop (d_1,d_2)=1}\sum_{a_1\le \frac{x_3N}{d_1d_2}}
\sum\limits_{b_1\le \frac{x_1N}{d_1} \atop (b_1,a_1d_2)=1}
\sum\limits_{c_1\le \frac{x_2N}{d_2} \atop (c_1,a_1d_1)=1}d_1d_2f(a_1,b_1,c_1).
\end{gather*}

В последнем равенстве мы воспользовались тождеством Джонсона (см ~\cite{Johnson})
$$f(da,db,c)=df(a,b,c).$$
Выражая функцию $f(a,b,c)$ через функцию Рёдсета, получим
\begin{gather*}
F=\sum\limits_{d_1d_2\le x_3N \atop (d_1,d_2)=1}d_1d_2\sum_{a_1\le \frac{x_3N}{d_1d_2}}
\sum\limits_{b_1\le \frac{x_1N}{d_1} \atop (b_1,a_1d_2)=1}
\sum\limits_{c_1\le \frac{x_2N}{d_2} \atop (c_1,a_1d_1)=1}\mathop{{\sum}^*}_{l=1}^{a_1}
\delta_{a_1}(b_1l-c_1)\rho_{l,a_1}(b_1,c_1)=\\=
\sum\limits_{d_1d_2\le x_3N \atop (d_1,d_2)=1}d_1d_2\sum_{a_1\le \frac{x_3N}{d_1d_2}}
\sum_{\delta_1|d_2a_1}\mu(\delta_1)\sum_{\delta_2|d_1a_1}\mu(\delta_2)\mathop{{\sum}^*}_{l=1}^{a_1}
\sum\limits_{b_1\le \frac{x_1N}{d_1} \atop \delta_1|b_1}
\sum\limits_{c_1\le \frac{x_2N}{d_2} \atop \delta_2|c_1}\delta_{a_1}(b_1l-c_1)\rho_{l,a_1}(b_1,c_1).
\end{gather*}
Используя лемму 2 из ~\cite[\S 5]{Ustinov-frob}, получаем
\begin{gather*}
S_{l,a_1}=\sum\limits_{b_1\le \frac{x_1N}{d_1} \atop \delta_1|b_1}\sum\limits_{c_1\le \frac{x_2N}{d_2} \atop \delta_2|c_1}
\delta_{a_1}(b_1l-c_1)\rho_{l,a_1}(b_1,c_1)=\\=\frac{(a_1,\delta_1,\delta_2)}{a_1\delta_1\delta_2}
\int_0^{\frac{x_1N}{d_1}}\int_0^{\frac{x_2N}{d_2}}\rho_{l,a_1}(t_1,t_2)dt_1dt_2+
O\left(\frac{x_1x_2}{d_1d_2}N^2s_1\left(\frac{l}{a_1}\right)\right).
\end{gather*}
Следовательно,
\begin{gather}
F=\sum\limits_{d_1d_2\le x_3N \atop (d_1,d_2)=1}d_1d_2\sum_{a_1\le \frac{x_3N}{d_1d_2}}
\sum_{\delta_1|d_2a_1}\mu(\delta_1)\sum_{\delta_2|d_1a_1}\mu(\delta_2)\frac{(a_1,\delta_1,\delta_2)}{a_1\delta_1\delta_2}
\int_0^{\frac{x_1N}{d_1}}\int_0^{\frac{x_2N}{d_2}}\mathop{{\sum}^*}_{l=1}^{a_1}\rho_{l,a_1}(t_1,t_2)dt_1dt_2+\notag\\+
O\left(x_1x_2N^2\sum\limits_{d_1d_2\le x_3N \atop (d_1,d_2)=1}\sum_{a_1\le \frac{x_3N}{d_1d_2}}
\sum_{\delta_1|d_2a_1}\sum_{\delta_2|d_1a_1}\mathop{{\sum}^*}_{l=1}^{a_1}s_1\left(\frac{l}{a_1}\right)\right).\label{FF}
\end{gather}
Для оценки остаточного члена воспользуемся леммой \ref{s1} и леммой \ref{tau}, получаем
\begin{gather*}
O\left(x_1x_2N^2\sum\limits_{d_1d_2\le x_3N \atop (d_1,d_2)=1}\sum_{a_1\le \frac{x_3N}{d_1d_2}}
\sum_{\delta_1|d_2a_1}\sum_{\delta_2|d_1a_1}\mathop{{\sum}^*}_{l=1}^{a_1}s_1\left(\frac{l}{a_1}\right)\right)=\\=
O\left(x_1x_2N^2\sum\limits_{d_1d_2\le x_3N \atop (d_1,d_2)=1}\sum_{a_1\le \frac{x_3N}{d_1d_2}}
\tau(d_1a_1)\tau(d_2a_1)a_1\log^2 a_1\right)=\\=O\left(x_1x_2N^2\sum\limits_{d_1d_2\le x_3N \atop (d_1,d_2)=1}
(d_1d_2)^{\varepsilon_1}\sum_{a_1\le \frac{x_3N}{d_1d_2}}a_1^{1+\varepsilon}\right)=
O\left(x_1x_2x_3^{2+\varepsilon}N^{4+\varepsilon}\right).
\end{gather*}
Подставляя полученную оценку в \eqref{FF} и используя формулы \eqref{1}---\eqref{4}, получаем
\begin{gather}\nonumber
F=\sum\limits_{d_1d_2\le x_3N \atop (d_1,d_2)=1}d_1d_2\sum_{a_1\le \frac{x_3N}{d_1d_2}}
\sum_{\delta_1|d_2a_1}\mu(\delta_1)\sum_{\delta_2|d_1a_1}\mu(\delta_2)\frac{(a_1,\delta_1,\delta_2)}{a_1\delta_1\delta_2}\\
\int_0^{\frac{x_1N}{d_1}}\int_0^{\frac{x_2N}{d_2}}\left(t_1\lambda^{*}\left(a_1,\frac{t_2}{t_1}\right)+
t_2\lambda^{*}\left(a_1,\frac{t_1}{t_2}\right)\right)dt_1dt_2+
O\left(x_1x_2x_3^{2+\varepsilon}N^{4+\varepsilon}\right).
\end{gather}
Меняем порядок суммирования
\begin{gather}\nonumber
F=\sum\limits_{d_1d_2\le x_3N \atop (d_1,d_2)=1}d_1d_2\sum_{\delta_1\le \frac{x_3N}{d_1}}\frac{\mu(\delta_1)}{\delta_1}
\sum_{\delta_2\le \frac{x_3N}{d_2}}\frac{\mu(\delta_2)}{\delta_2}
\sum\limits_{a_1\le \frac{x_3N}{d_1d_2} \atop \delta_1|d_2a_1, \delta_2|d_1a_1}(a_1,\delta_1,\delta_2)\\
\int_0^{\frac{x_1N}{d_1}}\int_0^{\frac{x_2N}{d_2}}\frac{t_1\lambda^{*}\left(a_1,\frac{t_2}{t_1}\right)+
t_2\lambda^{*}\left(a_1,\frac{t_1}{t_2}\right)}{a_1}dt_1dt_2+
O\left(x_1x_2x_3^{2+\varepsilon}N^{4+\varepsilon}\right).
\end{gather}
Заметим, что $(\delta_1,\delta_2)|(d_1a_1,d_2a_1)$, но $(d_1a_1,d_2a_1)=a_1(d_1,d_2)=a_1$. Следовательно, $(\delta_1,\delta_2)|a_1$, тогда
$$(a_1,\delta_1,\delta_2)=(a_1,(\delta_1,\delta_2))=(\delta_1,\delta_2).$$
Условия $\delta_1|d_2a_1, \delta_2|d_1a_1$ равносильны условию $\delta|a_1$, где
$$\delta=HOK\left(\frac{\delta_1}{(\delta_1,d_2)},\frac{\delta_2}{(\delta_2,d_1)}\right).$$
Тем самым лемма \ref{F} доказана.
\end{proof}
\end{Le}
\begin{Zam}
\label{G}
Справедливо следующее равенство
\begin{gather}\nonumber
G=
\sum\limits_{d_1d_2\le x_3N \atop (d_1,d_2)=1}d_1d_2\sum_{\delta_1\le \frac{x_3N}{d_1}}\sum_{\delta_2\le \frac{x_3N}{d_2}}
\frac{\mu(\delta_1)}{\delta_1}\frac{\mu(\delta_2)}{\delta_2}(\delta_1,\delta_2)\\\nonumber
\int_0^{\frac{x_1N}{d_1}}\int_0^{\frac{x_2N}{d_2}}\sum\limits_{a\le \frac{x_3N}{d_1d_2} \atop \delta|a}
\frac{8}{\pi}\frac{\varphi(a)}{\sqrt{a}}\sqrt{t_1t_2}dt_1dt_2+
O\left(x_1x_2x_3^{2+\varepsilon}N^{4+\varepsilon}\right).
\end{gather}
\begin{proof}
Доказывается аналогично лемме \ref{F}.
\end{proof}
\end{Zam}
Из леммы \ref{F} следует, что для доказательства теоремы \ref{Th1} необходимо исследовать сумму вида
$$\sum\limits_{a\le R \atop \delta|a}\frac{\lambda^{*}(a,\alpha)}{a}.$$
Используя формулу \eqref{ljmbda*} получаем
\begin{gather}
\sum\limits_{a\le R \atop \delta|a}\frac{\lambda^{*}(a,\alpha)}{a}=
\sum_{d_1d_2\le R}\sum\limits_{a\le R \atop \delta|a, d_1d_2|a}\mu(d_1)
\mu(d_2)\frac{d_2}{a}\lambda\left(\frac{a}{d_1d_2},\frac{d_1\alpha}{d_2}\right)=\notag\\=
\sum_{n|\delta}\sum\limits_{d_1d_2\le R \atop (d_1d_2,\delta)=n}
\sum\limits_{a\le R \atop  \frac{\delta}{n}d_1d_2|a}\mu(d_1)\mu(d_2)\frac{d_2}{a}
\lambda\left(\frac{a}{d_1d_2},\frac{d_1\alpha}{d_2}\right)=\notag\\=
\sum_{n|\delta}\sum\limits_{d_1d_2\le R \atop (d_1d_2,\delta)=n}\frac{\mu(d_1)\mu(d_2)}{d_1}
\sum\limits_{a\le \frac{R}{d_1d_2} \atop  \frac{\delta}{n}|a}\frac{\lambda\left(a,\frac{d_1}{d_2}\alpha\right)}{a}.
\label{lj-lj}
\end{gather}
Используя лемму \ref{abel}, получаем
\begin{equation}\label{ljabel}
\sum\limits_{a\le R \atop \delta|a}\frac{\lambda(a,\alpha)}{a}=\frac{1}{R}
\sum\limits_{a\le R \atop \delta|a}\lambda(a,\alpha)+\int_{1}^R\frac{1}{t^2}
\sum\limits_{a\le t \atop \delta|a}\lambda(a,\alpha)dt.
\end{equation}
Следовательно, задача свелась к нахождению асимптотической формулы для суммы
$$\sum\limits_{a\le R \atop \delta|a}\lambda(a,\alpha).$$
Нахождению этой формулы и посвящены последующие разделы.
\section{Разделение задачи на отдельные случаи}
Используя формулу \eqref{deltasum}, из формулы (\ref{ljmbda}) легко получить следующее равенство
\begin{gather}\nonumber
\delta\sum\limits_{a\le R \atop \delta|a}\lambda(a,\alpha)=
\sum_{z=1}^{\delta}\sum_{n\ge 1}\sum_{k=1}^n\sum_{Q’\ge 1}\sum_{Q\ge 0}\\\nonumber
\left[nQ’+kQ\le R, \alpha(n-k)<Q\le \alpha n\right]\exp\left(2\pi i\frac{nQ’+kQ}{\delta}z\right)(Q’+Q+\alpha k)=\\=
\nonumber
\sum_{z=1}^{\delta}\sum_{n\ge 1}\sum_{k=1}^n\sum_{Q’\ge 1}\sum_{Q\ge 0}
\left[nQ’+kQ\le R, \alpha(n-k)<Q\le \alpha n\right]\\
\exp\left(2\pi i\frac{nQ’}{\delta}z\right)\exp\left(2\pi i\frac{kQ}{\delta}z\right)(Q’+Q+\alpha k).
\end{gather}

Разобьем область суммирования по переменным $(n,k,\Q,Q)$
\begin{gather}
\left\{
  \begin{array}{ll}
    nQ’+kQ\le R, &  \\
    \alpha(n-k)<Q\le \alpha n, & \\
    1\le k\le n, &  \\
    0\le Q,1\le \Q.
    \end{array}
\right.\label{sistem}
\end{gather}
на пять подобластей. Тогда
\begin{equation}\label{summa}
\delta\sum\limits_{a\le R \atop \delta|a}\lambda(a,\alpha)=\Sigma_1+\Sigma_2+\Sigma_3+\Sigma_4+\Sigma_5,
\end{equation}
где $\Sigma_i$~--- сумма для i-го случая.
Для этого введем параметры
\begin{equation}
\label{U1U2R}
U_1=\sqrt{\frac{R}{\alpha}}, \qquad  U_2=\sqrt{R\alpha}.
\end{equation}
Очевидно, что выполнено
\begin{equation}
\label{U1U2}
U_1U_2=R \quad \mbox{и} \quad U_2=\alpha U_1.
\end{equation}
В Случаи 1 мы рассматриваем ситуацию, когда  $n\le U_1$. Следовательно, остальные случаи посвящены рассмотрению ситуации  $n> U_1$, но тогда $\Q\le U_2$. В Случаях 2 и 3 рассматривается ситуация, когда $k\le U_1$, а в Случаях 4 и 5 рассматривается ситуация, когда $k>U_1$.
\subsection{Случай 1}
Пусть $n\le U_1$ и внешнее суммирование ведется по области
$$\Omega_1=\left\{n\le U_1,  k\le n\right\}.$$
Тогда переменные $Q,\Q$ должны удовлетворять условиям
\par
\begin{picture}(250,160)
\put(0,0){\vector(0,1){150}}
\put(0,0){\vector(1,0){250}}
\put(250,0){Q}
\put(0,140){$\Q$}
\put(200,10){$\frac{R}{k}$}
\put(01,110){$\frac{R}{n}$}
\put(0,100){\line(2,-1){200}}
\put(40,0){\line(0,1){80}}
\put(120,0){\line(0,1){40}}
\put(30,-10){$\alpha (n-k)$}
\put(115,-10){$\alpha n$}
\put(80,70){$\frac{R-kQ}{n}$}
\end{picture}
\\[\bigskipamount]
$$
\alpha (n-k)<Q\le \min\left(\alpha n,\frac{R}{k}\right), \qquad \Q\le \frac{R-kQ}{n}.
$$
В силу выбора параметров получаем внутреннее суммирование ведется по области
$$
\Omega_{11}=\left\{\alpha (n-k)<Q\le \alpha n, \qquad 1\le\Q\le \frac{R-kQ}{n}\right\}.
$$
Следовательно,
\begin{gather*}\nonumber
\Sigma_1=\sum_{z=1}^{\delta}\sum_{(n,k)\in\Omega_1}\sum_{(Q,\Q)\in\Omega_{11}}
\exp\left(2\pi i\frac{zn}{\delta}\Q\right)\exp\left(2\pi i\frac{zk}{\delta}Q\right)(\Q+Q+\alpha k).
\end{gather*}
Следующее преобразование очевидно
\begin{gather*}
\sum_{z=1}^{\delta}\sum_{n\le U_1}\sum_{k\le n}=
\sum_{z=1}^{\delta}\sum\limits_{n\le U_1\atop \delta|zn}\sum\limits_{k\le n\atop \delta|zk}+
\sum_{z=1}^{\delta}\sum\limits_{n\le U_1\atop \delta|zn}\sum\limits_{k\le n\atop \delta\nmid zk}+
\sum_{z=1}^{\delta}\sum\limits_{n\le U_1\atop \delta\nmid zn}\sum\limits_{k\le n}.
\end{gather*}
Преобразуем первую сумму, учитывая, что мы суммируем функцию в которую $z$ и $\delta$ входят только в виде $\frac{z}{\delta}$.
\begin{gather*}
\sum_{z=1}^{\delta}\sum\limits_{n\le U_1\atop \delta|zn}\sum\limits_{k\le n\atop \delta|zk}=
\sum_{d|\delta}\sum\limits_{z=1\atop (\delta,z)=d}^{\delta}
\sum\limits_{n\le U_1\atop \frac{\delta}{d}|n}\sum\limits_{k\le n\atop \frac{\delta}{d}|k}=
\sum_{\tau|\delta}\mathop{{\sum}^*}_{z=1}^{\tau}
\sum\limits_{n\le U_1\atop \tau|n}\sum\limits_{k\le n\atop \tau|k}.
\end{gather*}
Аналогично для второй суммы получаем
\begin{gather*}
\sum_{z=1}^{\delta}\sum\limits_{n\le U_1\atop \delta|zn}\sum\limits_{k\le n\atop \delta\nmid zk}=
\sum_{\tau|\delta}\mathop{{\sum}^*}_{z=1}^{\tau}
\sum\limits_{n\le U_1\atop \tau|n}\sum\limits_{k\le n\atop \tau\nmid k}.
\end{gather*}
Следовательно,
\begin{gather}
\Sigma_1=\sum_{\tau|\delta}\mathop{{\sum}^*}_{z=1}^{\tau}\sum_{\substack{(n,k)\in\Omega_1\\ \tau|n, \tau|k}}
\sum_{(Q,\Q)\in\Omega_{11}}(\Q+Q+\alpha k)+\notag\\+
\sum_{\tau|\delta}\mathop{{\sum}^*}_{z=1}^{\tau}\sum_{\substack{(n,k)\in\Omega_1\\ \tau|n, \tau\nmid k}}
\sum_{(Q,\Q)\in\Omega_{11}}\exp\left(2\pi i\frac{zk}{\tau}Q\right)(\Q+Q+\alpha k)+\notag\\+
\sum_{z=1}^{\delta}\sum_{\substack{(n,k)\in\Omega_1\\\delta\nmid zn}}
\sum_{(Q,\Q)\in\Omega_{11}}\exp\left(2\pi i\frac{zk}{\delta}Q\right)
\exp\left(2\pi i\frac{zn}{\delta}\Q\right)(\Q+Q+\alpha k)=\notag\\
=\Sigma_{11}+\Sigma_{12}+\Sigma_{13}.\label{slu1}
\end{gather}

\subsection{Случай 2}
Пусть $n> U_1$, тогда $\Q\le U_2$ и пусть $k\le U_1$. Пусть также
$$U_1\le \frac{R}{\Q+\alpha k},$$
тогда внешнее суммирование ведется по области
$$
\Omega_{2}=\left\{k\le U_1, \quad 1\le\Q\le U_2-\alpha k\right\}.
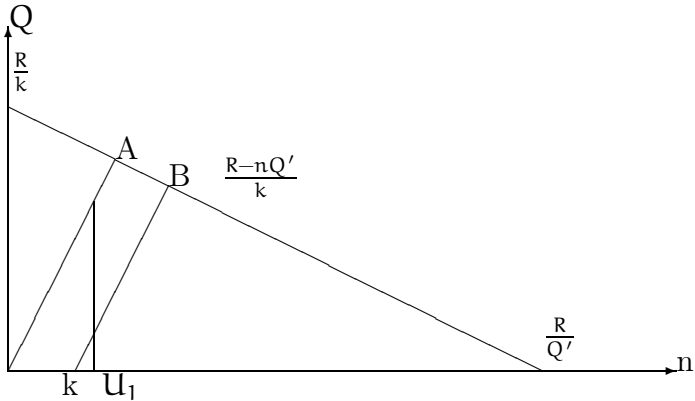
$$
\\[\bigskipamount]
\begin{picture}(250,160)
\put(0,15){\begin{picture}(250,150)
\put(0,0){\vector(0,1){130}}
\put(0,0){\vector(1,0){250}}
\put(250,0){n}
\put(0,130){$Q$}
\put(200,10){$\frac{R}{\Q}$}
\put(01,110){$\frac{R}{k}$}
\put(0,100){\line(2,-1){200}}
\put(0,0){\line(1,2){40}}
\put(25,0){\line(1,2){35}}
\put(20,-9){k}
\put(40,80){A}
\put(60,70){B}
\put(32,0){\line(0,1){64}}
\put(35,-9){$U_1$}
\put(80,70){$\frac{R-n\Q}{k}$}
\end{picture}}
\end{picture}
\\[\bigskipamount]
Пусть в плоскости $(n, Q)$ A---точка пересечения прямых $Q=\alpha n$ и $Q=\frac{R-n\Q}{k}$, B---точка пересечения прямых $Q=\alpha (n-k)$ и $Q=\frac{R-n\Q}{k}$.Тогда их координаты
$$
A\left(\frac{R}{\Q+\alpha k},\frac{R\alpha}{\Q+\alpha k}\right),\qquad
B\left(\frac{R+\alpha k^2}{\Q+\alpha k},\alpha\frac{R-\alpha\Q}{\Q+\alpha k}\right).
$$
Внутреннее суммирование ведется по объединению непересекающихся областей
$$
\Omega_{21}=\left\{U_1<n\le\frac{R}{\Q+\alpha k},\quad \alpha (n-k)<Q\le\alpha n \right\}
$$
и
$$
\Omega_{22}=\left\{\frac{R}{\Q+\alpha k}<n\le\frac{R+\alpha k^2}{\Q+\alpha k},\quad \alpha (n-k)<Q\le\frac{R-n\Q}{k}\right\}.
$$
Следовательно,
\begin{gather*}
\Sigma_2=\sum_{z=1}^{\delta}\sum_{(k,\Q)\in\Omega_{2}}\left(
\sum\limits_{(n,Q)\in\Omega_{21}}+\sum\limits_{(n,Q)\in\Omega_{22}}\right)
\exp\left(2\pi i\frac{zn}{\delta}\Q\right)\exp\left(2\pi i\frac{zk}{\delta}Q\right)(\Q+Q+\alpha k).
\end{gather*}
Проделывая преобразования аналогичные тем, что были сделаны для $\Sigma_1$ получаем
\begin{gather}
\label{slu2}
\Sigma_2=\sum_{\tau|\delta}\mathop{{\sum}^*}_{z=1}^{\tau}\sum_{\substack{(k,\Q)\in\Omega_{2}\\\tau\mid k,\tau\mid \Q}}
\left(\sum\limits_{(n,Q)\in\Omega_{21}}+\sum\limits_{(n,Q)\in\Omega_{22}}\right)(\Q+Q+\alpha k)+\notag\\r+
\sum_{\tau|\delta}\mathop{{\sum}^*}_{z=1}^{\tau}\sum_{\substack{(k,\Q)\in\Omega_{2}\\\tau\mid k,\tau\nmid \Q}}
\left(\sum\limits_{(n,Q)\in\Omega_{21}}+\sum\limits_{(n,Q)\in\Omega_{22}}\right)
\exp\left(2\pi i\frac{z\Q}{\tau}n\right)(\Q+Q+\alpha k)+\notag\\+
\sum_{z=1}^{\delta}\sum_{\substack{(k,\Q)\in\Omega_{2}\\\delta\nmid zk,}}
\left(\sum\limits_{(n,Q)\in\Omega_{21}}+\sum\limits_{(n,Q)\in\Omega_{22}}\right)
\exp\left(2\pi i\frac{z\Q}{\delta}n\right)\exp\left(2\pi i\frac{zk}{\delta}Q\right)(\Q+Q+\alpha k)=\notag\\
=\Sigma_{21}+\Sigma_{22}+\Sigma_{23}.\label{slu2}
\end{gather}

\subsection{Случай 3}
Пусть $n> U_1$,тогда $\Q\le U_2$ и пусть $k\le U_1$, а внешнее суммирование ведется по $k, \Q$ и выполнено
$$
\frac{R}{\Q+\alpha k}<U_1\le \frac{R+\alpha k^2}{\Q+\alpha k}.
$$
Тогда внешнее суммирование ведется по области
$$
\Omega_{3}=\left\{k\le U_1, \quad U_2-\alpha k<\Q\le U_2-\alpha k+\frac{\alpha k^2}{U_1}\right\},
$$
а внутреннее по области
$$
\Omega_{31}=\left\{U_1<n\le\frac{R+\alpha k^2}{\Q+\alpha k},\quad \alpha (n-k)<Q\le\frac{R-n\Q}{k} \right\}.
$$
Следовательно,
\begin{gather*}
\Sigma_3=\sum_{z=1}^{\delta}\sum_{(k,\Q)\in\Omega_{3}}\sum\limits_{(n,Q)\in\Omega_{31}}
\exp\left(2\pi i\frac{zn}{\delta}\Q\right)\exp\left(2\pi i\frac{zk}{\delta}Q\right)(\Q+Q+\alpha k).
\end{gather*}
Проделывая преобразования аналогичные тем, что были сделаны для $\Sigma_1$, получаем
\begin{gather}
\Sigma_3=
\sum_{\tau|\delta}\mathop{{\sum}^*}_{z=1}^{\tau}\sum_{\substack{(k,\Q)\in\Omega_{3}\\\tau\mid k,\tau\mid \Q}}
\sum\limits_{(n,Q)\in\Omega_{31}}(\Q+Q+\alpha k)+\notag\\+
\sum_{\tau|\delta}\mathop{{\sum}^*}_{z=1}^{\tau}\sum_{\substack{(k,\Q)\in\Omega_{3}\\
\tau\mid k,\tau\nmid \Q}}\sum\limits_{(n,Q)\in\Omega_{31}}\exp\left(2\pi i\frac{z\Q}{\tau}n\right)(\Q+Q+\alpha k)+\notag\\+
\sum_{z=1}^{\delta}\sum_{\substack{(k,\Q)\in\Omega_{3}\\\delta\nmid zk,}}
\sum\limits_{(n,Q)\in\Omega_{31}}
\exp\left(2\pi i\frac{z\Q}{\delta}n\right)\exp\left(2\pi i\frac{zk}{\delta}Q\right)(\Q+Q+\alpha k)=\notag\\
=\Sigma_{31}+\Sigma_{32}+\Sigma_{33}.\label{slu3}
\end{gather}

\subsection{Случай 4}
Пусть $n> U_1$,тогда $\Q\le U_2$. Пусть $k>U_1$, тогда $Q\le U_2$. Будем вести внешнее суммирование по $Q,\Q$.
\\[\bigskipamount]
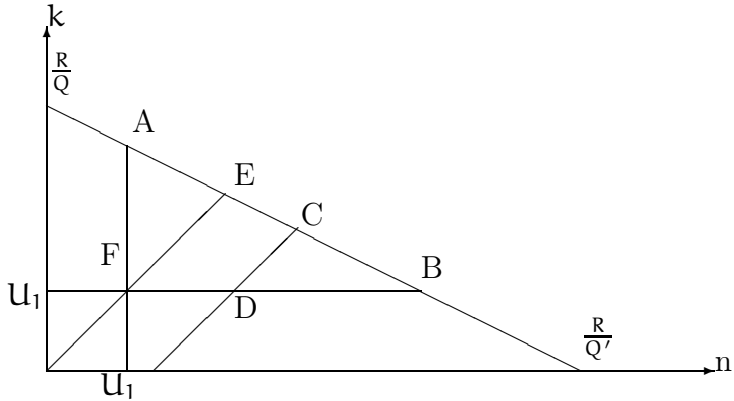
\begin{picture}(250,160)
\put(20,20){\begin{picture}(250,150)
\put(0,0){\vector(0,1){130}}
\put(0,0){\vector(1,0){250}}
\put(250,0){n}
\put(0,130){$k$}
\put(200,10){$\frac{R}{\Q}$}
\put(01,110){$\frac{R}{Q}$}
\put(0,100){\line(2,-1){200}}
\put(0,0){\line(1,1){67}}
\put(40,0){\line(1,1){54}}
\put(20,-9){$U_1$}
\put(-15,25){$U_1$}
\put(32,90){A}
\put(140,35){B}
\put(95,55){C}
\put(70,20){D}
\put(70,70){E}
\put(20,40){F}
\put(0,30){\line(1,0){140}}
\put(30,0){\line(0,1){85}}
\end{picture}}
\end{picture}
\\[\bigskipamount]
Пусть в плоскости $(n, k)$ A,B,E,C---точки пересечения прямой $k=\frac{R-n\Q}{Q}$ соответственно с прямыми $n=U_1,$ \mbox{$k=U_1,$} $k=n,$ $k=n-\frac{Q}{\alpha}$.
D---точка пересечения прямой $k=n-\frac{Q}{\alpha}$ с прямой $k=U_1$. F---точка пересечения прямой $k=U_1$ с прямой $n=U_1$.
Точки имеют следующие координаты
\begin{gather*}
A\left(U_1,\frac{R-U_1\Q}{Q}\right),\quad B\left(\frac{R-U_1Q}{\Q},U_1\right),\quad
C\left(\frac{R+Q^2/\alpha}{Q+\Q},\frac{R-Q\Q/\alpha}{Q+\Q}\right),\quad
D\left(U_1+\frac{Q}{\alpha},U_1\right),\\
F\left(U_1,U_1\right),\quad E\left(\frac{R}{Q+\Q},\frac{R}{Q+\Q}\right).
\end{gather*}
Чтобы $\triangle AFB$ был невырожденным необходимо
$$Q+\Q\le U_2.$$
Рассмотрим случай, когда C принадлежит отрезку AB, следовательно
$$
\frac{R+Q^2/\alpha}{Q+\Q}\le\frac{R-U_1Q}{\Q}.
$$
Получаем
$$
\Q\le\frac{R-U_1Q}{U_1+Q/\alpha}=U_2\frac{U_2-Q}{U_2+Q}
$$
и условие $Q+\Q\le U_2$ выполнено.
Тогда внешнее суммирование ведется по области
$$
\Omega_{4}=\left\{Q\le U_2,\quad\Q\le U_2\frac{U_2-Q}{U_2+Q}\right\},
$$
а внутреннее по области FECD, которую мы разбиваем следующим образом
$$FECD=\Omega_{41}+\Omega_{42}-\Omega_{43},$$
где
$$
\Omega_{41}=\left\{U_1<n\le\frac{R}{Q+\Q},\quad U_1<k\le n \right\},
$$
$$
\Omega_{42}=\left\{\frac{R}{Q+\Q}<n\le\frac{R+Q^2/\alpha}{Q+\Q},\quad U_1<k\le\frac{R-n\Q}{Q}\right\},
$$
$$
\Omega_{43}=\left\{U_1+\frac{Q}{\alpha}<n\le\frac{R+Q^2/\alpha}{Q+\Q},\quad U_1<k\le n-\frac{Q}{\alpha}\right\}.
$$
Следовательно,
\begin{gather*}
\Sigma_4=\sum_{z=1}^{\delta}\sum_{(Q,\Q)\in\Omega_{4}}
\left(\sum_{(n,k)\in\Omega_{41}}+\sum_{(n,k)\in\Omega_{42}}-\sum_{(n,k)\in\Omega_{43}}\right)
\exp\left(2\pi i\frac{z\Q}{\delta}n\right)\exp\left(2\pi i\frac{zQ}{\delta}k\right)(\Q+Q+\alpha k).
\end{gather*}
Проделывая преобразования аналогичные тем, что были сделаны для $\Sigma_1$ получаем
\begin{gather}
\label{slu4}
\Sigma_4=\sum_{\tau|\delta}\mathop{{\sum}^*}_{z=1}^{\tau}\sum_{\substack{(Q,\Q)\in\Omega_{4}\\\tau\mid Q,\tau\mid \Q}}
\left(\sum_{(n,k)\in\Omega_{41}}+\sum_{(n,k)\in\Omega_{42}}-\sum_{(n,k)\in\Omega_{43}}\right)(\Q+Q+\alpha k)+\notag\\+
\sum_{\tau|\delta}\mathop{{\sum}^*}_{z=1}^{\tau}\sum_{\substack{(Q,\Q)\in\Omega_{4}\\\tau\mid Q,\tau\nmid \Q}}
\left(\sum_{(n,k)\in\Omega_{41}}+\sum_{(n,k)\in\Omega_{42}}-\sum_{(n,k)\in\Omega_{43}}\right)
\exp\left(2\pi i\frac{z\Q}{\tau}n\right)(\Q+Q+\alpha k)+\notag\\+
\sum_{z=1}^{\delta}\sum_{\substack{(Q,\Q)\in\Omega_{4}\\\delta\nmid zk,}}
\left(\sum_{(n,k)\in\Omega_{41}}+\sum_{(n,k)\in\Omega_{42}}-\sum_{(n,k)\in\Omega_{43}}\right)
\exp\left(2\pi i\frac{z\Q}{\delta}n\right)\exp\left(2\pi i\frac{zQ}{\delta}k\right)(\Q+Q+\alpha k)=\notag\\
=\Sigma_{41}+\Sigma_{42}+\Sigma_{43}.\label{slu4}
\end{gather}

\subsection{Случай 5}
Пусть $n> U_1$,тогда $\Q\le U_2$. Пусть $k>U_1$, тогда $Q\le U_2$. Будем вести внешнее суммирование по $Q,\Q$.
Рассмотрим случай, когда C  не принадлежит отрезку AB, следовательно
$$
\frac{R+Q^2/\alpha}{Q+\Q}>\frac{R-U_1Q}{\Q}.
$$
Получаем
$$
\Q>\frac{R-U_1Q}{U_1+Q/\alpha}=U_2\frac{U_2-Q}{U_2+Q}
$$
и условие $Q+\Q\le U_2$ необходимо учитывать.
Тогда внешнее суммирование ведется по области
$$
\Omega_{5}=\left\{Q\le U_2,\quad U_2\frac{U_2-Q}{U_2+Q}<Q\le U_2-Q\right\},
$$
а внутреннее по $\triangle FBE$, который мы разбиваем следующим образом
$$FBE=\Omega_{51}+\Omega_{52},$$
где
$$
\Omega_{51}=\left\{U_1<n\le\frac{R}{Q+\Q},\quad U_1<k\le n \right\},
$$
$$
\Omega_{52}=\left\{\frac{R}{Q+\Q}<n\le\frac{R-U_1Q}{\Q},\quad U_1<k\le\frac{R-n\Q}{Q}\right\}.
$$
Следовательно,
\begin{gather*}
\Sigma_5=\sum_{z=1}^{\delta}\sum_{(Q,\Q)\in\Omega_{5}}
\left(\sum_{(n,k)\in\Omega_{51}}+\sum_{(n,k)\in\Omega_{52}}\right)
\exp\left(2\pi i\frac{z\Q}{\delta}n\right)\exp\left(2\pi i\frac{zQ}{\delta}k\right)(\Q+Q+\alpha k).
\end{gather*}
Проделывая преобразования аналогичные тем, что были сделаны для $\Sigma_1$ получаем
\begin{gather}
\Sigma_5=\sum_{\tau|\delta}\mathop{{\sum}^*}_{z=1}^{\tau}\sum_{\substack{(Q,\Q)\in\Omega_{5}\\\tau\mid Q,\tau\mid \Q}}
\left(\sum_{(n,k)\in\Omega_{51}}+\sum_{(n,k)\in\Omega_{52}}\right)(\Q+Q+\alpha k)+\notag\\+
\sum_{\tau|\delta}\mathop{{\sum}^*}_{z=1}^{\tau}\sum_{\substack{(Q,\Q)\in\Omega_{5}\\\tau\mid Q,\tau\nmid \Q}}
\left(\sum_{(n,k)\in\Omega_{51}}+\sum_{(n,k)\in\Omega_{52}}\right)
\exp\left(2\pi i\frac{z\Q}{\tau}n\right)(\Q+Q+\alpha k)+\notag\\+
\sum_{z=1}^{\delta}\sum_{\substack{(Q,\Q)\in\Omega_{4}\\\delta\nmid zk,}}
\left(\sum_{(n,k)\in\Omega_{51}}+\sum_{(n,k)\in\Omega_{52}}\right)
\exp\left(2\pi i\frac{z\Q}{\delta}n\right)\exp\left(2\pi i\frac{zQ}{\delta}k\right)(\Q+Q+\alpha k)=\notag\\
=\Sigma_{51}+\Sigma_{52}+\Sigma_{53}.\label{slu5}
\end{gather}
\section{Вычисление сумм первого типа}
В этом параграфе мы вычислим $\Sigma_{11}, \Sigma_{21}, \Sigma_{31}, \Sigma_{41}, \Sigma_{51}.$
\subsection{Случай 1}
\begin{Le}
\label{11}
Справедлива следующая асимптотическая формула
\begin{gather*}
\Sigma_{11}=
\frac{53}{150}R^{5/2}\sqrt{\alpha}\sum_{\tau|\delta}\frac{\varphi(\tau)}{\tau^2}+\sum_{\tau|\delta}\varphi(\tau)
O\left(\frac{\alpha R^2}{\tau}\left(1+\log \frac{U_1}{\tau}\right)+\frac{R^2}{\tau^2}\log \frac{U_1}{\tau}\right),
\end{gather*}
где $\Sigma_{11}$ определена в  \eqref{slu1}.
\begin{proof}
Производя суммирование по переменной $\Q$, получаем
\begin{gather*}
\sum_{\Q\le\frac{R-kQ}{n}}(\Q+Q+\alpha k)=(Q+\alpha k)\frac{R-kQ}{n}+\frac{(R-kQ)^2}{2n^2}+
O\left(Q+\alpha k\right)+O\left(\frac{R-kQ}{n}\right)=\\=
Q^2\left(\frac{k^2}{2n^2}-\frac{k}{n}\right)+Q\left(\frac{R}{n}-\frac{\alpha k^2}{n}-R\frac{k}{n^2}\right)+
\left(\alpha R\frac{k}{n}+\frac{R^2}{2n^2}\right)+O\left(Q+\alpha k\right)+O\left(\frac{R-kQ}{n}\right).
\end{gather*}
Далее необходимо просуммировать полученное выражение по переменной $Q$. Учитывая, что
$$(n,k)\in\Omega_1=\left\{n\le U_1,  k\le n\right\}\quad \mbox{и}\quad U_1=\sqrt{\frac{R}{\alpha}},$$
получаем следующие соотношения, необходимые для получения оценки остаточного члена
\begin{gather*}
\left|\frac{k^2}{2n^2}-\frac{k}{n}\right|\alpha^2 n^2\ll\alpha^2 kn,\quad
\left|\frac{R}{n}-\frac{\alpha k^2}{n}-R\frac{k}{n^2}\right|\alpha n\ll\alpha R,\quad
\left|\alpha R\frac{k}{n}+\frac{R^2}{2n^2}\right|\ll\frac{R^2}{n^2},\\
\sum_{\alpha (n-k)<Q\le \alpha n} \left(Q+\alpha k\right)\ll\alpha^2 nk,\qquad
\sum_{\alpha (n-k)<Q\le \alpha n} \frac{R-kQ}{n}\ll\frac{R}{n}\alpha k,\\
\alpha^2 kn+\alpha R+\frac{R}{n}\alpha k\ll\frac{R^2}{n^2},
\end{gather*}
Напомним, что
$$
\Omega_{11}=\left\{\alpha (n-k)<Q\le \alpha n, \qquad 1\le\Q\le \frac{R-kQ}{n}\right\}.
$$
Применяя следствие \ref{sl} (a) (см. Приложение), получаем
\begin{gather*}
\sum_{(Q,\Q)\in\Omega_{11}}(\Q+Q+\alpha k)=\frac{\alpha^3}{3}\left(\frac{k^2}{2n^2}-\frac{k}{n}\right)(n^3-(n-k)^3)+\\+
\frac{\alpha^2}{2}\left(\frac{R}{n}-\frac{\alpha k^2}{n}-\frac{Rk}{n^2}\right)(n^2-(n-k)^2)+
\alpha k\left(\alpha R\frac{k}{n}+\frac{R^2}{2n^2}\right)+O\left(\frac{R^2}{n^2}\right)=\\=
\alpha^3\left(\frac{k^3}{2}-\frac{k^4}{3n}+\frac{k^5}{6n^2}-nk^2\right)+
\alpha^2\left(Rk-R\frac{k^2}{2n}+R\frac{k^3}{2n^2}\right)+
\alpha R^2\frac{k}{2n^2}+O\left(\frac{R^2}{n^2}\right).
\end{gather*}
Полученное выражение необходимо просуммировать по переменным $$(n,k)\in\Omega_1, \tau|n, \tau|k,$$ получаем
\begin{gather*}
\sum_{\substack{(n,k)\in\Omega_1\\ \tau|n, \tau|k}}
\sum_{(Q,\Q)\in\Omega_{11}}(\Q+Q+\alpha k)=\\=
\sum_{n\le \frac{U_1}{\tau}}\sum_{k\le n}\alpha^3\tau^3\left(\frac{k^3}{2}-\frac{k^4}{3n}+\frac{k^5}{6n^2}-nk^2\right)+
\alpha^2\tau\left(Rk-R\frac{k^2}{2n}+R\frac{k^3}{2n^2}\right)+
\alpha R^2\frac{k}{2n^2\tau}+O\left(\frac{R^2}{\tau^2n^2}\right).
\end{gather*}
Применяя следствие \ref{sl} (a) (см. Приложение) сначала при суммировании по $k$, а потом при суммировании по $n$, получаем
\begin{gather*}
\sum_{\substack{(n,k)\in\Omega_1\\ \tau|n, \tau|k}}
\sum_{(Q,\Q)\in\Omega_{11}}(\Q+Q+\alpha k)=\\=
\sum_{n\le \frac{U_1}{\tau}}\Biggl(
\alpha^3\tau^3\left(\frac{n^4}{8}-\frac{n^4}{15}+\frac{n^4}{36}-\frac{n^4}{3}+O(n^3)\right)+
\alpha^2\tau\left(R\frac{n^2}{2}-R\frac{n^2}{6}+R\frac{n^2}{8}+O(Rn)\right)+\\+
\alpha R^2\frac{1}{4\tau}+O\left(\frac{\alpha R^2}{\tau n}\right)+O\left(\frac{R^2}{\tau^2n}\right)
\Biggl)=\\=
\sum_{n\le \frac{U_1}{\tau}}\Biggl(
-\frac{89}{360}\alpha^3\tau^3n^4+\frac{11}{24}\alpha^2\tau Rn^2+\frac{\alpha R^2}{4\tau}+
O\left(\alpha^3\tau^3n^3+\alpha^2\tau Rn+\frac{\alpha R^2}{\tau n}\right)+O\left(\frac{R^2}{\tau^2n}\right)\Biggl)=\\=
\frac{R^{5/2}\sqrt{\alpha}}{\tau^2}\left(-\frac{89}{5\cdot360}+\frac{11}{3\cdot24}+\frac{1}{4}\right)+
O\left(\frac{\alpha R^2}{\tau}\left(1+\log \frac{U_1}{\tau}\right)+\frac{R^2}{\tau^2}\log \frac{U_1}{\tau}\right)=\\=
\frac{R^{5/2}\sqrt{\alpha}}{\tau^2}\frac{53}{150}+
O\left(\frac{\alpha R^2}{\tau}\left(1+\log \frac{U_1}{\tau}\right)+\frac{R^2}{\tau^2}\log \frac{U_1}{\tau}\right).
\end{gather*}
Используя определение функции Эйлера
$$\varphi(\tau)=\mathop{{\sum}^*}_{z=1}^{\tau}1,$$
получаем асимптотическую формулу для  $\Sigma_{11}$
\begin{gather*}
\Sigma_{11}=\frac{53}{150}R^{5/2}\sqrt{\alpha}\sum_{\tau|\delta}\frac{\varphi(\tau)}{\tau^2}+\sum_{\tau|\delta}\varphi(\tau)
O\left(\frac{\alpha R^2}{\tau}\left(1+\log \frac{U_1}{\tau}\right)+\frac{R^2}{\tau^2}\log \frac{U_1}{\tau}\right).
\end{gather*}
Тем самым лемма доказана.
\end{proof}
\end{Le}

\subsection{Случай 2}
\begin{Le}
\label{21}
Справедлива следующая асимптотическая формула
\begin{gather*}
\Sigma_{21}=\frac{19}{75}R^{5/2}\sqrt{\alpha}\sum_{\tau|\delta}\frac{\varphi(\tau)}{\tau^2}+\sum_{\tau|\delta}\varphi(\tau)
O\left(\frac{R^2}{\tau}\left(1+\alpha\right)+\frac{R^2}{\tau}(1+\alpha)\log \frac{U_1}{\tau}\right),
\end{gather*}
где $\Sigma_{21}$ определена в \eqref{slu2}.
\begin{proof}
Вычислим отдельно суммы по областям
$$
\Omega_{21}=\left\{U_1<n\le\frac{R}{\Q+\alpha k},\quad \alpha (n-k)<Q\le\alpha n \right\}
$$
и
$$
\Omega_{22}=\left\{\frac{R}{\Q+\alpha k}<n\le\frac{R+\alpha k^2}{\Q+\alpha k},\quad \alpha (n-k)<Q\le\frac{R-n\Q}{k}\right\}.
$$
\begin{enumerate}
  \item Вычислим сумму $\sum\limits_{(n,Q)\in\Omega_{21}}(\Q+Q+\alpha k)$.
\begin{gather*}
\sum_{U_1<n\le\frac{R}{\Q+\alpha k}}\sum_{\alpha (n-k)<Q\le\alpha n}(\Q+Q+\alpha k)=\\
\sum_{U_1<n\le\frac{R}{\Q+\alpha k}}\left(\alpha^2nk+\alpha^2\frac{k^2}{2}+\alpha k\Q+
O\left(\alpha n\right)+O\left(\Q+\alpha k\right)\right)=\\=
\left(\alpha^2\frac{k^2}{2}+\alpha k\Q\right)\left(\frac{R}{\Q+\alpha k}-U_1\right)+
\frac{\alpha^2k}{2}\left(\frac{R^2}{(\Q+\alpha k)^2}-U_1^2\right)+
O\left(\alpha^2k\frac{R}{\Q+\alpha k}\right)+\\+
O\left(\alpha^2k^2+\alpha k\Q\right)+
O\left(\alpha\frac{R^2}{(\Q+\alpha k)^2}\right)+
O\left(R\right)=\\=
\alpha k\left(R-U_1(\Q+\alpha k)\right)-\frac{\alpha^2k^2}{2}\left(\frac{R}{\Q+\alpha k}-U_1\right)+\\+
\frac{\alpha^2k}{2}\left(\frac{R^2}{(\Q+\alpha k)^2}-U_1^2\right)+
O\left(\alpha^2k\frac{R}{\Q+\alpha k}\right)+O\left(\alpha\frac{R^2}{(\Q+\alpha k)^2}\right).
\end{gather*}
Мы воспользовались тем, что из
$$\Q\le U_2-\alpha k$$
следует
$$\alpha k(\Q+\alpha k)\le\alpha^2k\frac{R}{\Q+\alpha k},$$
$$R\le\alpha\frac{R^2}{(\Q+\alpha k)^2}.$$
 \item Вычислим сумму $\sum\limits_{(n,Q)\in\Omega_{22}}(\Q+Q+\alpha k)$.
\begin{gather}
\sum_{\alpha (n-k)<Q\le\frac{R-n\Q}{k}}(\Q+Q+\alpha k)=(\Q+\alpha k)\left(\frac{R-n\Q}{k}-\alpha (n-k)\right)+
O\left(\Q+\alpha k\right)+\notag\\+\frac{1}{2}\left(\frac{(R-n\Q)^2}{k^2}-\alpha^2 (n-k)^2\right)+
O\left(\frac{R-n\Q}{k}\right)=\notag\\=\frac{R+\alpha k^2}{k}\left(\frac{R+\alpha k^2}{2k}+\Q\right)-
n\frac{\Q}{k}\left(\frac{R+\alpha k^2}{k}+\Q+\alpha k\right)+\frac{n^2}{2}\frac{\Q^2-\alpha^2k^2}{k^2}+\notag\\
O\left(\Q+\alpha k\right)+O\left(\frac{R-n\Q}{k}\right).\label{sum22}
\end{gather}
При суммировании по n воспользуемся следствием \ref{ostatok} (см. Приложение). В данном случае определим функцию $f(x)$ следующим образом
\begin{equation}
\label{f2}
f(x)=\frac{R+\alpha k^2}{k}\left(\frac{R+\alpha k^2}{2k}+\Q\right)-
x\frac{\Q}{k}\left(\frac{R+\alpha k^2}{k}+\Q+\alpha k\right)+
\frac{x^2}{2}\frac{\Q^2-\alpha^2k^2}{k^2}.
\end{equation}
Используя предыдущее равенство, легко получить другое выражение для функции $f(x)$
\begin{equation}
\label{f2v2}
f(x)=(\Q+\alpha k)\left(\frac{R-x\Q}{k}-\alpha (x-k)\right)+
\frac{1}{2}\left(\frac{(R-x\Q)^2}{k^2}-\alpha^2 (x-k)^2\right)
\end{equation}
Дифференцируя формулу \eqref{f2} получаем
$$
f’(x)=-\frac{\Q}{k}\left(\frac{R}{k}+\alpha k\right)-\frac{\Q}{k}\left(\Q+\alpha k\right)+
x\left(\frac{\Q}{k}-\alpha\right)\frac{\Q+\alpha k}{k}
$$
Следовательно,
\begin{gather*}
f’\left(\frac{R}{\Q+\alpha k}\right)=-\frac{\Q}{k}\left(\frac{R}{k}+\alpha k\right)-\frac{\Q}{k}\left(\Q+\alpha k\right)+\frac{R}{k}\left(\frac{\Q}{k}-\alpha\right)=\\=-\Q\alpha-
\frac{\Q}{k}\left(\Q+\alpha k\right)-\frac{R\alpha}{k}<0
\end{gather*}
и
\begin{gather*}
f’\left(\frac{R+\alpha k^2}{\Q+\alpha k}\right)=-\frac{\Q}{k}\left(\frac{R}{k}+\alpha k\right)
-\frac{\Q}{k}\left(\Q+\alpha k\right)+\frac{R+\alpha k^2}{k}\left(\frac{\Q}{k}-\alpha\right)=\\=
-\frac{\Q}{k}\left(\Q+\alpha k\right)-\frac{R+\alpha k^2}{k}\alpha<0.
\end{gather*}
Получаем, что $f(x)$ монотонно убывает  на промежутке
$$\frac{R}{\Q+\alpha k}<x\le\frac{R+\alpha k^2}{\Q+\alpha k}.$$
Из формулы \eqref{f2v2} получаем
$$f\left(\frac{R+\alpha k^2}{\Q+\alpha k}\right)=0.$$
Следовательно,
\begin{gather*}
\sum\limits_{(n,Q)\in\Omega_{22}}(\Q+Q+\alpha k)=\int_{\frac{R}{\Q+\alpha k}}^{\frac{R+\alpha k^2}{\Q+\alpha k}}f(x)dx
+O\left(f\left(\frac{R}{\Q+\alpha k}\right)\right)+\\+
\sum_{\frac{R}{\Q+\alpha k}<n\le\frac{R+\alpha k^2}{\Q+\alpha k}}\Biggl(
O\left(\Q+\alpha k\right)+O\left(\frac{R-n\Q}{k}\right)\Biggl).
\end{gather*}
Вычислим интеграл, используя формулу \eqref{f2}. Интегрируя каждое слагаемое формулы \eqref{f2}, получаем
\begin{gather*}
\int_{\frac{R}{\Q+\alpha k}}^{\frac{R+\alpha k^2}{\Q+\alpha k}}f(x)dx=
\frac{\alpha k^2}{\Q+\alpha k}\frac{R+\alpha k^2}{k}\left(\frac{R+\alpha k^2}{2k}+\Q\right)-\\-
\frac{\Q}{k}\left(\frac{R+\alpha k^2}{k}+\Q+\alpha k\right)\frac{(R+\alpha k^2)^2-R^2}{2(\Q+\alpha k)^2}+
\frac{\Q^2-\alpha^2k^2}{k^2}\frac{(R+\alpha k^2)^3-R^3}{6(\Q+\alpha k)^3}.
\end{gather*}
Преобразовывая отдельно каждое слагаемое, получаем
\begin{gather*}
\int_{\frac{R}{\Q+\alpha k}}^{\frac{R+\alpha k^2}{\Q+\alpha k}}f(x)dx=\alpha k\frac{R+\alpha k^2}{\Q+\alpha k}
\left(\Q+\alpha k+\frac{R-\alpha k^2}{2k}\right)-\\-\left(\frac{\Q+\alpha k}{k}-\alpha\right)
\left(\frac{R+\alpha k^2}{k}+\Q+\alpha k\right)\frac{(R+\alpha k^2)^2-R^2}{2(\Q+\alpha k)^2}+\\+
\frac{(R+\alpha k^2)^3-R^3}{6(\Q+\alpha k)^2k^2}\left(\Q+\alpha k-2\alpha k\right).
\end{gather*}
Раскрывая скобки и группируя слагаемые по одинаковым степеням $\left(\Q+\alpha k\right)$, получаем
\begin{gather*}
\int_{\frac{R}{\Q+\alpha k}}^{\frac{R+\alpha k^2}{\Q+\alpha k}}f(x)dx=
\alpha k(R+\alpha k^2)-\frac{(R+\alpha k^2)^2-R^2}{2k}+\\+\frac{1}{\Q+\alpha k}
\Biggl(
\frac{\alpha(R+\alpha k^2)(R-\alpha k^2)}{2}-\left(\frac{R+\alpha k^2}{k^2}-\alpha\right)\frac{(R+\alpha k^2)^2-R^2}{2}
+\frac{(R+\alpha k^2)^3-R^3}{6k^2}
\Biggl)+\\+
\frac{1}{(\Q+\alpha k)^2}
\Biggl(
\alpha\frac{R+\alpha k^2}{k}\frac{(R+\alpha k^2)^2-R^2}{2}-\alpha\frac{(R+\alpha k^2)^3-R^3}{3k}
\Biggl)=\\=
\frac{\alpha^2k^3}{2}-\frac{1}{\Q+\alpha k}\frac{\alpha^3k^4}{3}+\frac{1}{(\Q+\alpha k)^2}
\left(R\frac{\alpha^3k^3}{2}+\frac{\alpha^4k^5}{6}\right).
\end{gather*}
Вычислим асимптотический член суммы $\sum\limits_{(n,Q)\in\Omega_{22}}(\Q+Q+\alpha k)$.
\begin{enumerate}
  \item
Используя формулу \eqref{f2v2}, получаем
\begin{gather*}
f\left(\frac{R}{\Q+\alpha k}\right)=
(\Q+\alpha k)\alpha k+\frac{\alpha k}{2}\left(\frac{2R\alpha}{\Q+\alpha k}-\alpha k\right)=\\=
\frac{\alpha^2k^2}{2}+\alpha k\Q+\frac{R\alpha^2k}{\Q+\alpha k},
\end{gather*}
  \item
$$\sum_{\frac{R}{\Q+\alpha k}<n\le\frac{R+\alpha k^2}{\Q+\alpha k}}(\Q+\alpha k)=O(\alpha k^2),$$
  \item
\begin{gather*}
\sum_{\frac{R}{\Q+\alpha k}<n\le\frac{R+\alpha k^2}{\Q+\alpha k}}O\left(\frac{R-n\Q}{k}\right)=
O\left(\frac{R\alpha k}{\Q+\alpha k}\right).
\end{gather*}
\end{enumerate}
Учитывая неравенства  $k\le U_1, \Q\le U_2-\alpha k$  и формулы \eqref{U1U2} и \eqref{U1U2R}, получаем
$$
\frac{\alpha^2k^2}{2}+\alpha k\Q+\frac{R\alpha^2k}{\Q+\alpha k}\le \alpha k\left(\Q+\alpha k+
\frac{R\alpha}{\Q+\alpha k}\right)\le\frac{R\alpha^2k}{\Q+\alpha k}
$$
и
$$\alpha k^2\le\frac{R\alpha k}{\Q+\alpha k}.$$
Следовательно,
\begin{gather*}
\sum\limits_{(n,Q)\in\Omega_{22}}(\Q+Q+\alpha k)=\frac{\alpha^2k^3}{2}-\frac{1}{\Q+\alpha k}\frac{\alpha^3k^4}{3}
+\frac{1}{(\Q+\alpha k)^2}\left(R\frac{\alpha^3k^3}{2}+\frac{\alpha^4k^5}{6}\right)+\\+
O\left(\frac{R\alpha k}{(\Q+\alpha k)}+\frac{R\alpha^2k}{\Q+\alpha k}\right).
\end{gather*}
\end{enumerate}
Получаем
\begin{gather*}
\left(\sum\limits_{(n,Q)\in\Omega_{21}}+\sum\limits_{(n,Q)\in\Omega_{22}}\right)(\Q+Q+\alpha k)=
\alpha k\left(R-U_1(\Q+\alpha k)\right)-\frac{\alpha^2k^2}{2}\left(\frac{R}{\Q+\alpha k}-U_1\right)+\\
\frac{\alpha^2k}{2}\left(\frac{R^2}{(\Q+\alpha k)^2}-U_1^2\right)+
\frac{\alpha^2k^3}{2}-\frac{1}{\Q+\alpha k}\frac{\alpha^3k^4}{3}
+\frac{1}{(\Q+\alpha k)^2}\left(R\frac{\alpha^3k^3}{2}+\frac{\alpha^4k^5}{6}\right)+\\+
O\left(\frac{\alpha R^2}{(\Q+\alpha k)^2}+\frac{R\alpha^2k}{\Q+\alpha k}\right).
\end{gather*}
Полученное выражение необходимо просуммировать по
$$(k,\Q)\in\Omega_{2}=\left\{k\le U_1, \quad \Q\le U_2-\alpha k\right\},\quad\tau\mid k,\quad\tau\mid \Q.$$
Получаем
\begin{gather*}
\sum_{k\le\frac{U_1}{\tau}}\sum_{\Q\le\alpha\left(\frac{U_1}{\tau}-k\right)}\Biggl(
\alpha\tau Rk+\frac{\alpha^2\tau^2U_1}{2}k^2-\frac{\alpha^2\tau U_1^2}{2}k+\frac{\alpha^2\tau^3}{2}k^3-
\alpha\tau^2U_1k(\Q+\alpha k)-\\-\frac{1}{\Q+\alpha k}\left(\frac{R\alpha^2\tau}{2}k^2+\frac{\alpha^3\tau^3}{3}k^4\right)+
\frac{1}{(\Q+\alpha k)^2}\left(\frac{R^2\alpha^2}{2\tau}k+\frac{R\alpha^3\tau}{2}k^3+\frac{\alpha^4\tau^3}{6}k^5\right)+\\+
O\left(\frac{\alpha R^2}{\tau^2(\Q+\alpha k)^2}+\frac{R\alpha^2k}{\Q+\alpha k}\right)
\Biggl).
\end{gather*}
Используя  следствие \ref{ostatok} (a) (см. Приложение), получаем
\begin{gather*}
\sum_{\Q\le\alpha\left(\frac{U_1}{\tau}-k\right)}(\Q+\alpha k)=\frac{\alpha^2}{2}\left(\frac{U_1^2}{\tau^2}-k^2\right)+
O\left(\frac{\alpha U_1}{\tau}\right),\\
\sum_{\Q\le\alpha\left(\frac{U_1}{\tau}-k\right)}\frac{1}{\Q+\alpha k}=\log\frac{U_1}{\tau k}+
O\left(\frac{1}{\alpha k}\right),\\
\sum_{\Q\le\alpha\left(\frac{U_1}{\tau}-k\right)}\frac{1}{(\Q+\alpha k)^2}=\frac{1}{\alpha k}-\frac{\tau}{\alpha U_1}+
O\left(\frac{1}{\alpha^2 k^2}\right),\\
\sum_{\Q\le\alpha\left(\frac{U_1}{\tau}-k\right)}1=\frac{\alpha U_1}{\tau}-\alpha k+O(1).
\end{gather*}
Применяя эти соотношения, получаем
\begin{gather*}
\sum_{k\le\frac{U_1}{\tau}}\Biggl(
\frac{R^2\alpha}{2\tau}+k\left(\alpha^2RU_1-\frac{\alpha^3U_1^3}{2}-\frac{\alpha^3U_1^3}{2}-\frac{R^2\alpha}{2U_1}\right)
+k^2\left(\frac{\alpha^3\tau U_1^2}{2}-\alpha^2\tau R+\frac{\alpha^3\tau U_1^2}{2}+\frac{R\alpha^2\tau}{2}\right)+\\+
k^3\left(\frac{\alpha^3\tau^2 U_1}{2}-\frac{\alpha^3\tau^2 U_1}{2}+\frac{\alpha^3\tau^2 U_1}{2}-
\frac{R\alpha^2\tau^2 }{2}\right)+k^4\left(-\frac{\alpha^3\tau^3}{2}+\frac{\alpha^3\tau^3}{6}\right)-\\
-\frac{\alpha^3\tau^4}{6U_1}k^5-
\left(\frac{R\alpha^2\tau}{2}k^2+\frac{\alpha^3\tau^3}{3}k^4\right)\log\frac{U_1}{\tau k}
\Biggl)+O\left(\frac{R^2}{\tau}(1+\alpha)+\frac{R^2}{\tau}\log\frac{U_1}{\tau}\right).
\end{gather*}
При суммировании асимптотических слагаемых мы воспользовались пунктами (a) и (b) следствия \ref{sl} (см. Приложение) и оценкой
$$\sum_{k\le\frac{U_1}{\tau}}\frac{1}{k}=O(\log\frac{U_1}{\tau}).$$
Используя пункты (a) и (b) следствия \ref{sl} (см. Приложение), получаем
\begin{gather*}
\sum_{k\le\frac{U_1}{\tau}}\Biggl(\frac{R^2\alpha}{2\tau}-\frac{\alpha^3U_1^3}{2}k+\frac{\alpha^3\tau U_1^2}{2}k^2-
\frac{\alpha^3\tau^3}{3}k^4-\frac{\alpha^3\tau^4}{6U_1}k^5-
\left(\frac{R\alpha^2\tau}{2}k^2+\frac{\alpha^3\tau^3}{3}k^4\right)\log\frac{U_1}{\tau k}
\Biggl)+\\+O\left(\frac{R^2}{\tau}(1+\alpha)+\frac{R^2}{\tau}\log\frac{U_1}{\tau}\right)=
\frac{R^{5/2}\sqrt{\alpha}}{\tau^2}\left(\frac{1}{2}-\frac{1}{4}+\frac{1}{6}-\frac{1}{15}
-\frac{1}{36}-\frac{1}{18}-\frac{1}{75}\right)\\+O\left(\frac{R^2}{\tau}(1+\alpha)+\frac{R^2}{\tau}\log\frac{U_1}{\tau}\right).
\end{gather*}
Следовательно,
\begin{gather*}
\sum_{\substack{(k,\Q)\in\Omega_{2}\\\tau\mid k,\tau\mid \Q}}
\left(\sum\limits_{(n,Q)\in\Omega_{21}}+\sum\limits_{(n,Q)\in\Omega_{22}}\right)(\Q+Q+\alpha k)=\\=
\frac{19}{75}\frac{R^{5/2}\sqrt{\alpha}}{\tau^2}+
O\left(\frac{R^2}{\tau}(1+\alpha)+\frac{R^2}{\tau}(1+\alpha)\log\frac{U_1}{\tau}\right).
\end{gather*}
Подставляя полученное выражение в формулу
$$\Sigma_{21}=\sum_{\tau|\delta}\mathop{{\sum}^*}_{z=1}^{\tau}\sum_{\substack{(k,\Q)\in\Omega_{2}\\\tau\mid k,\tau\mid \Q}}
\left(\sum\limits_{(n,Q)\in\Omega_{21}}+\sum\limits_{(n,Q)\in\Omega_{22}}\right)(\Q+Q+\alpha k),$$
получаем утверждение леммы.
\end{proof}
\end{Le}

\subsection{Случай 3}
\begin{Le}
\label{31}
Справедлива следующая асимптотическая формула
\begin{gather*}
\Sigma_{31}=\left(\frac{2\pi}{15}-\frac{2\log 2}{5}-\frac{91}{900}\right)R^{5/2}\sqrt{\alpha}
\sum_{\tau|\delta}\frac{\varphi(\tau)}{\tau^2}+\sum_{\tau|\delta}\varphi(\tau)
O\left(\frac{R^2}{\tau}\left(1+\alpha\right)\right),
\end{gather*}
где $\Sigma_{31}$ определена в \eqref{slu3}.
\begin{proof}
Напомним, что внутреннее суммирование ведется по области
$$
\Omega_{31}=\left\{U_1<n\le\frac{R+\alpha k^2}{\Q+\alpha k},\quad \alpha (n-k)<Q\le\frac{R-n\Q}{k} \right\}.
$$
Заметим, что внутреннее суммирование по $Q$ аналогично одному из случаев в лемме \ref{21}. По-этому
\begin{gather*}
\sum_{\alpha (n-k)<Q\le\frac{R-n\Q}{k}}(\Q+Q+\alpha k)=f(n)+
O\left(\Q+\alpha k\right)+O\left(\frac{R-n\Q}{k}\right),
\end{gather*}
где $f(n)$ определяется формулой \eqref{f2}. При суммировании по n воспользуемся \mbox{следствием \ref{ostatok}}
(см. Приложение). Определим функцию $F_1(k,\Q)$
\begin{gather}
F_1(k,\Q)=\frac{\Q^2-\alpha^2k^2}{6k^2}\left(\frac{(R+\alpha k^2)^3}{(\Q+\alpha k)^3}-U_1^3\right)-\notag\\-
\frac{\Q}{2k}\left(\frac{R+\alpha k^2}{k}+\Q+\alpha k\right)\left(\frac{(R+\alpha k^2)^2}{(\Q+\alpha k)^2}-U_1^2\right)+
\notag\\+\frac{R+\alpha k^2}{k}\left(\frac{R+\alpha k^2}{2k}+\Q\right)\left(\frac{R+\alpha k^2}{\Q+\alpha k}-U_1\right),\label{kQ}
\end{gather}
тогда
\begin{gather*}
\sum\limits_{(n,Q)\in\Omega_{31}}(\Q+Q+\alpha k)=F_1(k,\Q)+O\left(f(U_1)\right)+
\sum_{U_1<n\le\frac{R+\alpha k^2}{\Q+\alpha k}}\left(O\left(\Q+\alpha k\right)+O\left(\frac{R-n\Q}{k}\right)\right).
\end{gather*}
Полученное выражение необходимо просуммировать по
$$
(k,\Q)\in\Omega_{3}=\left\{k\le U_1, \quad U_2-\alpha k<\Q\le U_2-\alpha k+\frac{\alpha k^2}{U_1}\right\},
\quad\tau\mid k,\quad\tau\mid \Q.
$$
Учитывая эти ограничения и формулы \eqref{U1U2}, \eqref{U1U2R}, получаем
\begin{gather*}
\sum_{U_1<n\le\frac{R+\alpha k^2}{\Q+\alpha k}}\left(O\left(\Q+\alpha k\right)+O\left(\frac{R-n\Q}{k}\right)\right)\ll
O\left(R\right)+O\left(\frac{R^2}{k(\Q+\alpha k)}\right),
\end{gather*}
\begin{gather*}
f(U_1)=\frac{R+\alpha k^2}{k}\left(\frac{R+\alpha k^2}{2k}+\Q\right)-
U_1\frac{\Q}{k}\left(\frac{R+\alpha k^2}{k}+\Q+\alpha k\right)+\\+
\frac{U_1^2}{2}\left(\frac{\Q^2-\alpha^2k^2}{k^2}\right)\ll\frac{R}{k}\frac{R}{k}+U_1\frac{\Q}{k}\frac{R}{k}+
\frac{U_1^2U_2^2}{k^2}=\frac{R^2}{k^2}.
\end{gather*}
Получаем
\begin{gather*}
O\left(f(U_1)\right)+
\sum_{u_1<n\le\frac{R+\alpha k^2}{\Q+\alpha k}}\left(O\left(\Q+\alpha k\right)+O\left(\frac{R-n\Q}{k}\right)\right)\le
O\left(\frac{R^2}{k(\Q+\alpha k)}\right)+O\left(\frac{R^2}{k^2}\right).
\end{gather*}
Определим функцию $F_2(k,\Q)$ следующим образом. Разложим в $F_1(k,\Q)$ слагаемые по степеням $(\Q+\alpha k)$
\begin{gather*}
F_1(k,\Q)=\frac{(R+\alpha k^2)^3}{k^2(\Q+\alpha k)^2}\left(\frac{\Q-\alpha k}{6}-
\frac{\Q}{2}\right)+\frac{(R+\alpha k^2)^3}{2k^2(\Q+\alpha k)}+\\+
\frac{(R+\alpha k^2)^3}{2k}\frac{\Q}{\Q+\alpha k}-\frac{\Q^2-\alpha^2k^2}{6k^2}U_1^3+\\
+\frac{U_1^2}{2k}\left(\Q(\Q+\alpha k)+\frac{R+\alpha k^2}{k}\Q\right)-
U_1\frac{R+\alpha k^2}{k}\left(\frac{R+\alpha k^2}{2k}+\Q\right)=\\=
-\frac{(R+\alpha k^2)^3}{3k^2(\Q+\alpha k)}+\frac{(R+\alpha k^2)^3\alpha k}{6k^2(\Q+\alpha k)^2}+
\frac{(R+\alpha k^2)^3}{2k^2(\Q+\alpha k)}+\\+\frac{(R+\alpha k^2)^2}{2k}-\frac{(R+\alpha k^2)^2\alpha k}{2k(\Q+\alpha k)}-
\frac{\Q^2-\alpha^2k^2}{6k^2}U_1^3+\\+
\frac{U_1^2}{2k}\left((\Q+\alpha k)^2+\frac{R}{k}(\Q+\alpha k)-(R+\alpha k^2)\alpha\right)-
U_1\frac{R+\alpha k^2}{k}\left(\frac{R-\alpha k^2}{2k}+(\Q+\alpha k)\right).
\end{gather*}
и пусть
$$F_2(k,\Q)=F_1(\tau k,\tau\Q).$$
Тогда, проделав элементарные преобразования с $F_1(k,\Q)$, получаем
\begin{gather*}
F_2(k,\Q)=\frac{(R+\alpha \tau^2k^2)^3}{6\tau^3k^2}\frac{1}{\Q+\alpha k}+\\+
\frac{(R+\alpha \tau^2k^2)^3\alpha}{6\tau^3k}\frac{1}{(\Q+\alpha k)^2}-
\frac{(R+\alpha \tau^2k^2)^2\alpha}{2\tau}\frac{1}{\Q+\alpha k}+
\frac{(R+\alpha \tau^2k^2)^2}{2\tau k}-\frac{\Q^2-\alpha^2k^2}{6k^2}U_1^3+\\+
\frac{U_1^2}{2\tau k}\left(\tau^2(\Q+\alpha k)^2+\frac{R}{k}(\Q+\alpha k)-(R+\alpha\tau^2k^2)\alpha\right)-
U_1\frac{R+\alpha\tau^2k^2}{\tau k}\left(\frac{R-\alpha\tau^2k^2}{2\tau k}+\tau(\Q+\alpha k)\right).
\end{gather*}
Следовательно,
\begin{gather*}
\sum_{\substack{(k,\Q)\in\Omega_{3}\\\tau\mid k,\tau\mid \Q}}\sum\limits_{(n,Q)\in\Omega_{31}}(\Q+Q+\alpha k)=\\=
\sum_{k\le\frac{U_1}{\tau}}\sum_{\frac{U_2}{\tau}-\alpha k<\Q\le \frac{U_2}{\tau}-\alpha k+\frac{\alpha \tau k^2}{U_1}}
\Biggl(F_2(k,\Q)+O\left(\frac{R^2}{\tau^2k(\Q+\alpha k)}\right)+O\left(\frac{R^2}{\tau^2k^2}\right)\Biggl).
\end{gather*}
Легко проверить, что функция $F_2(k,\Q)$ не является монотонной по переменной $\Q$. По-этому непосредственное применение следствия \ref{ostatok} невозможно. Однако, в данном случаи утверждение следствия \ref{ostatok} оказывается верным. Заметим, что уравнение
$$\frac{\partial F_2(k,\Q)}{\partial\Q}=0$$
имеет конечное число решений на нашем полуинтервале. Следовательно,
\begin{gather*}
\left|\int\limits_{\frac{U_2}{\tau}-\alpha k}^{\frac{U_2}{\tau}-\alpha k+\frac{\alpha \tau k^2}{U_1}}
\rho(x)F’_2(k,x)dx\right|\ll\max_{\frac{U_2}{\tau}-\alpha k<\Q\le \frac{U_2}{\tau}-\alpha k+\frac{\alpha \tau k^2}{U_1}} \left|F_2(k,\Q)\right|.
\end{gather*}
Функция $F_1(k,\Q)$ является главным членом в асимптотической формуле для
$$\sum\limits_{(n,Q)\in\Omega_{31}}(\Q+Q+\alpha k).$$
Поэтому $F_1(k,\Q)>0$, следовательно, и $F_2(k,\Q)>0$. Так как
$$U_2<\Q+\alpha k\le 2U_2,$$
а область суммирования максимальна при $\Q=U_2-\alpha k$, то
$$\max_{U_2-\alpha k<\Q\le U_2-\alpha k+\frac{\alpha k^2}{U_1}} F_1(k,\Q)\ll F_1\left(k,U_2-\alpha k\right).$$
Следовательно,
\begin{gather*}
\max_{\frac{U_2}{\tau}-\alpha k<\Q\le \frac{U_2}{\tau}-\alpha k+\frac{\alpha \tau k^2}{U_1}} \left|F_2(k,\Q)\right|
\ll F_2\left(k,\frac{U_2}{\tau}-\alpha k\right).
\end{gather*}
Из формулы \eqref{kQ} следует,  что
$$F_2\left(k,\frac{U_2}{\tau}-\alpha k+\frac{\alpha \tau k^2}{U_1}\right)=F_1\left(\tau k,
U_2-\alpha\tau k+\frac{\alpha \tau^2 k^2}{U_1}\right)=0,$$
получаем
\begin{gather*}
\sum_{\frac{U_2}{\tau}-\alpha k<\Q\le \frac{U_2}{\tau}-\alpha k+\frac{\alpha \tau k^2}{U_1}}
F_2(k,\Q)=\int\limits_{\frac{U_2}{\tau}-\alpha k}^{\frac{U_2}{\tau}-\alpha k+\frac{\alpha \tau k^2}{U_1}}F_2(k,\Q)d\Q+
O\left(F_2(k,\frac{U_2}{\tau}-\alpha k)\right).
\end{gather*}
Следовательно,
\begin{gather*}
\sum_{\substack{(k,\Q)\in\Omega_{3}\\\tau\mid k,\tau\mid \Q}}\sum\limits_{(n,Q)\in\Omega_{31}}(\Q+Q+\alpha k)=
\sum_{k\le\frac{U_1}{\tau}}\Biggl(
\int\limits_{\frac{U_2}{\tau}-\alpha k}^{\frac{U_2}{\tau}-\alpha k+\frac{\alpha \tau k^2}{U_1}}F_2(k,\Q)d\Q+
O\left(F_2(k,\frac{U_2}{\tau}-\alpha k)\right)\Biggl)+\\+
\sum_{k\le\frac{U_1}{\tau}}\sum_{\frac{U_2}{\tau}-\alpha k<\Q\le \frac{U_2}{\tau}-\alpha k+\frac{\alpha \tau k^2}{U_1}}
\Biggl(O\left(\frac{R^2}{\tau^2k(\Q+\alpha k)}\right)+O\left(\frac{R^2}{\tau^2k^2}\right)\Biggl).
\end{gather*}
Просуммируем полученные асимптотические слагаемые.
\begin{enumerate}
  \item
Используя соотношение $F_2(k,\Q)=F_1(\tau k,\tau\Q)$ и формулу \eqref{kQ}, в которой мы заменим
$$
\frac{\Q}{2k}\left(\frac{R+\alpha k^2}{k}+\Q+\alpha k\right)=\frac{1}{2k}
\left((\Q+\alpha k)^2+\frac{R}{k}\Q-\alpha^2k^2\right),
$$
получаем
\begin{gather*}
F_2\left(k,\frac{U_2}{\tau}-\alpha k\right)=\frac{1}{6k^2}\frac{U_2}{\tau}(\frac{U_2}{\tau}-2\alpha k)
\left(\frac{(R+\alpha\tau^2k^2)^3}{U_2^3}-U_1^3\right)-\\-\frac{1}{2}\left(\frac{R}{\tau k^2}(\frac{U_2}{\tau}-\alpha k)-\alpha^2\tau k+\frac{U_2^2}{\tau k}\right)\left(\frac{(R+\alpha\tau^2k^2)^2}{U_2^2}-U_1^2\right)+\\+
\left(\frac{R^2-\alpha^2\tau^4k^4}{2\tau^2k^2}+\frac{R+\alpha\tau^2k^2}{\tau k}U_2\right)
\left(\frac{R+\alpha\tau^2k^2}{U_2}-U_1\right)
\end{gather*}
Применяя формулы разности квадратов и кубов и формулы \eqref{U1U2R}, \eqref{U1U2}, получаем
\begin{gather*}
F_2\left(k,\frac{U_2}{\tau}-\alpha k\right)=\frac{\alpha^2\tau}{6}\left(\frac{U_1}{\tau}-2k\right)
\left(\left(U_1+\frac{\tau^2k^2}{U_1}\right)^2+\frac{R}{\alpha}+\tau^2k^2+U_1^2\right)-\\-
\left(\frac{R\alpha}{\tau k^2}\left(\frac{U_1}{\tau}-k\right)+\frac{\alpha^2}{\tau k}(U_1^2-\tau^2k^2)\right)
\left(\tau^2k^2+\frac{\tau^4k^4}{2U_1^2}\right)+
\frac{R^2-\alpha^2\tau^4k^4}{2U_1}+\left(R+\alpha\tau^2k^2\right)\alpha\tau k.
\end{gather*}
Так как $k\le U_1/\tau$, легко получить
$$O\left(F_2(k,\frac{U_2}{\tau}-\alpha k)\right)<O(RU_1\alpha).$$
Следовательно,
$$\sum_{k\le\frac{U_1}{\tau}}O\left(F_2(k,\frac{U_2}{\tau}-\alpha k)\right)=O\left(\frac{R^2}{\tau}\right).$$
  \item
Используя асимптотические формулу
$$
\sum_{\frac{U_2}{\tau}-\alpha k<\Q\le \frac{U_2}{\tau}-\alpha k+\frac{\alpha \tau k^2}{U_1}}\frac{1}{\Q+\alpha k}=
\log\left(1+\frac{\alpha\tau^2}{R}k^2\right)+O\left(\frac{\tau}{U_2}\right),
$$
$$\sum_{n\le a}\frac{1}{n}\log\left(1+\frac{n^2}{a^2}\right)=O(1),$$
получаем
\begin{gather*}
\sum_{k\le\frac{U_1}{\tau}}\sum_{\frac{U_2}{\tau}-\alpha k<\Q\le \frac{U_2}{\tau}-\alpha k+\frac{\alpha \tau k^2}{U_1}}
O\left(\frac{R^2}{\tau^2k(\Q+\alpha k)}\right)=\\=\sum_{k\le\frac{U_1}{\tau}}
O\left(\frac{R^2}{\tau^2k}\log\left(1+\frac{\alpha\tau^2}{R}k^2\right)\right)=O\left(\frac{R^2}{\tau^2}\right).
\end{gather*}
  \item
Очевидно, что
\begin{gather*}
\sum_{k\le\frac{U_1}{\tau}}\sum_{\frac{U_2}{\tau}-\alpha k<\Q\le \frac{U_2}{\tau}-\alpha k+\frac{\alpha \tau k^2}{U_1}}
O\left(\frac{R^2}{\tau^2k^2}\right)=O\left(\frac{R^2\alpha}{\tau^2}\right).
\end{gather*}
\end{enumerate}
Следовательно,
\begin{gather*}
\sum_{\substack{(k,\Q)\in\Omega_{3}\\\tau\mid k,\tau\mid \Q}}\sum\limits_{(n,Q)\in\Omega_{31}}(\Q+Q+\alpha k)=
\sum_{k\le\frac{U_1}{\tau}}
\int\limits_{\frac{U_2}{\tau}-\alpha k}^{\frac{U_2}{\tau}-\alpha k+\frac{\alpha \tau k^2}{U_1}}F_2(k,\Q)d\Q+
O\left(\frac{R^2}{\tau}(1+\alpha)\right).
\end{gather*}
Вычислим интеграл.
\begin{equation}\label{F2int1}
\int\limits_{\frac{U_2}{\tau}-\alpha k}^{\frac{U_2}{\tau}-\alpha k+\frac{\alpha \tau k^2}{U_1}}\frac{1}{(\Q+\alpha k)}d\Q=\log\left(1+\frac{\alpha\tau^2}{R}k^2\right),
\end{equation}
\begin{equation}\label{F2int2}
\int\limits_{\frac{U_2}{\tau}-\alpha k}^{\frac{U_2}{\tau}-\alpha k+\frac{\alpha \tau k^2}{U_1}}\frac{1}{(\Q+\alpha k)^2}d\Q=\frac{\tau}{U_2}-\frac{U_1\tau}{R+\alpha\tau^2 k^2}=\frac{\alpha\tau^3k^2}{U_2(R+\alpha\tau^2 k^2)},
\end{equation}
\begin{gather}
\int\limits_{\frac{U_2}{\tau}-\alpha k}^{\frac{U_2}{\tau}-\alpha k+\frac{\alpha \tau k^2}{U_1}}
\Biggl(-\frac{\Q^2-\alpha^2k^2}{6k^2}U_1^3+
\frac{U_1^2}{2}\left(\frac{\tau}{k}(\Q+\alpha k)^2+\frac{R}{\tau k^2}(\Q+\alpha k)-(R+\alpha\tau^2k^2)\frac{\alpha}{\tau k}\right)-\notag\\-
U_1\frac{R+\alpha\tau^2k^2}{\tau k}\left(\frac{R-\alpha\tau^2k^2}{2\tau k}+\tau(\Q+\alpha k)\right)\Biggl)d\Q=
\notag\\=
-\frac{U_1^3}{6k^2}\Biggl(\frac{\alpha^3U_1}{\tau}k^2+\frac{\tau\alpha^3}{U_1}k^4+\frac{\tau^3\alpha^3}{3U_1^3}k^6-
\frac{\tau^2\alpha^3}{U_1^2}k^5-2\alpha^3k^3
\Biggl)+\frac{U_1^2}{2}\Biggl(
\frac{R}{2\tau k^2}\left(\frac{\tau^2\alpha^2}{U_1^2}k^4+2\alpha^2k^2\right)+\notag\\+
\frac{\tau}{3k}\left(\frac{3}{\tau}\alpha^3U_1k^2+3\tau\frac{\alpha^3}{U_1}k^4+\frac{\tau^3\alpha^3}{U_1^3}k^6\right)-
(R+\alpha\tau^2k^2)\frac{\alpha}{\tau k}\frac{\alpha\tau}{U_1}k^2\Biggl)
-\notag\\-
U_1\left(\frac{R^2-\alpha^2\tau^4k^4}{2\tau^2k^2}\frac{\alpha\tau}{U_1}k^2+
\frac{R+\alpha\tau^2k^2}{2k}\left(\frac{\tau^2\alpha^2}{U_1^2}k^4+2\alpha^2k^2\right)\right)
=\notag\\=
-\frac{R^2\alpha}{6\tau}-\frac{2}{3}R^{3/2}\alpha^{3/2}k+
\frac{1}{12}R\alpha^2\tau k^2-\frac{4}{3}R^{1/2}\alpha^{5/2}\tau^2k^3+
\frac{4}{9}\alpha^3\tau^3k^4-\frac{1}{3}R^{-1/2}\alpha^{7/2}\tau^4k^5.
\label{F2int3}
\end{gather}
Используя полученные соотношения \eqref{F2int1}, \eqref{F2int2}, \eqref{F2int3}, получаем
\begin{gather*}
\int\limits_{\frac{U_2}{\tau}-\alpha k}^{\frac{U_2}{\tau}-\alpha k+\frac{\alpha \tau k^2}{U_1}}F_2(k,\Q)d\Q=
\frac{(R+\alpha\tau^2 k^2)^3}{6\tau^3k^2}\log\left(1+\frac{\alpha\tau^2}{R}k^2\right)+
\frac{(R+\alpha\tau^2 k^2)^2}{6U_2}\alpha^2k-\\-
\frac{(R+\alpha \tau^2k^2)^2\alpha}{2\tau}\log\left(1+\frac{\alpha\tau^2}{R}k^2\right)+
\frac{(R+\alpha\tau^2k^2)^2}{2U_1}\alpha k-\frac{R^2\alpha}{6\tau}-\frac{2}{3}R^{3/2}\alpha^{3/2}k+\\+
\frac{1}{12}R\alpha^2\tau k^2-\frac{4}{3}R^{1/2}\alpha^{5/2}\tau^2k^3+
\frac{4}{9}\alpha^3\tau^3k^4-\frac{1}{3}R^{-1/2}\alpha^{7/2}\tau^4k^5.
\end{gather*}
Осталось просуммировать по k.
\begin{enumerate}
  \item
Обозначим $a=\frac{U_1}{\tau}$, тогда   $a^2=\frac{R}{\alpha\tau^2}.$ Следовательно,
$$\sum_{k\le\frac{U_1}{\tau}}\frac{(R+\alpha\tau^2 k^2)^3}{6\tau^3k^2}\log\left(1+\frac{\alpha\tau^2}{R}k^2\right)=
\frac{R^3}{6\tau^3}\sum_{k\le a}\frac{1}{k^2}\left(1+\frac{k^2}{a^2}\right)^3\log\left(1+\frac{k^2}{a^2}\right).$$
Применяя лемму \ref{euler} (см. Приложение) и формулы
\begin{gather*}
\int\frac{\log(1+x^2)}{x^2}dx=-\frac{\log(1+x^2)}{x}+2\arctan{x},\\
\int\log(1+x^2)dx=x\log(1+x^2)-2x+2\arctan{x},\\
\int x^2\log(1+x^2)dx=\frac{1}{3}\left(x^3\log(1+x^2)-\frac{2}{3}x^3+2x-2\arctan{x}\right),\\
\int x^4\log(1+x^2)dx=\frac{1}{5}\left(x^5\log(1+x^2)-\frac{2}{5}x^5+\frac{2}{3}x^3-2x+2\arctan{x}\right),
\end{gather*}
получаем
\begin{gather*}
\sum_{k\le a}\frac{1}{k^2}\left(1+\frac{k^2}{a^2}\right)^3\log\left(1+\frac{k^2}{a^2}\right)=
\frac{1}{a}\left(\frac{16}{5}\log2+\frac{8\pi}{5}-\frac{376}{75}\right)+O\left(\frac{1}{a^2}\right).
\end{gather*}
Следовательно,
\begin{gather*}
\sum_{k\le\frac{U_1}{\tau}}\frac{(R+\alpha\tau^2 k^2)^3}{6\tau^3k^2}\log\left(1+\frac{\alpha\tau^2}{R}k^2\right)=
\frac{R^{5/2}\sqrt{\alpha}}{\tau^2}\left(\frac{8}{15}\log2+\frac{4\pi}{15}-\frac{188}{225}\right)+
O\left(\frac{R^2\alpha}{\tau}\right).
\end{gather*}
  \item
Применяя пункт (a) следствия \ref{sl} (см. Приложение), получаем
\begin{gather*}
\sum_{k\le\frac{U_1}{\tau}}\frac{(R+\alpha\tau^2 k^2)^2}{6U_2}\alpha^2k =
\frac{R^2\alpha^2}{6U_2}\sum_{k\le\frac{U_1}{\tau}}\left(1+\frac{\alpha\tau^2}{R}k^2\right)^2k=
\frac{7}{36}\frac{R^{5/2}\sqrt{\alpha}}{\tau^2}+O\left(\frac{R^2\alpha}{\tau}\right).
\end{gather*}
  \item
Аналогично первому пункту получаем
\begin{gather*}
\sum_{k\le\frac{U_1}{\tau}}\frac{(R+\alpha \tau^2k^2)^2\alpha}{2\tau}\log\left(1+\frac{\alpha\tau^2}{R}k^2\right)=
\frac{R^2\alpha}{2\tau}\sum_{k\le a}\left(1+\frac{k^2}{a^2}\right)^2\log\left(1+\frac{k^2}{a^2}\right)=\\=
\frac{R^{5/2}\sqrt{\alpha}}{\tau^2}\left(\frac{14}{15}\log2+\frac{2\pi}{15}-\frac{164}{225}\right)+
O\left(\frac{R^2\alpha}{\tau}\right).
\end{gather*}
  \item
Применяя пункт (a) следствия \ref{sl} (см. Приложение), получаем
\begin{gather*}
\sum_{k\le\frac{U_1}{\tau}}\frac{(R+\alpha\tau^2k^2)^2}{2U_1}\alpha k=\frac{R^2\alpha}{2U_1}\sum_{k\le\frac{U_1}{\tau}}
\left(1+\frac{\alpha\tau^2}{R}k^2\right)^2k=
\frac{7}{12}\frac{R^{5/2}\sqrt{\alpha}}{\tau^2}+O\left(\frac{R^2\alpha}{\tau}\right).
\end{gather*}
  \item
Применяя пункт (a) следствия \ref{sl} (см. Приложение), получаем
\begin{gather*}
\sum_{k\le\frac{U_1}{\tau}}\Biggl(-\frac{R^2\alpha}{6\tau}-\frac{2}{3}R^{3/2}\alpha^{3/2}k+
\frac{1}{12}R\alpha^2\tau k^2-\frac{4}{3}R^{1/2}\alpha^{5/2}\tau^2k^3+\\+
\frac{4}{9}\alpha^3\tau^3k^4-\frac{1}{3}R^{-1/2}\alpha^{7/2}\tau^4k^5 \Biggl)=
-\frac{139}{180}\frac{R^{5/2}\sqrt{\alpha}}{\tau^2}+O\left(\frac{R^2\alpha}{\tau}\right).
\end{gather*}
\end{enumerate}
Следовательно,
\begin{gather*}
\sum_{\substack{(k,\Q)\in\Omega_{3}\\\tau\mid k,\tau\mid \Q}}\sum\limits_{(n,Q)\in\Omega_{31}}(\Q+Q+\alpha k)=
\left(\frac{2\pi}{15}-\frac{2\log 2}{5}-\frac{91}{900}\right)\frac{R^{5/2}\sqrt{\alpha}}{\tau^2}+
O\left(\frac{R^2}{\tau}(1+\alpha)\right).
\end{gather*}
Тем самым лемма доказана.
\end{proof}
\end{Le}

\subsection{Случай 4}
\begin{Le}
\label{41}
Справедлива следующая асимптотическая формула
\begin{gather*}
\Sigma_{41}=\left(\frac{11\pi}{120}-\frac{251}{60}\log 2+\frac{14}{5}\right)R^{5/2}\sqrt{\alpha}
\sum_{\tau|\delta}\frac{\varphi(\tau)}{\tau^2}+\\+\sum_{\tau|\delta}\varphi(\tau)
O\left(\frac{R^2}{\tau}\left(1+\log\frac{U_2}{\tau}\right)+\frac{R^2\alpha}{\tau^2}\log\frac{U_2}{\tau}\right),
\end{gather*}
где $\Sigma_{41}$ определена в \eqref{slu4}.
\begin{proof}
Вычислим отдельно суммы по областям
$$
\Omega_{41}=\left\{U_1<n\le\frac{R}{Q+\Q},\quad U_1<k\le n \right\},
$$
$$
\Omega_{42}=\left\{\frac{R}{Q+\Q}<n\le\frac{R+Q^2/\alpha}{Q+\Q},\quad U_1<k\le\frac{R-n\Q}{Q}\right\},
$$
$$
\Omega_{43}=\left\{U_1+\frac{Q}{\alpha}<n\le\frac{R+Q^2/\alpha}{Q+\Q},\quad U_1<k\le n-\frac{Q}{\alpha}\right\}.
$$
\begin{enumerate}
  \item Вычисление первой суммы. Используя следствие \ref{ostatok}(см. Приложение)
\begin{gather*}
\sum_{(n,k)\in\Omega_{41}}(Q+\Q+\alpha k)=\sum_{U_1<n\le \frac{R}{Q+\Q}}\Biggl(
(Q+\Q)(n-U_1)+\frac{\alpha}{2}(n^2-U_1^2)+\\+O(Q+\Q)+O(\alpha n)\Biggl)=
\int_{U_1}^{\frac{R}{Q+\Q}}\Biggl((Q+\Q)(n-U_1)+\frac{\alpha}{2}(n^2-U_1^2)\Biggl)dn+\\+
O\left(\alpha\frac{R^2}{(Q+\Q)^2}\right)+O(R)=
\frac{Q+\Q}{2}\left(\frac{R}{Q+\Q}-U_1\right)^2+\frac{\alpha}{6}\left(\frac{R^3}{(Q+\Q)^3}-U_1^3\right)-\\-
\frac{\alpha U_1^2}{2}\left(\frac{R}{Q+\Q}-U_1\right)+O\left(\alpha\left(\frac{R^2}{(Q+\Q)^2}\right)\right).
\end{gather*}
При вычислении асимптотики мы учли, что
$$Q+\Q\le U_2<\alpha n.$$
  \item Вычисление второй суммы.  Используя следствие \ref{ostatok}(см. Приложение)
\begin{gather*}
\sum_{(n,k)\in\Omega_{42}}(Q+\Q+\alpha k)=\sum_{\frac{R}{Q+\Q}<n\le \frac{R+Q^2/\alpha}{Q+\Q}}\Biggl(
(Q+\Q)\left(\frac{R-n\Q}{Q}-U_1\right)+\\+\frac{\alpha}{2}\left(\frac{(R-n\Q)^2}{Q^2}-U_1^2\right)+
O(Q+\Q)+O\left(\alpha\frac{R-n\Q}{Q}\right)
\Biggl)=\\=
\int\limits_{\frac{R}{Q+\Q}}^{\frac{R+Q^2/\alpha}{Q+\Q}}\Biggl(
(Q+\Q)\left(\frac{R-n\Q}{Q}-U_1\right)+\frac{\alpha}{2}\left(\frac{(R-n\Q)^2}{Q^2}-U_1^2\right)\Biggl)+\\+
O(R)+O\left(\frac{\alpha R^2}{(Q+\Q)^2}\right)+O\left(\frac{Q^2}{\alpha}\right)+O\left(\frac{RQ}{Q+\Q}\right).
\end{gather*}
Следовательно,
\begin{gather*}
\sum_{(n,k)\in\Omega_{42}}(Q+\Q+\alpha k)
\frac{(Q+\Q)}{2Q\Q}\left(\left(R-\Q\frac{R}{Q+\Q}\right)^2-\left(R-\Q\frac{R+Q^2/\alpha}{Q+\Q}\right)^2\right)
-U_1\frac{Q^2}{\alpha}+\\+
\frac{\alpha}{6\Q Q^2}\left(\left(R-\Q\frac{R}{Q+\Q}\right)^3-\left(R-\Q\frac{R+Q^2/\alpha}{Q+\Q}\right)^3\right)
-U_1^2\frac{Q^2}{2(Q+\Q)}+
O\left(\frac{\alpha R^2}{(Q+\Q)^2}\right)
=\\=
\frac{Q}{2\Q(Q+\Q)}\left(R^2-\left(R-\frac{Q\Q}{\alpha}\right)^2\right)-U_1\frac{Q^2}{\alpha}+\\+
\frac{\alpha Q}{6\Q(Q+\Q)^3}\left(R^3-\left(R-\frac{Q\Q}{\alpha}\right)^3\right)-U_1^2\frac{Q^2}{2(Q+\Q)}+
O\left(\frac{\alpha R^2}{(Q+\Q)^2}\right).
\end{gather*}
  \item Вычисление третьей суммы. Используя следствие \ref{ostatok}(см. Приложение)
\begin{gather*}
\sum_{(n,k)\in\Omega_{43}}(Q+\Q+\alpha k)=\sum_{U_1+\frac{Q}{\alpha}<n\le \frac{R+Q^2/\alpha}{Q+\Q}}\Biggl(
(Q+\Q)\left(n-\frac{Q}{\alpha}-U_1\right)+\\+\frac{\alpha}{2}\left(\left(n-\frac{Q}{\alpha}\right)^2-U_1^2\right)+
O(Q+\Q)+O\left(\alpha\left(n-\frac{Q}{\alpha}\right)\right)
\Biggl)=\\=
\int\limits_{U_1+\frac{Q}{\alpha}}^{ \frac{R+Q^2/\alpha}{Q+\Q}}\Biggl(
(Q+\Q)\left(n-\frac{Q}{\alpha}-U_1\right)+\frac{\alpha}{2}\left(\left(n-\frac{Q}{\alpha}\right)^2-U_1^2\right)
\Biggl)+\\+O(R)+O\left(\frac{\alpha R^2}{(Q+\Q)^2}\right)
=
\frac{Q+\Q}{2}\left(\frac{R-\frac{Q\Q}{\alpha}}{Q+\Q}-U_1\right)^2+
\frac{\alpha}{6}\left(\frac{\left(R-\frac{Q\Q}{\alpha}\right)^3}{(Q+\Q)^3}-U_1^3\right)-\\-
\frac{\alpha U_1^2}{2}\left(\frac{R-\frac{Q\Q}{\alpha}}{Q+\Q}-U_1\right)+O\left(\frac{\alpha R^2}{(Q+\Q)^2}\right).
\end{gather*}
\end{enumerate}
Следовательно,
\begin{gather*}
\left(\sum_{(n,k)\in\Omega_{41}}+\sum_{(n,k)\in\Omega_{42}}-\sum_{(n,k)\in\Omega_{43}}\right)(\Q+Q+\alpha k)=
\frac{Q+\Q}{2}\left(\frac{R^2-\left(R-\frac{Q\Q}{\alpha}\right)^2}{(Q+\Q)^2}-2U_1\frac{\frac{Q\Q}{\alpha}}{Q+\Q}\right)+\\+
\frac{\alpha}{6}\frac{R^3-\left(R-\frac{Q\Q}{\alpha}\right)^3}{(Q+\Q)^3}-
\frac{\alpha U_1^2}{2}\frac{\frac{Q\Q}{\alpha}}{Q+\Q}+\frac{Q}{2\Q(Q+\Q)}\left(R^2-
\left(R-\frac{Q\Q}{\alpha}\right)^2\right)-U_1\frac{Q^2}{\alpha}+\\+
\frac{\alpha Q}{6\Q(Q+\Q)^3}\left(R^3-\left(R-\frac{Q\Q}{\alpha}\right)^3\right)-U_1^2\frac{Q^2}{2(Q+\Q)}
+O\left(\frac{\alpha R^2}{(Q+\Q)^2}\right)
=\\=
\frac{1}{2\Q}\left(R^2-\left(R-\frac{Q\Q}{\alpha}\right)^2\right)+\\+
\frac{\alpha}{6\Q(Q+\Q)^2}\left(R^3-\left(R-\frac{Q\Q}{\alpha}\right)^3\right)-\frac{U_1}{\alpha}Q(Q+\Q)-
\frac{U_1^2}{2}Q+O\left(\frac{\alpha R^2}{(Q+\Q)^2}\right).
\end{gather*}
Полученное выражение необходимо просуммировать по
$$(Q,\Q)\in\Omega_{4}=\left\{Q\le U_2 \quad \Q\le U_2\frac{U_2-Q}{U_2+Q}\right\},\quad\tau\mid Q,\quad\tau\mid \Q.$$
\begin{enumerate}
  \item Суммируем первое слагаемое.
\begin{gather*}
\sum_{\substack{(Q,\Q)\in\Omega_{4}\\\tau\mid Q,\tau\mid \Q}}
\frac{1}{2\Q}\left(R^2-\left(R-\frac{Q\Q}{\alpha}\right)^2\right)=
\sum_{Q\le\frac{U_2}{\tau}}\sum_{\Q\le \frac{U_2}{\tau}\frac{U_2-\tau Q}{U_2+\tau Q}}
\left(\frac{R\tau}{\alpha}Q-\frac{\tau^3}{2\alpha^2}Q^2\Q\right)=\\=
\sum_{Q\le\frac{U_2}{\tau}}\Biggl(
\frac{RU_2}{\alpha}Q\frac{U_2-\tau Q}{U_2+\tau Q}-\frac{U_2^2\tau}{4\alpha^2}Q^2\frac{(U_2-\tau Q)^2}{(U_2+\tau Q)^2}
+O\left(\frac{R\tau}{\alpha}Q\right)
\Biggl).
\end{gather*}
Применяя пункты (c),(d) следствия \ref{sl} (см. Приложение), получаем
\begin{gather*}
\sum_{\substack{(Q,\Q)\in\Omega_{4}\\\tau\mid Q,\tau\mid \Q}}
\frac{1}{2\Q}\left(R^2-\left(R-\frac{Q\Q}{\alpha}\right)^2\right)=
\frac{R^{5/2}\sqrt{\alpha}}{\tau^2}\left(\log2-\frac{7}{12}\right)+O\left(\frac{R^2}{\tau}\right).
\end{gather*}
  \item Суммируем второе слагаемое.
\begin{gather*}
\sum_{\substack{(Q,\Q)\in\Omega_{4}\\\tau\mid Q,\tau\mid \Q}}
\frac{\alpha}{6\Q(Q+\Q)^2}\left(R^3-\left(R-\frac{Q\Q}{\alpha}\right)^3\right)=\\=
\sum_{Q\le\frac{U_2}{\tau}}\sum_{\Q\le \frac{U_2}{\tau}\frac{U_2-\tau Q}{U_2+\tau Q}}\frac{1}{(Q+\Q)^2}\Biggl(
\frac{R^2}{2\tau}Q-\frac{R\tau}{2\alpha}Q^2\Q+\frac{\tau^3}{6\alpha^2}Q^3\Q^2
\Biggl).
\end{gather*}
Для упрощения записи в следующих формулах $a=\frac{U_2}{\tau}$.
\begin{enumerate}
 \item[(a)]
Для первого слагаемого получаем.
\begin{gather*}
\sum_{\Q\le a\frac{a-Q}{a+Q}}\frac{1}{(Q+\Q)^2}=
\frac{a(a-Q)}{Q(Q^2+a^2)}+O\left(\frac{1}{Q^2}\right).
\end{gather*}
Используя пункт (e) следствия \ref{sl} (см. Приложение), получаем
\begin{gather*}
\sum_{Q\le\frac{U_2}{\tau}}\sum_{\Q\le \frac{U_2}{\tau}\frac{U_2-\tau Q}{U_2+\tau Q}}
\frac{R^2}{2\tau}\frac{Q}{(Q+\Q)^2}=
\frac{R^{5/2}\sqrt{\alpha}}{\tau^2}\left(\frac{\pi}{8}-\frac{\log2}{4}\right)+
O\left(\frac{R^2}{\tau}\left(1+\log\frac{U_2}{\tau}\right)\right).
\end{gather*}
 \item[(b)]
 Для второго слагаемого получаем.
\begin{gather*}
\sum_{\Q\le a\frac{a-Q}{a+Q}}\frac{\Q}{(Q+\Q)^2}=
\log\left(1+\frac{a}{Q}\frac{a-Q}{a+Q}\right)-\frac{a(a-Q)}{(Q^2+a^2)}+O\left(\frac{1}{Q}\right).
\end{gather*}
Используя пункты (f), (g) следствия \ref{sl} (см. Приложение), получаем
\begin{gather*}
\sum_{Q\le\frac{U_2}{\tau}}\sum_{\Q\le \frac{U_2}{\tau}\frac{U_2-\tau Q}{U_2+\tau Q}}
\frac{R\tau}{2\alpha}\frac{Q^2\Q}{(Q+\Q)^2}=
\frac{R^{5/2}\sqrt{\alpha}}{\tau^2}\left(\frac{\pi}{24}-\frac{5\log2}{12}+\frac{1}{6}\right)+
O\left(\frac{R^2}{\tau}\left(1+\log\frac{U_2}{\tau}\right)\right).
\end{gather*}
 \item[(c)]
 Для третьего слагаемого получаем.
\begin{gather*}
\sum_{\Q\le a\frac{a-Q}{a+Q}}\frac{\Q^2}{(Q+\Q)^2}=a\frac{a-Q}{a+Q}-
2Q\log\left(1+\frac{a}{Q}\frac{a-Q}{a+Q}\right)+Q\frac{a(a-Q)}{(Q^2+a^2)}+O(1).
\end{gather*}
Используя пункты (h), (i), (j) следствия \ref{sl} (см. Приложение), получаем
\begin{gather*}
\sum_{Q\le\frac{U_2}{\tau}}\sum_{\Q\le \frac{U_2}{\tau}\frac{U_2-\tau Q}{U_2+\tau Q}}
\frac{\tau^3}{6\alpha^2}\frac{Q^3\Q^2}{(Q+\Q)^2}=
\frac{R^{5/2}\sqrt{\alpha}}{\tau^2}\left(\frac{\pi}{120}-\frac{7\log2}{20}+\frac{13}{60}\right)+
O\left(\frac{R^2}{\tau}\left(1+\log\frac{U_2}{\tau}\right)\right).
\end{gather*}
 \end{enumerate}
Следовательно,
\begin{gather*}
\sum_{\substack{(Q,\Q)\in\Omega_{4}\\\tau\mid Q,\tau\mid \Q}}
\frac{\alpha}{6\Q(Q+\Q)^2}\left(R^3-\left(R-\frac{Q\Q}{\alpha}\right)^3\right)=\\=
\frac{R^{5/2}\sqrt{\alpha}}{\tau^2}\left(\frac{11\pi}{120}-\frac{11\log2}{60}+\frac{1}{20}\right)+
O\left(\frac{R^2}{\tau}\left(1+\log\frac{U_2}{\tau}\right)\right).
\end{gather*}
  \item Суммируем третье слагаемое.
\begin{gather*}
\sum_{\substack{(Q,\Q)\in\Omega_{4}\\\tau\mid Q,\tau\mid \Q}}
\frac{U_1}{\alpha}Q(Q+\Q)=
\sum_{Q\le\frac{U_2}{\tau}}\sum_{\Q\le \frac{U_2}{\tau}\frac{U_2-\tau Q}{U_2+\tau Q}}\frac{U_1\tau^2}{\alpha}Q(Q+\Q).
\end{gather*}
\begin{enumerate}
 \item[(a)]
Для первого слагаемого получаем.
\begin{gather*}
\sum_{\Q\le a\frac{a-Q}{a+Q}}Q^2=aQ^2\frac{a-Q}{a+Q}+O(Q^2).
\end{gather*}
Используя пункт (k) следствия \ref{sl} (см. Приложение), получаем
\begin{gather*}
\sum_{Q\le\frac{U_2}{\tau}}\sum_{\Q\le \frac{U_2}{\tau}\frac{U_2-\tau Q}{U_2+\tau Q}}\frac{U_1\tau^2}{\alpha}Q^2=
\frac{R^{5/2}\sqrt{\alpha}}{\tau^2}\left(2\log 2-\frac{4}{3}\right)+O\left(\frac{R^2}{\tau}\right).
\end{gather*}
 \item[(b)]
 Для второго слагаемого получаем.
\begin{gather*}
\sum_{\Q\le a\frac{a-Q}{a+Q}}Q\Q=\frac{a^2}{2}Q\frac{(a-Q)^2}{(a+Q)^2}+ O\left(aQ\frac{a-Q}{a+Q}\right).
\end{gather*}
Используя пункты (c), (l), следствия \ref{sl} (см. Приложение), получаем
\begin{gather*}
\sum_{Q\le\frac{U_2}{\tau}}\sum_{\Q\le \frac{U_2}{\tau}\frac{U_2-\tau Q}{U_2+\tau Q}}\frac{U_1\tau^2}{\alpha}Q\Q=
\frac{R^{5/2}\sqrt{\alpha}}{\tau^2}\left(4\log 2-\frac{11}{4}\right)+O\left(\frac{R^2}{\tau}\right).
\end{gather*}
Следовательно,
\begin{gather*}
\sum_{\substack{(Q,\Q)\in\Omega_{4}\\\tau\mid Q,\tau\mid \Q}}\frac{U_1}{\alpha}Q(Q+\Q)=
\frac{R^{5/2}\sqrt{\alpha}}{\tau^2}\left(6\log 2-\frac{49}{12}\right)+O\left(\frac{R^2}{\tau}\right).
\end{gather*}
 \end{enumerate}
  \item Суммируем четвертое слагаемое. Используя пункт (c), следствия \ref{sl} (см. Приложение), получаем
\begin{gather*}
\sum_{\substack{(Q,\Q)\in\Omega_{4}\\\tau\mid Q,\tau\mid \Q}}
\frac{U_1^2}{2}Q=
\sum_{Q\le\frac{U_2}{\tau}}\sum_{\Q\le \frac{U_2}{\tau}\frac{U_2-\tau Q}{U_2+\tau Q}}\frac{U_1^2}{2}\tau Q =
\frac{R^{5/2}\sqrt{\alpha}}{\tau^2}\left(\frac{3}{4}-\log 2\right)+O\left(\frac{R^2}{\tau}\right).
\end{gather*}
\end{enumerate}
Следовательно,
\begin{gather*}
\sum_{\substack{(Q,\Q)\in\Omega_{4}\\\tau\mid Q,\tau\mid \Q}} \left(\sum_{(n,k)\in\Omega_{41}}+\sum_{(n,k)\in\Omega_{42}}-\sum_{(n,k)\in\Omega_{43}}\right)(\Q+Q+\alpha k)=\\=
\frac{R^{5/2}\sqrt{\alpha}}{\tau^2}\left(\frac{11}{120}\pi-\frac{251}{60}\log2+\frac{14}{5}\right)+
O\left(\frac{R^2}{\tau}\left(1+\log\frac{U_2}{\tau}\right)\right)+O\left(\frac{R^2\alpha}{\tau^2}\log\frac{U_2}{\tau}\right).
\end{gather*}
Тем самым лемма доказана.
\end{proof}
\end{Le}
\subsection{Случай 5}
\begin{Le}
\label{51}
Справедлива следующая асимптотическая формула
\begin{gather*}
\Sigma_{51}=\left(\frac{\pi}{24}+\frac{55}{12}\log 2-\frac{119}{36}\right)R^{5/2}\sqrt{\alpha}
\sum_{\tau|\delta}\frac{\varphi(\tau)}{\tau^2}+\sum_{\tau|\delta}\varphi(\tau)
O\left(\frac{R^2}{\tau}\left(1+\log\frac{U_2}{\tau}\right)+\frac{\alpha R^2}{\tau^2}\right),
\end{gather*}
где $\Sigma_{51}$ определена в \eqref{slu5}.
\begin{proof}
Вычислим отдельно суммы по областям
$$
\Omega_{51}=\left\{U_1<n\le\frac{R}{Q+\Q},\quad U_1<k\le n \right\},
$$
$$
\Omega_{52}=\left\{\frac{R}{Q+\Q}<n\le\frac{R-U_1Q}{\Q},\quad U_1<k\le\frac{R-n\Q}{Q}\right\}.
$$
\begin{enumerate}
  \item Вычисление первой суммы аналогично первому случаю в лемме \ref{41}
\begin{gather*}
\sum_{(n,k)\in\Omega_{51}}(Q+\Q+\alpha k)=\sum_{U_1<n\le \frac{R}{Q+\Q}}\Biggl(
(Q+\Q)(n-U_1)+\frac{\alpha}{2}(n^2-U_1^2)+\\+O(Q+\Q+\alpha n)\Biggl)=
\frac{Q+\Q}{2}\left(\frac{R}{Q+\Q}-U_1\right)^2+\frac{\alpha}{6}\left(\frac{R^3}{(Q+\Q)^3}-U_1^3\right)-\\-
\frac{\alpha U_1^2}{2}\left(\frac{R}{Q+\Q}-U_1\right)+O\left(\alpha\left(\frac{R^2}{(Q+\Q)^2}\right)\right).
\end{gather*}
При вычислении асимптотики мы учли, что
$$Q+\Q\le U_2<\alpha n.$$
  \item Вычисление второй суммы. Используя следствие \ref{ostatok}(см. Приложение)
\begin{gather*}
\sum_{(n,k)\in\Omega_{52}}(Q+\Q+\alpha k)=\sum_{\frac{R}{Q+\Q}<n\le \frac{R-U_1Q}{\Q}}\Biggl(
(Q+\Q)\left(\frac{R-n\Q}{Q}-U_1\right)+\\+\frac{\alpha}{2}\left(\frac{(R-n\Q)^2}{Q^2}-U_1^2\right)+O(Q+\Q)+
O\left(\alpha\frac{R-n\Q}{Q}\right)\Biggl)=\\=
\int\limits_{\frac{R}{Q+\Q}}^{\frac{R-U_1Q}{\Q}}\Biggl(
(Q+\Q)\left(\frac{R-n\Q}{Q}-U_1\right)+\frac{\alpha}{2}\left(\frac{(R-n\Q)^2}{Q^2}-U_1^2\right)\Biggl)+
O\left(\alpha\frac{R^2}{(Q+\Q)^2}\right)+\\+O(R)+O\left((Q+\Q)\left(\frac{R-U_1Q}{\Q}-\frac{R}{Q+\Q}\right)\right)
+O\left(\frac{\alpha}{2}\frac{Q}{\Q}\left(\frac{R^2}{(Q+\Q)^2}-U_1^2\right)\right).
\end{gather*}
Следовательно,
\begin{gather*}
\sum_{(n,k)\in\Omega_{52}}(Q+\Q+\alpha k)=
(Q+\Q)\frac{Q}{2\Q}\left(\frac{R}{Q+\Q}-U_1\right)^2+\frac{\alpha}{6}\frac{Q}{\Q}\left(\frac{R^3}{(Q+\Q)^3}-U_1^3\right)-\\-
\frac{\alpha}{2}U_1^2\left(\frac{R-U_1Q}{\Q}-\frac{R}{Q+\Q}\right)+
O\left((Q+\Q)\left(\frac{R-U_1Q}{\Q}-\frac{R}{Q+\Q}\right)\right)+\\+O\left(\frac{\alpha}{2}\frac{Q}{\Q}\left(\frac{R^2}{(Q+\Q)^2}-U_1^2\right)\right).
\end{gather*}
Используя $(Q,\Q)\in\Omega_{5}$, легко получить
\begin{gather*}
\Q>U_2\frac{U_2-Q}{U_2+Q}=\frac{R-U_1Q}{U_1+Q/\alpha}.
\end{gather*}
Следовательно,
\begin{equation}\label{nerav5}
\frac{R-U_1Q}{\Q}<\frac{R+Q^2/\alpha}{Q+\Q}.
\end{equation}
Поэтому
\begin{gather*}
O\left((Q+\Q)\left(\frac{R-U_1Q}{\Q}-\frac{R}{Q+\Q}\right)\right)<O\left(\frac{Q^2}{\alpha}\right)<
O\left(\alpha\frac{R^2}{(Q+\Q)^2}\right),
\end{gather*}
Из неравенства \eqref{nerav5} следует
\begin{equation}
R-\frac{Q\Q}{\alpha}<U_1(Q+\Q).
\end{equation}
Поэтому
\begin{gather*}
O\left(\frac{\alpha}{2}\frac{Q}{\Q}\left(\frac{R^2}{(Q+\Q)^2}-U_1^2\right)\right)<
O\left(\frac{Q\alpha}{2\Q(Q+\Q)^2}\left(R^2-\left(R-\frac{Q\Q}{\alpha}\right)^2\right)\right)<\\<
O\left(\alpha\frac{R^2}{(Q+\Q)^2}\right).
\end{gather*}
\end{enumerate}
Следовательно,
\begin{gather*}
\left(\sum_{(n,k)\in\Omega_{51}}+\sum_{(n,k)\in\Omega_{52}}\right)(\Q+Q+\alpha k)=
\frac{1}{2\Q}(Q+\Q)^2\left(\frac{R}{Q+\Q}-U_1\right)^2+\\+
\frac{\alpha}{6}\frac{Q+\Q}{\Q}\left(\frac{R^3}{(Q+\Q)^3}-U_1^3\right)-
\frac{\alpha}{2}U_1^2\left(\frac{R-U_1Q}{\Q}-U_1\right)
+O\left(\frac{\alpha R^2}{(Q+\Q)^2}\right).
\end{gather*}
Полученное выражение необходимо просуммировать по
$$(Q,\Q)\in\Omega_{5}=\left\{Q\le U_2 \quad U_2\frac{U_2-Q}{U_2+Q}<\Q\le U_2-Q\right\},\quad\tau\mid Q,\quad\tau\mid \Q.$$
Для упрощения записи в формулах часто используется обозначение $a=\frac{U_2}{\tau}$.
Определим функцию $G(Q,\Q)$
\begin{gather*}
G(Q,\Q)=\left(-\frac{2}{3}RU_1Q+\frac{U_1^2\tau}{2}Q^2\right)\frac{1}{\Q}+\\+\frac{\alpha R^3}{6\tau^3\Q(Q+\Q)^2}+
\left(U_1^2\tau Q-\frac{2}{3}RU_1\right)+\frac{U_1^2\tau}{2}\Q.
\end{gather*}
Используя соотношение
\begin{gather*}
\frac{1}{2\Q}(Q+\Q)^2\left(\frac{R}{Q+\Q}-U_1\right)^2=
\left(\frac{R^2}{2}-RU_1Q+\frac{U_1^2}{2}Q^2\right)\frac{1}{\Q}-RU_1+U_1^2Q+\frac{U_1^2}{2}\Q,
\end{gather*}
после тривиальных преобразований легко получить
\begin{gather*}
\sum_{\substack{(Q,\Q)\in\Omega_{5}\\\tau\mid Q,\tau\mid \Q}}
\left(\sum_{(n,k)\in\Omega_{51}}+\sum_{(n,k)\in\Omega_{52}}\right)(\Q+Q+\alpha k)=\\=
\sum_{Q\le a}\sum_{a\frac{a-Q}{a+Q}<\Q\le a-Q}\Biggl(
G(Q,\Q)+O\left(\frac{\alpha R^2}{\tau^2(Q+\Q)^2}\right)
\Biggl).
\end{gather*}
Заметим
$$2(\sqrt{2}-1)U_2<Q+\Q\le U_2,$$
по-этому, используя те же рассуждения, что и в случаи 3, получаем
\begin{gather*}
\sum_{\substack{(Q,\Q)\in\Omega_{5}\\\tau\mid Q,\tau\mid \Q}}
\left(\sum_{(n,k)\in\Omega_{51}}+\sum_{(n,k)\in\Omega_{52}}\right)(\Q+Q+\alpha k)=\\=\sum_{Q\le a}
\left(\int\limits_{a\frac{a-Q}{a+Q}}^{a-Q}G(Q,\Q)d\Q+G\left(Q,a\frac{a-Q}{a+Q}\right)\right)+
O\left(\frac{\alpha R^2}{\tau^2}\right).
\end{gather*}
Используя следующие оценки
\begin{gather*}
\sum_{Q\le a}Q\frac{a+Q}{a-Q}=O(a^2\log a),\quad
\sum_{Q\le a}Q^2\frac{a+Q}{a-Q}=O(a^3\log a),\\
\sum_{Q\le a}\frac{(a+Q)^3}{(a-Q)(a^2+Q^2)^2}=O\left(\frac{\log a}{a}\right),
\end{gather*}
легко получить
$$
\sum_{Q\le a}G\left(Q,a\frac{a-Q}{a+Q}\right)=O\left(\frac{R^2}{\tau}\left(1+\log\frac{U_2}{\tau}\right)\right).
$$
Следовательно,
\begin{gather*}
\sum_{\substack{(Q,\Q)\in\Omega_{5}\\\tau\mid Q,\tau\mid \Q}}
\left(\sum_{(n,k)\in\Omega_{51}}+\sum_{(n,k)\in\Omega_{52}}\right)(\Q+Q+\alpha k)=\\=\sum_{Q\le a}\Biggl(
\left(-\frac{2}{3}RU_1Q+\frac{U_1^2\tau}{2}Q^2\right)\log\left(1+\frac{Q}{a}\right)+
\frac{\alpha R^3}{6\tau^3}\left(\frac{1}{Q^2}\log\left(1+\frac{Q^2}{a^2}\right)-\frac{1}{a}\frac{a-Q}{a^2+Q^2}\right)+\\+
\left(U_1^2\tau Q-\frac{2}{3}RU_1\right)Q\frac{a-Q}{a+Q}+\frac{U_1^2\tau}{4}(a-Q)^2\frac{(a+Q)^2-a^2}{(a+Q)^2}
\Biggl)+\\+
O\left(\frac{\alpha R^2}{\tau^2}\right)+O\left(\frac{R^2}{\tau}\left(1+\log\frac{U_2}{\tau}\right)\right).
\end{gather*}
\begin{enumerate}
  \item Суммируем первое слагаемое. Используя пункты (m), (n), следствия \ref{sl} (см. Приложение), получаем
\begin{gather*}
\sum_{Q\le a}\left(-\frac{2}{3}RU_1Q+\frac{U_1^2\tau}{2}Q^2\right)\log\left(1+\frac{Q}{a}\right)=
\frac{R^{5/2}\sqrt{\alpha}}{\tau^2}\left(\frac{1}{3}\log2-\frac{11}{36}\right)+O\left(\frac{R^2}{\tau}\right).
\end{gather*}
  \item Суммируем второе слагаемое.Используя пункты (o), (e) следствия \ref{sl} (см. Приложение), получаем
\begin{gather*}
\sum_{Q\le a}
\frac{\alpha R^3}{6\tau^3}\left(\frac{1}{Q^2}\log\left(1+\frac{Q^2}{a^2}\right)-\frac{1}{a}\frac{a-Q}{a^2+Q^2}\right)=
\frac{R^{5/2}\sqrt{\alpha}}{\tau^2}\left(\frac{\pi}{24}-\frac{1}{12}\log2\right)+O\left(\frac{R^2}{\tau}\right).
\end{gather*}
  \item Суммируем третье слагаемое. Используя пункты (c), (k) следствия \ref{sl} (см. Приложение), получаем
\begin{gather*}
\sum_{Q\le a}
\left(U_1^2\tau Q-\frac{2}{3}RU_1\right)Q\frac{a-Q}{a+Q}=
\frac{R^{5/2}\sqrt{\alpha}}{\tau^2}\left(\frac{10}{3}\log2-\frac{7}{3}\right)+O\left(\frac{R^2}{\tau}\right).
\end{gather*}
   \item Суммируем четвертое слагаемое. Используя пункты (p) следствия \ref{sl} (см. Приложение), получаем
\begin{gather*}
\sum_{Q\le a}\frac{U_1^2\tau}{4}(a-Q)^2\frac{(a+Q)^2-a^2}{(a+Q)^2}=
\frac{R^{5/2}\sqrt{\alpha}}{\tau^2}\left(\log2-\frac{2}{3}\right)+O\left(\frac{R^2}{\tau}\right).
\end{gather*}
\end{enumerate}
Следовательно,
\begin{gather*}
\sum_{\substack{(Q,\Q)\in\Omega_{5}\\\tau\mid Q,\tau\mid \Q}}
\left(\sum_{(n,k)\in\Omega_{51}}+\sum_{(n,k)\in\Omega_{52}}\right)(\Q+Q+\alpha k)=\\=
\frac{R^{5/2}\sqrt{\alpha}}{\tau^2}\left(\frac{\pi}{24}+\frac{55}{12}\log2-\frac{119}{36}\right)+
O\left(\frac{\alpha R^2}{\tau^2}\right)+O\left(\frac{R^2}{\tau}\left(1+\log\frac{U_2}{\tau}\right)\right).
\end{gather*}
Тем самым лемма доказана.
\end{proof}
\end{Le}
\section{Вычисление сумм второго типа}
В этом параграфе мы вычислим $\Sigma_{12}, \Sigma_{22}, \Sigma_{32}, \Sigma_{42}, \Sigma_{52}.$
\subsection{Случай 1}
\begin{Le}
\label{12}
Справедлива следующая асимптотическая формула
\begin{gather*}
\Sigma_{12}=\sum_{\tau|\delta}\mathop{{\sum}^*}_{z=1}^{\tau}\sum_{\substack{(n,k)\in\Omega_1\\ \tau|n, \tau\nmid k}}
\sum_{(Q,\Q)\in\Omega_{11}}\exp\left(2\pi i\frac{zk}{\tau}Q\right)(\Q+Q+\alpha k)=\\=
\sum_{\tau|\delta}\varphi(\tau)
O\left(\frac{R^2}{\tau}\left(1+\log \frac{U_1}{\tau}\log\tau\right)\right),
\end{gather*}
где $\Sigma_{12}$ определена в \eqref{slu1}.
\begin{proof}
Суммирование по $\Q$ было проведено в лемме \ref{11}. Используя лемму \ref{trigsum} и следствие \ref{trigsum12} получаем
\begin{gather*}
\sum_{(Q,\Q)\in\Omega_{11}}\exp\left(2\pi i\frac{zk}{\tau}Q\right)(\Q+Q+\alpha k)\ll \frac{1}{\|\frac{zk}{\tau}\|}
\Biggl(
\left|\frac{k^2}{2n^2}-\frac{k}{n}\right|\alpha^2 n^2+\\+
\left|\frac{R}{n}-\frac{\alpha k^2}{n}-R\frac{k}{n^2}\right|\alpha n+
\left|\alpha R\frac{k}{n}+\frac{R^2}{2n^2}\right|
\Biggl)+\\+
\sum_{\alpha (n-k)<Q\le \alpha n}O\left(Q+\alpha k\right)+
\sum_{\alpha (n-k)<Q\le \alpha n}O\left(\frac{R-kQ}{n}\right).
\end{gather*}
Применяя оценки из лемме \ref{11}, получаем
\begin{gather*}
\sum_{(Q,\Q)\in\Omega_{11}}\exp\left(2\pi i\frac{zk}{\tau}Q\right)(\Q+Q+\alpha k)\ll \frac{1}{\|\frac{zk}{\tau}\|}
\Biggl(\alpha^2kn+R\alpha+\frac{R^2}{n^2}
\Biggl)+O\left(R\alpha\right)\ll\\\ll
\frac{1}{\|\frac{zk}{\tau}\|}\frac{R^2}{n^2}+O\left(R\alpha\right).
\end{gather*}
Полученное выражение необходимо просуммировать по переменным $$(n,k)\in\Omega_1, \tau|n, \tau\nmid k,$$ используя лемму \ref{dozelogo}
\begin{gather*}
\sum_{\substack{(n,k)\in\Omega_1\\ \tau|n, \tau\nmid k}}
\left(\frac{1}{\|\frac{zk}{\tau}\|}\frac{R^2}{n^2}+O\left(R\alpha\right)\right)\ll
\sum_{n\le\frac{U_1}{\tau}}\frac{R^2}{\tau^2n^2}\tau(n+1)\log\tau+O\left(\frac{R^2}{\tau}\right)=\\=
O\left(\frac{R^2}{\tau}\left(1+\log \frac{U_1}{\tau}\log\tau\right)\right).
\end{gather*}
Тем самым лемма доказана.
\end{proof}
\end{Le}

\subsection{Случай 2}
\begin{Le}
\label{22}
Справедлива следующая асимптотическая формула
\begin{gather*}
\Sigma_{22}=\sum_{\tau|\delta}\varphi(\tau)
O\left(\frac{R^2}{\tau}\left(1+\log \frac{U_1}{\tau}\right)+\frac{R^2\alpha}{\tau}\log\tau+
R^{3/2}\sqrt{\alpha}\log\tau\right),
\end{gather*}
где $\Sigma_{22}$ определена в \eqref{slu2}.
\begin{proof}
Вычислим отдельно суммы по областям $\Omega_{21}$ и $\Omega_{22}$.
\begin{enumerate}
  \item
Используя вычисления произведенные в лемме \ref{21} в соответствующем пункте, а так же  лемму \ref{trigsum} и следствие \ref{trigsum12} получаем
\begin{gather*}
\sum\limits_{(n,Q)\in\Omega_{21}}\exp\left(2\pi i\frac{z\Q}{\tau}n\right)(\Q+Q+\alpha k)=\\=
\sum_{U_1<n\le\frac{R}{\Q+\alpha k}}\exp\left(2\pi i\frac{z\Q}{\tau}n\right)
\left(\alpha^2nk+\alpha^2\frac{k^2}{2}+\alpha k\Q+O\left(\alpha n\right)+O\left(\Q+\alpha k\right)\right)\ll\\\ll
\frac{1}{\|\frac{z\Q}{\tau}\|}\left(\frac{R\alpha^2k}{\Q+\alpha k}+\alpha k\Q+\frac{\alpha^2k^2}{2}\right)+
O\left(\alpha\frac{R^2}{(\Q+\alpha k)^2}\right)\ll\\\ll
\frac{1}{\|\frac{z\Q}{\tau}\|}\frac{R\alpha^2k}{\Q+\alpha k}+O\left(\alpha\frac{R^2}{(\Q+\alpha k)^2}\right).
\end{gather*}
  \item
Используя вычисления произведенные в лемме \ref{21} в соответствующем пункте, а так же  лемму \ref{trigsum} и следствие \ref{trigsum12} получаем
\begin{gather*}
\sum\limits_{(n,Q)\in\Omega_{22}}\exp\left(2\pi i\frac{z\Q}{\tau}n\right)(\Q+Q+\alpha k)=\\=
\sum_{U_1<n\le\frac{R}{\Q+\alpha k}}\exp\left(2\pi i\frac{z\Q}{\tau}n\right)f(n)+
O\left(\frac{R\alpha^2k^2}{(\Q+\alpha k)^2}\right),
\end{gather*}
где функция $f(n)$ определена в лемме \ref{21}.\par
Применим следствие \ref{trigsum12} с учетом того, что
$$
f\left(\frac{R+\alpha k^2}{\Q+\alpha k}\right)=0.
$$
Заметим, что удобнее использовать следующее преобразование
\begin{gather*}
f\left(\frac{R}{\Q+\alpha k}\right)=f\left(\frac{R}{\Q+\alpha k}\right)-f\left(\frac{R+\alpha k^2}{\Q+\alpha k}\right)
=\frac{R\alpha^2k}{\Q+\alpha k}+\Q\alpha k+\frac{\alpha^2k^2}{2}.
\end{gather*}
Тогда легко получить следующую оценку
\begin{gather*}
\sum\limits_{(n,Q)\in\Omega_{22}}\exp\left(2\pi i\frac{z\Q}{\tau}n\right)(\Q+Q+\alpha k)\ll\\\ll
\frac{1}{\|\frac{z\Q}{\tau}\|}\left(\frac{R\alpha^2k}{\Q+\alpha k}+\Q\alpha k+\frac{\alpha^2k^2}{2}\right)+
O\left(\frac{R\alpha^2k^2}{(\Q+\alpha k)^2}\right).
\end{gather*}
\end{enumerate}
Следовательно,
\begin{gather*}
\left(\sum\limits_{(n,Q)\in\Omega_{21}}+\sum\limits_{(n,Q)\in\Omega_{22}}\right)
\exp\left(2\pi i\frac{z\Q}{\tau}n\right)(\Q+Q+\alpha k)\ll\\\ll
\frac{1}{\|\frac{z\Q}{\tau}\|}\frac{R\alpha^2k}{\Q+\alpha k}+
O\left(\frac{\alpha R^2}{(\Q+\alpha k)^2}\right).
\end{gather*}
Полученное выражение необходимо просуммировать по
$$(k,\Q)\in\Omega_{2}=\left\{k\le U_1 \quad \Q\le U_2-\alpha k\right\},\quad\tau\mid k,\quad\tau\nmid \Q.$$
\begin{enumerate}
  \item
Применяя лемму \ref{abel} и лемму \ref{trigsum}, получаем
\begin{gather*}
\sum_{\substack{\Q\le U_2-\alpha k\\\tau\nmid \Q }}\frac{1}{\|\frac{z\Q}{\tau}\|}\frac{1}{\Q+\alpha k}\ll
\left(\frac{U_2-\alpha k}{U_2}+\log\frac{U_2}{\alpha k}+\frac{\tau}{\alpha k}\right)\log\tau.
\end{gather*}
Далее, используя следствия \ref{sl}, получаем
\begin{gather*}
\sum_{\substack{k\le U_1\\\tau\mid k}}R\alpha^2k
\left(\frac{U_1-k}{U_1}+\log\frac{U_1}{k}+\frac{\tau}{\alpha k}\right)\log\tau\ll\\\ll
\frac{R^2\alpha}{\tau}\log\tau+\frac{R^2\alpha}{\tau}\log\tau+R^{3/2}\sqrt{\alpha}\log\tau.
\end{gather*}
  \item
Применяя лемму \ref{euler}, получаем
\begin{gather*}
\sum_{\substack{(k,\Q)\in\Omega_{2}\\\tau\mid k,\tau\nmid \Q}}O\left(\frac{\alpha R^2}{(\Q+\alpha k)^2}\right)=
O\left(\frac{R^2}{\tau}\left(1+\log \frac{U_1}{\tau}\right)\right).
\end{gather*}
\end{enumerate}
Тем самым лемма доказана.
\end{proof}
\end{Le}

\subsection{Случай 3}
\begin{Le}
\label{32}
Справедлива следующая асимптотическая формула
\begin{gather*}
\Sigma_{32}=\sum_{\tau|\delta}\varphi(\tau)
O\left(\frac{R^2}{\tau}+\frac{R^2}{\tau}(1+\alpha)\log\tau\right),
\end{gather*}
где $\Sigma_{32}$ определена в \eqref{slu3}.
\begin{proof}
Используя вычисления произведенные в лемме \ref{31}, получаем
\begin{gather*}
\sum\limits_{(n,Q)\in\Omega_{31}}\exp\left(2\pi i\frac{z\Q}{\tau}n\right)(\Q+Q+\alpha k)=\\=
\sum_{U_1<n\le\frac{R+\alpha k^2}{\Q+\alpha k}}\exp\left(2\pi i\frac{z\Q}{\tau}n\right)f(n)+
O\left(\frac{R^2}{k(\Q+\alpha k)}\right),
\end{gather*}
где функция $f(n)$ определена в лемме \ref{21}.\par
Применим следствие \ref{trigsum12} с учетом того, что
\begin{gather*}
f\left(\frac{R+\alpha k^2}{\Q+\alpha k}\right)=0,\\
f(U_1)=\Q^2\left(\frac{U_1^2}{2k^2}-\frac{U_1}{k}\right)-
\Q\left(U_1\alpha+\frac{R+\alpha k^2}{k^2}U_1-\frac{R+\alpha k^2}{k}\right)+\frac{(R+\alpha k^2)^2}{k^2}-
\frac{U_1^2\alpha^2}{2},
\end{gather*}
получаем
\begin{gather*}
\sum_{U_1<n\le\frac{R+\alpha k^2}{\Q+\alpha k}}\exp\left(2\pi i\frac{z\Q}{\tau}n\right)f(n)\ll
\frac{1}{\|\frac{z\Q}{\tau}\|}f(U_1)\ll\\\ll
\frac{1}{\|\frac{z\Q}{\tau}\|}\left(\Q^2\frac{U_1^2}{k^2}+\Q\frac{R+\alpha k^2}{k^2}U_1+\frac{R^2}{k^2}\right).
\end{gather*}
Полученное выражение необходимо просуммировать по
$$
(k,\Q)\in\Omega_{3}=\left\{k\le U_1, \quad U_2-\alpha k<\Q\le U_2-\alpha k+\frac{\alpha k^2}{U_1}\right\},
\quad\tau\mid k,\quad\tau\nmid \Q.
$$
Оценим каждое слагаемое, используя лемму\ref{dozelogo}.
\begin{enumerate}
  \item Суммируем первое слагаемое.
\begin{gather*}
\sum_{\substack{U_2-\alpha k<\Q\le U_2-\alpha k+\frac{\alpha k^2}{U_1}\\\tau\nmid \Q}}
\frac{\Q^2}{\|\frac{z\Q}{\tau}\|}\frac{U_1^2}{k^2}\ll
\frac{U_1^2}{k^2}\left( U_2-\alpha k+\frac{\alpha k^2}{U_1}\right)^2\left(\frac{\alpha k^2}{U_1}+\tau\right)\log\tau\ll\\\ll
\frac{U_1^4\alpha^2}{k^2}\left(\frac{\alpha k^2}{U_1}+\tau\right)\log\tau\ll\\
\left(U_1^3\alpha^2+\frac{U_1^4\alpha^2\tau}{k^2}\right)\log\tau.
\end{gather*}
Следовательно,
\begin{gather*}
\sum_{\substack{(k,\Q)\in\Omega_{3}\\\tau\mid k,\tau\nmid \Q}}\frac{\Q^2}{\|\frac{z\Q}{\tau}\|}\frac{U_1^2}{k^2}\ll
\frac{R^2}{\tau}\log\tau.
\end{gather*}
  \item Суммируем второе слагаемое.
\begin{gather*}
\sum_{\substack{U_2-\alpha k<\Q\le U_2-\alpha k+\frac{\alpha k^2}{U_1}\\\tau\nmid \Q}}
\frac{\Q}{\|\frac{z\Q}{\tau}\|}\frac{R+\alpha k^2}{k^2}U_1\ll
\frac{R+\alpha k^2}{k^2}U_1^2\alpha\left(\frac{\alpha k^2}{U_1}+\tau\right)\log\tau\ll\\\ll
\left(RU_1\alpha^2+\frac{R^2}{k^2}\tau\right)\log\tau.
\end{gather*}
Следовательно,
\begin{gather*}
\sum_{\substack{(k,\Q)\in\Omega_{3}\\\tau\mid k,\tau\nmid \Q}}
\frac{\Q}{\|\frac{z\Q}{\tau}\|}\frac{R+\alpha k^2}{k^2}U_1\ll
\frac{R^2}{\tau}(1+\alpha)\log\tau.
\end{gather*}
  \item Суммируем третье слагаемое.
\begin{gather*}
\sum_{\substack{(k,\Q)\in\Omega_{3}\\\tau\mid k,\tau\nmid \Q}}
\frac{1}{\|\frac{z\Q}{\tau}\|}\frac{R^2}{k^2}\ll
\sum_{\substack{k\le U_1\\\tau\mid k}}\frac{R^2}{k^2}\left(\frac{\alpha k^2}{U_1}+\tau\right)\log\tau\ll\\\ll
\frac{R^2}{\tau}(1+\alpha)\log\tau.
\end{gather*}
  \item Суммируем асимптотическое слагаемое.
\begin{gather*}
\sum_{\substack{(k,\Q)\in\Omega_{3}\\\tau\mid k,\tau\nmid \Q}}
O\left(\frac{R^2}{k(\Q+\alpha k)}\right)=\sum_{\substack{k\le U_1\\\tau\mid k}}
O\left(\frac{R^2}{k}\log\left(1+\frac{k^2}{U_1^2}\right)\right)\ll\\\ll
\sum_{\substack{k\le U_1\\\tau\mid k}}O\left(\frac{R^2k}{U_1^2}\right)=O\left(\frac{R^2}{\tau}\right).
\end{gather*}
\end{enumerate}
Тем самым лемма доказана.
\end{proof}
\end{Le}

\subsection{Случай 4}
\begin{Le}
\label{42}
Справедлива следующая асимптотическая формула
\begin{gather*}
\Sigma_{42}=\sum_{\tau|\delta}\varphi(\tau)
O\left(\frac{R^2\alpha}{\tau}(1+\log\tau)\left(1+\log\frac{U_2}{\tau}\right)\right),
\end{gather*}
где $\Sigma_{42}$ определена в \eqref{slu4}.
\begin{proof}
Вычислим каждую из трех внутренних сумм.
\begin{enumerate}
  \item Вычисление первой суммы.Определим функцию
$$g_1(n)=\frac{\alpha}{2}n^2+(Q+\Q)n-U_1(Q+\Q)-\frac{\alpha}{2}U_1^2,$$
тогда
\begin{gather*}
\sum_{U_1<k\le n}(Q+\Q+\alpha k)=g_1(n)+O(Q+\Q)+O(\alpha n).
\end{gather*}
Заметим, что  $g_1(n)$ возрастает на промежутке $U_1<n\le \frac{R}{Q+\Q}$. Очевидно, что
$$g_1(U_1)=0,\quad g_1\left(\frac{R}{Q+\Q}\right)\ll\frac{R^2\alpha}{(Q+\Q)^2}.$$
Применим следствие \ref{trigsum12}, получаем
\begin{gather*}
\sum_{(n,k)\in\Omega_{41}}\exp\left(2\pi i\frac{z\Q}{\tau}n\right)(\Q+Q+\alpha k)\ll
\frac{1}{\|\frac{z\Q}{\tau}\|}\frac{R^2\alpha}{(Q+\Q)^2}.
\end{gather*}
  \item Вычисление второй суммы.
\begin{gather*}
\sum_{U_1<k\le \frac{R-n\Q}{Q}}(Q+\Q+\alpha k)=
(Q+\Q)\left(\frac{R-n\Q}{Q}-U_1\right)+\frac{\alpha}{2}\left(\frac{(R-n\Q)^2}{Q^2}-U_1^2\right)+\\+O(Q+\Q)+Q\left(\alpha \frac{R-n\Q}{Q}\right)
=\frac{\alpha\Q^2}{2Q^2}n^2-
\left(Q+\Q+\frac{R\alpha}{Q}\right)\frac{\Q}{Q}n+(Q+\Q)\left(\frac{R}{Q}-U_1\right)+\\+
\frac{\alpha}{2}\left(\frac{R^2}{Q^2}-U_1^2\right)+O(Q+\Q)+O\left(\alpha \frac{R-n\Q}{Q}\right)=
g_2(n)+O(Q+\Q)+O\left(\alpha \frac{R-n\Q}{Q}\right).
\end{gather*}
Заметим, что  $g_2(n)$ убывает на промежутке
$$\frac{R}{Q+\Q}<n\le\frac{R+\frac{Q^2}{\alpha}}{Q+\Q}.$$
Действительно,
\begin{gather*}
g_2(n)’=\frac{\alpha\Q^2}{Q^2}n-\left(Q+\Q+\frac{R\alpha}{Q}\right)\frac{\Q}{Q}=
\frac{\alpha\Q}{Q^2}(n\Q-R)-\frac{\Q}{Q}(Q+\Q)<0.
\end{gather*}
Мы воспользовались тем, что $n\Q\le R$ (см. \eqref{sistem}).
Очевидно, что
$$g_2\left(\frac{R}{Q+\Q}\right)\ll\frac{R^2\alpha}{(Q+\Q)^2}.$$
Применим следствие \ref{trigsum12}, получаем
\begin{gather*}
\sum_{(n,k)\in\Omega_{42}}\exp\left(2\pi i\frac{z\Q}{\tau}n\right)(\Q+Q+\alpha k)\ll
\frac{1}{\|\frac{z\Q}{\tau}\|}\frac{R^2\alpha}{(Q+\Q)^2}.
\end{gather*}
  \item Вычисление третьей суммы.
\begin{gather*}
\sum_{U_1<k\le n-\frac{Q}{\alpha}}(Q+\Q+\alpha k)=\frac{\alpha}{2}n^2+n\Q-(Q+\Q)\left(\frac{Q}{\alpha}+U_1\right)+\\+
\frac{\alpha}{2}\left(\frac{Q^2}{\alpha^2}-U_1^2\right)+O(Q+\Q)+O(\alpha n-Q)=
g_3(n)+O(Q+\Q)+O(\alpha n-Q).
\end{gather*}
Заметим, что  $g_3(n)$ возрастает  на промежутке $U_1+\frac{Q}{\alpha}<\frac{R+\frac{Q^2}{\alpha}}{Q+\Q}$. Очевидно, что
$$
g_3\left(U_1+\frac{Q}{\alpha}\right)=0,\quad
g_3\left(\frac{R+\frac{Q^2}{\alpha}}{Q+\Q}\right)\ll\frac{R^2\alpha}{(Q+\Q)^2}.$$
Применим следствие \ref{trigsum12}, получаем
\begin{gather*}
\sum_{(n,k)\in\Omega_{43}}\exp\left(2\pi i\frac{z\Q}{\tau}n\right)(\Q+Q+\alpha k)\ll
\frac{1}{\|\frac{z\Q}{\tau}\|}\frac{R^2\alpha}{(Q+\Q)^2}.
\end{gather*}
\end{enumerate}
Следовательно,
\begin{gather*}
\left(\sum_{(n,k)\in\Omega_{41}}+\sum_{(n,k)\in\Omega_{42}}-\sum_{(n,k)\in\Omega_{43}}\right)
\exp\left(2\pi i\frac{z\Q}{\tau}n\right)(\Q+Q+\alpha k)\ll\\\ll
\frac{1}{\|\frac{z\Q}{\tau}\|}\frac{R^2\alpha}{(Q+\Q)^2}+O\left(\frac{R^2\alpha}{(Q+\Q)^2}\right)+O\left(R\right).
\end{gather*}

Полученное выражение необходимо просуммировать по
$$(Q,\Q)\in\Omega_{4}=\left\{Q\le U_2 \quad \Q\le U_2\frac{U_2-Q}{U_2+Q}\right\},\quad\tau\mid Q,\quad\tau\nmid \Q.$$
Увеличим область суммирования по $\Q$ до
$$\Q\le U_2-Q.$$
Используя лемму \ref{abel} и лемму \ref{dozelogo}, получаем
\begin{gather*}
\sum_{\substack{\Q\le U_2-Q\\\tau\nmid \Q}}\frac{1}{\|\frac{z\Q}{\tau}\|}\frac{R^2\alpha}{(Q+\Q)^2}\ll
R^2\alpha\left(\frac{U_2-Q}{U_2^2}+\frac{(U_2-Q)^2}{QU_2^2}+\frac{\tau}{Q^2}\right)\log\tau.
\end{gather*}
Легко получить
\begin{gather*}
R^2\alpha\sum_{\substack{Q\le U_2\\\tau\mid \Q}}
\left(\frac{U_2-Q}{U_2^2}+\frac{(U_2-Q)^2}{QU_2^2}+\frac{\tau}{Q^2}\right)\log\tau=
O\left(\frac{R^2\alpha}{\tau}\left(1+\log\frac{U_2}{\tau}\right)\log\tau\right).
\end{gather*}
Осталось просуммировать асимптотические слагаемые.
\begin{gather*}
\sum_{\substack{Q\le U_2,\Q\le U_2-Q\\\tau\mid \Q}}
\left(O\left(\frac{R^2\alpha}{(Q+\Q)^2}\right)+O\left(R\right)\right)=
O\left(\frac{R^2\alpha}{\tau}\left(1+\log\frac{U_2}{\tau}\right)\right).
\end{gather*}
Тем самым лемма доказана.
\end{proof}
\end{Le}

\subsection{Случай 5}
\begin{Le}
\label{52}
Справедлива следующая асимптотическая формула
\begin{gather*}
\Sigma_{52}=\sum_{\tau|\delta}\varphi(\tau)
O\left(\frac{R^2\alpha}{\tau}\right),
\end{gather*}
где $\Sigma_{52}$ определена в \eqref{slu5}.
\begin{proof}
Вычислим каждую из трех внутренних сумм.
\begin{enumerate}
  \item
Заметим, что первая сумма совпадает с первой суммой леммы \ref{42}. Так как в \mbox{лемме \ref{42}} мы увеличивали область суммирования по $\Q$ до
$$\Q\le U_2-Q,$$
то тем самым мы уже все сосчитали.
   \item Вычисление второй суммы.
\begin{gather*}
\sum_{U_1<k\le \frac{R-n\Q}{Q}}(Q+\Q+\alpha k)=(Q+\Q)\left(\frac{R-n\Q}{Q}-U_1\right)+
\frac{\alpha}{2}\left(\frac{(R-n\Q)^2}{Q^2}-U_1^2\right)+\\
+O(Q+\Q)+Q\left(\alpha \frac{R-n\Q}{Q}\right)=
g_2(n)+O(Q+\Q)+Q\left(\alpha \frac{R-n\Q}{Q}\right).
\end{gather*}
Заметим, что  $g_2(n)$ убывает на промежутке
$$\frac{R}{Q+\Q}<n\le\frac{R-U_1Q}{\Q}$$ и
$$g_2\left(\frac{R}{Q+\Q}\right)\ll\frac{R^2\alpha}{(Q+\Q)^2}.$$
Применим следствие \ref{trigsum12}, получаем
\begin{gather*}
\sum_{(n,k)\in\Omega_{52}}\exp\left(2\pi i\frac{z\Q}{\tau}n\right)(\Q+Q+\alpha k)\ll
\frac{1}{\|\frac{z\Q}{\tau}\|}\frac{R^2\alpha}{(Q+\Q)^2}.
\end{gather*}
Мы получили ту же оценку что и в лемме \ref{42}, а так как в лемме \ref{42} мы увеличивали область суммирования по $\Q$ до
$$\Q\le U_2-Q,$$
то тем самым мы уже сосчитали суммы, не входящие в асимптотические слагаемые. Суммируем  асимптотические слагаемые по n
\begin{gather*}
\sum_{\frac{R}{Q+\Q}<n\le\frac{R-U_1Q}{\Q}}\left(O(Q+\Q)+O\left(\alpha \frac{R-n\Q}{Q}\right)\right)\ll
O\left(R\frac{Q}{\Q}\right)+O\left(\alpha\frac{R^2Q}{\Q(Q+\Q)^2}\right).
\end{gather*}
\end{enumerate}
Полученное выражение необходимо просуммировать по
$$(Q,\Q)\in\Omega_{5}=\left\{Q\le U_2, \quad U_2\frac{U_2-Q}{U_2+Q}<\Q\le U_2-Q\right\},\quad\tau\mid Q,\quad\tau\nmid \Q.$$
\begin{enumerate}
  \item Суммируем первое слагаемое.
  \begin{gather*}
\sum_{U_2\frac{U_2-Q}{U_2+Q}<\Q\le U_2-Q}O\left(R\frac{Q}{\Q}\right)=
O\left(RQ\log\left(1+\frac{Q}{U_2}\right)\right)\ll O\left(R\frac{Q^2}{U_2}\right).
\end{gather*}
Следовательно,
\begin{gather*}
\sum_{\substack{(Q,\Q)\in\Omega_{5}\\\tau\mid Q,\tau\nmid \Q}}O\left(R\frac{Q}{\Q}\right)\ll
O\left(\frac{R^2\alpha}{\tau}\right).
\end{gather*}
  \item Суммируем второе слагаемое.
\begin{gather*}
\sum_{U_2\frac{U_2-Q}{U_2+Q}<\Q\le U_2-Q}\frac{1}{\Q(Q+\Q)^2}\ll
\int\limits_{\frac{Q^2+U_2^2}{Q+U_2}}^{U_2}\frac{1}{(x-Q)x^2}dx=\\=
\int\limits_{\frac{Q^2+U_2^2}{Q+U_2}}^{U_2}\left(-\frac{1}{Q^2x}-\frac{1}{Qx^2}+\frac{1}{Q^2(x-Q)}\right)dx=
\frac{1}{Q^2}\log\left(1+\frac{Q^2}{U_2^2}\right)-\frac{U_2-Q}{U_2(U_2^2+Q^2)}.
\end{gather*}
Следовательно,
\begin{gather*}
\sum_{U_2\frac{U_2-Q}{U_2+Q}<\Q\le U_2-Q}O\left(\alpha\frac{R^2Q}{\Q(Q+\Q)^2}\right)\ll
O\left(\alpha\frac{R^2Q}{U_2^2}\right)+O\left(\alpha R^2\frac{U_2-Q}{U_2^2+Q^2}\right).
\end{gather*}
Используя пункт (e) следствия \ref{sl}, получаем
\begin{gather*}
\sum_{\substack{(Q,\Q)\in\Omega_{5}\\\tau\mid Q,\tau\nmid \Q}}O\left(\alpha\frac{R^2Q}{\Q(Q+\Q)^2}\right)\ll
O\left(\frac{R^2\alpha}{\tau}\right).
\end{gather*}
\end{enumerate}
Тем самым лемма доказана.
\end{proof}
\end{Le}
\section{Вычисление сумм третьего типа}
В этом параграфе мы вычислим $\Sigma_{13}, \Sigma_{23}, \Sigma_{33}, \Sigma_{43}, \Sigma_{53}.$
\subsection{Случай 1}
\begin{Le}
\label{13}
Справедлива следующая асимптотическая формула
\begin{gather*}
\Sigma_{13}=O\left(R^2\sum_{d\mid\delta}d\log d\right),
\end{gather*}
где $\Sigma_{13}$ определена в \eqref{slu1}.
\begin{proof}
Оценим внутреннюю сумму, используя следствие \ref{trigsum12}, получаем
\begin{gather*}
\sum_{\Q\le \frac{R-kQ}{n}}\exp\left(2\pi i\frac{zn}{\tau}\Q\right)(\Q+Q+\alpha k)\ll \frac{1}{\|\frac{zn}{\delta}\|}
\left(\frac{R-kQ}{n}+Q+\alpha k\right)\ll\frac{1}{\|\frac{zn}{\delta}\|}\frac{R}{n}.
\end{gather*}
Следовательно,
\begin{gather*}
\sum_{(Q,\Q)\in\Omega_{11}}\exp\left(2\pi i\frac{zn}{\tau}\Q\right)(\Q+Q+\alpha k)\ll
\frac{R}{n}\alpha k\frac{1}{\|\frac{zn}{\delta}\|}.
\end{gather*}
и
\begin{gather*}
\sum_{k\le n}\frac{R}{n}\alpha k\frac{1}{\|\frac{zn}{\delta}\|}\ll R\alpha\frac{n}{\|\frac{zn}{\delta}\|}.
\end{gather*}
Отсюда получаем
\begin{gather*}
\Sigma_{13}\ll\sum_{z=1}^{\delta}\sum_{\substack{n\le U_1\\\delta\nmid zn}}R\alpha\frac{n}{\|\frac{zn}{\delta}\|}=
\sum_{n\le U_1}R\alpha n\sum_{\substack{z\le\delta\\\delta\nmid zn}}\frac{1}{\|\frac{zn}{\delta}\|}=
\sum_{d\mid\delta}\sum_{\substack{n\le U_1/d\\(n,\delta/d)=1}}R\alpha nd
\sum_{\substack{z\le\delta\\\frac{\delta}{d}\nmid zn}}\frac{1}{\|\frac{zn}{\delta/d}\|}.
\end{gather*}
Так как
\begin{gather*}
\sum_{\substack{z\le\delta\\\frac{\delta}{d}\nmid zn}}\frac{1}{\|\frac{zn}{\delta/d}\|}\ll\delta\log\frac{\delta}{d},
\end{gather*}
то получаем
\begin{gather*}
\Sigma_{13}\ll\sum_{d\mid\delta}d\delta\log\frac{\delta}{d}\sum_{n\le\frac{U_1}{d}}R\alpha n\ll
R^2\sum_{d\mid\delta}d\log d.
\end{gather*}
Тем самым лемма доказана.
\end{proof}
\end{Le}

\subsection{Случай 2}
\begin{Le}
\label{23}
Справедлива следующая асимптотическая формула
\begin{gather*}
\Sigma_{23}=O\left(R^2\sum_{d\mid\delta}d\log d\left(1+\log\frac{U_1}{\delta/d}\right)\right),
\end{gather*}
где $\Sigma_{23}$ определена в \eqref{slu2}.
\begin{proof}
Вычислим отдельно суммы по областям $\Omega_{21}$ и $\Omega_{22}$.
\begin{enumerate}
  \item
Используя следствие \ref{trigsum12}, получаем
\begin{gather*}
\sum_{\alpha(n-k)<Q\le\alpha n}\exp\left(2\pi i\frac{zk}{\delta}Q\right)(\Q+Q+\alpha k)\ll \frac{1}{\|\frac{zk}{\delta}\|}(\Q+\alpha n+\alpha n)\ll
\frac{1}{\|\frac{zk}{\delta}\|}\alpha n.
\end{gather*}
Следовательно,
\begin{gather*}
\sum\limits_{(n,Q)\in\Omega_{21}}\exp\left(2\pi i\frac{zk}{\delta}Q\right)(\Q+Q+\alpha k)\ll
\frac{R^2\alpha}{(\Q+\alpha k)^2}\frac{1}{\|\frac{zk}{\delta}\|}.
\end{gather*}
Суммируем по $\Q$, получаем
\begin{gather*}
\sum_{\Q\le U_2-\alpha k}\sum\limits_{(n,Q)\in\Omega_{21}}\exp\left(2\pi i\frac{zk}{\delta}Q\right)(\Q+Q+\alpha k)\ll
\frac{R^2}{k}\frac{1}{\|\frac{zk}{\delta}\|}.
\end{gather*}
  \item
Используя следствие \ref{trigsum12}, получаем
\begin{gather*}
\sum_{\alpha(n-k)<Q\le\frac{R-n\Q}{k}}\exp\left(2\pi i\frac{zk}{\delta}Q\right)(\Q+Q+\alpha k)\ll \frac{1}{\|\frac{zk}{\delta}\|}\left(\frac{R-n\Q}{k}+\Q+\alpha n\right)\ll
\frac{1}{\|\frac{zk}{\delta}\|}\frac{R}{k}.
\end{gather*}
Следовательно,
\begin{gather*}
\sum\limits_{(n,Q)\in\Omega_{22}}\exp\left(2\pi i\frac{zk}{\delta}Q\right)(\Q+Q+\alpha k)\ll
\frac{R\alpha k}{\Q+\alpha k}\frac{1}{\|\frac{zk}{\delta}\|}.
\end{gather*}
Суммируем по $\Q$, получаем
\begin{gather*}
\sum_{\Q\le U_2-\alpha k}\sum\limits_{(n,Q)\in\Omega_{22}}\exp\left(2\pi i\frac{zk}{\delta}Q\right)(\Q+Q+\alpha k)\ll
R\alpha k\log\frac{U_1}{k}\frac{1}{\|\frac{zk}{\delta}\|}.
\end{gather*}
\end{enumerate}
Следовательно,
\begin{gather*}
\Sigma_{23}\ll\sum_{z=1}^{\delta}\sum_{\substack{k\le U_1\\\delta\nmid zk}}
\left(\frac{R^2}{k}+R\alpha k\log\frac{U_1}{k}\right)\frac{1}{\|\frac{zn}{\delta}\|}\ll
\sum_{d\mid\delta}\delta\log\frac{\delta}{d}\sum_{k\le\frac{U_1}{d}}
\left(\frac{R^2}{dk}+R\alpha dk\log\frac{U_1}{dk}\right).
\end{gather*}
Используя пункт (b), следствия \ref{sl}, получаем
\begin{gather*}
\Sigma_{23}\ll\sum_{d\mid\delta}\delta\log\frac{\delta}{d}
\left(\frac{R^2}{d}\log\frac{U_1}{d}+\frac{R^2}{d}\right).
\end{gather*}
Тем самым лемма доказана.
\end{proof}
\end{Le}

\subsection{Случай 3}
\begin{Le}
\label{33}
Справедлива следующая асимптотическая формула
\begin{gather*}
\Sigma_{33}=O\left(R^2\sum_{d\mid\delta}d\log d\right),
\end{gather*}
где $\Sigma_{33}$ определена в \eqref{slu3}.
\begin{proof}
Используя следствие \ref{trigsum12}, получаем
\begin{gather*}
\sum_{\alpha(n-k)<Q\le\frac{R-n\Q}{k}}\exp\left(2\pi i\frac{zk}{\delta}Q\right)(\Q+Q+\alpha k)\ll \frac{1}{\|\frac{zk}{\delta}\|}\left(\frac{R-n\Q}{k}+\Q+\alpha k\right)\ll
\frac{1}{\|\frac{zk}{\delta}\|}\frac{R}{k}.
\end{gather*}
Следовательно,
\begin{gather*}
\sum\limits_{(n,Q)\in\Omega_{31}}\exp\left(2\pi i\frac{zk}{\delta}Q\right)(\Q+Q+\alpha k)\ll
\frac{R(R+\alpha k^2)}{k(\Q+\alpha k)}\frac{1}{\|\frac{zk}{\delta}\|}.
\end{gather*}
Суммируем по $\Q$, получаем
\begin{gather*}
\sum_{U_2-\alpha k<\Q\le U_2-\alpha k+\frac{\alpha k^2}{U_1}}\sum\limits_{(n,Q)\in\Omega_{31}}
\exp\left(2\pi i\frac{zk}{\delta}Q\right)(\Q+Q+\alpha k)\ll\\\ll
\frac{R}{k}(R+\alpha k^2)\log\left(1+\frac{k^2}{U_1^2}\right)\frac{1}{\|\frac{zk}{\delta}\|}\ll
\frac{R^2}{k}\log\left(1+\frac{k^2}{U_1^2}\right)\frac{1}{\|\frac{zk}{\delta}\|}.
\end{gather*}
Следовательно,
\begin{gather*}
\Sigma_{33}\ll\sum_{z=1}^{\delta}\sum_{\substack{k\le U_1\\\delta\nmid zk}}
\frac{R^2}{k}\log\left(1+\frac{k^2}{U_1^2}\right)\frac{1}{\|\frac{zk}{\delta}\|}\ll
\sum_{d\mid\delta}\delta\log\frac{\delta}{d}\sum_{k\le\frac{U_1}{d}}\frac{R^2}{dk}\log\left(1+\frac{d^2k^2}{U_1^2}\right)\ll\\\ll
\sum_{d\mid\delta}\delta\log\frac{\delta}{d}\sum_{k\le\frac{U_1}{d}}\frac{R^2}{U_1^2}dk\ll R^2\sum_{d\mid\delta}d\log d.
\end{gather*}
Тем самым лемма доказана.
\end{proof}
\end{Le}

\subsection{Случай 4}
\begin{Le}
\label{43}
Справедлива следующая асимптотическая формула
\begin{gather*}
\Sigma_{43}=O\left(R^2\alpha\sum_{d\mid\delta}d\log d\left(1+\log\frac{U_2}{\delta/d}\right)\right),
\end{gather*}
где $\Sigma_{43}$ определена в \eqref{slu4}.
\begin{proof}
Вычислим отдельно суммы по трем областям.
\begin{enumerate}
  \item
Используя следствие \ref{trigsum12}, получаем
\begin{gather*}
\sum_{(n,k)\in\Omega_{41}}\exp\left(2\pi i\frac{zQ}{\delta}k\right)(\Q+Q+\alpha k)\ll
\frac{1}{\|\frac{zQ}{\delta}\|}\sum_{U_1<n\le\frac{R}{Q+\Q}}(\Q+Q+\alpha n)\ll\\\ll
\frac{1}{\|\frac{zQ}{\delta}\|}\sum_{U_1<n\le\frac{R}{Q+\Q}}\alpha n\ll
\frac{1}{\|\frac{zQ}{\delta}\|}\frac{R^2\alpha}{(Q+\Q)^2}.
\end{gather*}
Суммируем по $\Q$, получаем
\begin{gather*}
\sum_{\Q \le U_2\frac{U_2-Q}{U_2+Q}}
\frac{1}{\|\frac{zQ}{\delta}\|}\frac{R^2\alpha}{(Q+\Q)^2}\ll
\frac{1}{\|\frac{zQ}{\delta}\|}\frac{R^2\alpha}{Q}.
\end{gather*}
  \item
Используя следствие \ref{trigsum12}, получаем
\begin{gather*}
\sum_{(n,k)\in\Omega_{42}}\exp\left(2\pi i\frac{zQ}{\delta}k\right)(\Q+Q+\alpha k)\ll
\frac{1}{\|\frac{zQ}{\delta}\|}\sum_{\frac{R}{Q+\Q}<n\le\frac{R+Q^2/\alpha}{Q+\Q}}
\left(\Q+Q+\alpha\frac{R-n\Q}{Q}\right)\ll\\\ll
\frac{1}{\|\frac{zQ}{\delta}\|}\sum_{\frac{R}{Q+\Q}<n\le\frac{R+Q^2/\alpha}{Q+\Q}}\frac{R\alpha}{Q}\ll
\frac{1}{\|\frac{zQ}{\delta}\|}\frac{RQ}{Q+\Q}.
\end{gather*}
Суммируем по $\Q$, получаем
\begin{gather*}
\sum_{\Q\le U_2\frac{U_2-Q}{U_2+Q}}\frac{1}{\|\frac{zQ}{\delta}\|}\frac{RQ}{Q+\Q}\ll
\frac{1}{\|\frac{zQ}{\delta}\|}RQ\log\left(\frac{U_2^2+Q^2}{Q(U_2+Q)}\right).
\end{gather*}
  \item
Используя следствие \ref{trigsum12}, получаем
\begin{gather*}
\sum_{(n,k)\in\Omega_{43}}\exp\left(2\pi i\frac{zQ}{\delta}k\right)(\Q+Q+\alpha k)\ll
\frac{1}{\|\frac{zQ}{\delta}\|}\sum_{U_1+\frac{Q}{\alpha}<n\le\frac{R+Q^2/\alpha}{Q+\Q}}\alpha n\ll\\\ll
\frac{1}{\|\frac{zQ}{\delta}\|}\alpha\frac{(R+Q^2/\alpha)^2}{(Q+\Q)^2}.
\end{gather*}
Суммируем по $\Q$, получаем
\begin{gather*}
\sum_{\Q\le U_2\frac{U_2-Q}{U_2+Q}}\frac{1}{\|\frac{zQ}{\delta}\|}\alpha\frac{(R+Q^2/\alpha)^2}{(Q+\Q)^2}\ll
\frac{1}{\|\frac{zQ}{\delta}\|}\alpha\frac{(R+Q^2/\alpha)^2}{Q}\ll
\frac{1}{\|\frac{zQ}{\delta}\|}\frac{R^2\alpha}{Q}.
\end{gather*}
\end{enumerate}
Следовательно,
\begin{gather*}
\Sigma_{43}\ll\sum_{z=1}^{\delta}\sum_{\substack{Q\le U_2\\\delta\nmid zQ}}\frac{1}{\|\frac{zQ}{\delta}\|}
\left(\frac{R^2\alpha}{Q}+RQ\log\left(\frac{U_2^2+Q^2}{Q(U_2+Q)}\right)+\frac{R^2\alpha}{Q}\right)\ll\\\ll
\sum_{d\mid\delta}\delta\log\frac{\delta}{d}\sum_{Q\le\frac{U_2}{d}}
\left(\frac{R^2\alpha}{dQ}+RdQ\log\left(\frac{U_2^2+d^2Q^2}{dQ(U_2+dQ)}\right)\right).
\end{gather*}
Используя пункт (r), следствия \ref{sl}, получаем
\begin{gather*}
\Sigma_{43}\ll
\sum_{d\mid\delta}\frac{\delta}{d}\log\frac{\delta}{d}\left(R^2\alpha\log\frac{U_2}{d}+R^2\alpha\right).
\end{gather*}
Тем самым лемма доказана.
\end{proof}
\end{Le}

\subsection{Случай 5}
\begin{Le}
\label{53}
Справедлива следующая асимптотическая формула
\begin{gather*}
\Sigma_{53}=O\left(R^2\alpha\sum_{d\mid\delta}d\log d\right),
\end{gather*}
где $\Sigma_{53}$ определена в \eqref{slu5}.
\begin{proof}
Вычислим отдельно суммы по двум областям.
\begin{enumerate}
  \item
Используя следствие \ref{trigsum12}, получаем
\begin{gather*}
\sum_{(n,k)\in\Omega_{51}}\exp\left(2\pi i\frac{zQ}{\delta}k\right)(\Q+Q+\alpha k)\ll
\frac{1}{\|\frac{zQ}{\delta}\|}\sum_{U_1<n\le\frac{R}{Q+\Q}}(\Q+Q+\alpha n)\ll\\\ll
\frac{1}{\|\frac{zQ}{\delta}\|}\sum_{U_1<n\le\frac{R}{Q+\Q}}\alpha n\ll
\frac{1}{\|\frac{zQ}{\delta}\|}\frac{R^2\alpha}{(Q+\Q)^2}.
\end{gather*}
Суммируем по $\Q$, получаем
\begin{gather*}
\sum_{U_2\frac{U_2-Q}{U_2+Q}<\Q\le U_2-Q}
\frac{1}{\|\frac{zQ}{\delta}\|}\frac{R^2\alpha}{(Q+\Q)^2}\ll
\frac{1}{\|\frac{zQ}{\delta}\|}R^2\alpha\frac{U_2+Q}{U_2^2+Q^2}.
\end{gather*}
  \item
Используя следствие \ref{trigsum12}, получаем
\begin{gather*}
\sum_{(n,k)\in\Omega_{52}}\exp\left(2\pi i\frac{zQ}{\delta}k\right)(\Q+Q+\alpha k)\ll
\frac{1}{\|\frac{zQ}{\delta}\|}\sum_{\frac{R}{Q+\Q}<n\le\frac{R-U_1Q}{\Q}}(\Q+Q+\alpha\frac{R-n\Q}{Q})\ll\\\ll
\frac{1}{\|\frac{zQ}{\delta}\|}\sum_{\frac{R}{Q+\Q}<n\le\frac{R-U_1Q}{\Q}}\frac{R\alpha}{Q}\ll
\frac{1}{\|\frac{zQ}{\delta}\|}\frac{R^2\alpha}{Q}\left(\frac{1}{\Q}-\frac{1}{Q+\Q}\right).
\end{gather*}
Суммируем по $\Q$, получаем
\begin{gather*}
\sum_{U_2\frac{U_2-Q}{U_2+Q}<\Q\le U_2-Q}\frac{R^2\alpha}{Q}\left(\frac{1}{\Q}-\frac{1}{Q+\Q}\right)
\frac{1}{\|\frac{zQ}{\delta}\|}\ll
\frac{1}{\|\frac{zQ}{\delta}\|}\frac{R^2\alpha}{Q}\log\left(1+\frac{Q^2}{U_2^2}\right)\ll\\\ll
\frac{1}{\|\frac{zQ}{\delta}\|}\frac{R^2Q\alpha}{U_2^2}.
\end{gather*}
\end{enumerate}
Следовательно, используя пункт (s), следствия \ref{sl}, получаем
\begin{gather*}
\Sigma_{53}\ll\sum_{z=1}^{\delta}\sum_{\substack{Q\le U_2\\\delta\nmid zQ}}\frac{1}{\|\frac{zQ}{\delta}\|}
\left(R^2\alpha\frac{U_2+Q}{U_2^2+Q^2}+\frac{R^2Q\alpha}{U_2^2}\right)\ll\\\ll
\sum_{d\mid\delta}\delta\log\frac{\delta}{d}\sum_{Q\le\frac{U_2}{d}}
\left(R^2\alpha\frac{U_2+dQ}{U_2^2+d^2Q^2}+\frac{R^2dQ\alpha}{U_2^2}\right)\ll
R^2\alpha\sum_{d\mid\delta}\frac{\delta}{d}\log\frac{\delta}{d}.
\end{gather*}
Тем самым лемма доказана.
\end{proof}
\end{Le}

\section{Доказательство основного результата}
\begin{Th}\label{Th2}
Справедлива оценка
\begin{gather*}
F-G=
O\left(N^{4+\varepsilon}\left(
x_1^{1+\varepsilon}\frac{x_3^{\varepsilon}}{x_2^{\varepsilon}}+
x_2^{1+\varepsilon}\frac{x_3^{\varepsilon}}{x_1^{\varepsilon}}
+x_3^{1+\varepsilon}
\right)\right),
\end{gather*}
где $x_1>0,x_2>0, x_3>0, \varepsilon>0$~---действительные числа. F и G определены в \eqref{FG}.
\begin{proof}
Подставляя в \eqref{summa} полученные в леммах \ref{11}---\ref{53} асимптотические формулы и используя равенство
\begin{gather*}
\frac{53}{150}+\frac{19}{75}+\left(\frac{2\pi}{15}-\frac{2\log 2}{5}-\frac{91}{900}\right)+\\\nonumber+
\left(\frac{11\pi}{120}-\frac{251}{60}\log 2+\frac{14}{5}\right)+
\left(\frac{\pi}{24}+\frac{55}{12}\log 2-\frac{119}{36}\right)=\frac{4\pi}{15},
\end{gather*}
получаем
\begin{gather}\nonumber
\delta\sum\limits_{a\le R \atop \delta|a}\lambda(a,\alpha)=\sum_{j=1}^3\sum_{i=1}^5\Sigma_{ij}
=\frac{4\pi}{15}R^{5/2}\sqrt{\alpha}\sum_{\tau|\delta}\frac{\varphi(\tau)}{\tau^2}+\\\nonumber+
\sum_{\tau|\delta}\varphi(\tau)\Biggl(
\frac{R^2}{\tau}(1+\log\tau)\left(1+\log\frac{U_1}{\tau}+\log\frac{U_2}{\tau}\right)(1+\alpha)+
R^{3/2}\sqrt{\alpha}\log\tau
\Biggl)+\\+
O\left(R^2(1+\alpha)\sum_{d\mid\delta}d\log d\right)+
O\left(R^2\sum_{d\mid\delta}d\log d\left(\log\frac{U_1}{\delta/d}+\alpha\log\frac{U_2}{\delta/d}\right)\right).
\end{gather}
Для оценки полученных выражений воспользуемся тривиальной оценкой функции Эйлера, логарифма и леммой \ref{tau}.
\begin{enumerate}
  \item Оценим первое асимптотическое слагаемое. Так как
\begin{gather*}
\frac{\varphi(\tau)}{\tau}(1+\log\tau)\left(1+\log\frac{U_1}{\tau}+\log\frac{U_2}{\tau}\right)\ll
(1+\tau^{\varepsilon})\left(1+R^{\varepsilon}\tau^{\varepsilon}(\alpha^{\varepsilon}+\alpha^{-\varepsilon})\right)\ll\\\ll
R^{\varepsilon}\tau^{\varepsilon}(\alpha^{\varepsilon}+\alpha^{-\varepsilon}),
\end{gather*}
\begin{gather*}
R^{\varepsilon}\tau^{\varepsilon}(\alpha^{\varepsilon}+\alpha^{-\varepsilon})(1+\alpha)\ll
R^{\varepsilon}\tau^{\varepsilon}(\alpha^{1+\varepsilon}+\alpha^{-\varepsilon})
\end{gather*}
и
\begin{gather*}
\varphi(\tau)\sqrt{\alpha}\log\tau\ll\sqrt{\alpha}\tau^{1+\varepsilon}
\end{gather*}
Следовательно,
\begin{gather*}
\sum_{\tau|\delta}\varphi(\tau)\Biggl(
\frac{R^2}{\tau}(1+\log\tau)\left(1+\log\frac{U_1}{\tau}+\log\frac{U_2}{\tau}\right)(1+\alpha)+
R^{3/2}\sqrt{\alpha}\log\tau
\Biggl)\ll\\\ll
R^{2+\varepsilon}\delta^{\varepsilon}(\alpha^{1+\varepsilon}+\alpha^{-\varepsilon})+R^{3/2}\sqrt{\alpha}\delta^{1+\varepsilon}
\end{gather*}
  \item Оценим второе асимптотическое слагаемое.
\begin{gather*}
R^2(1+\alpha)\sum_{d\mid\delta}d\log d\ll R^2(1+\alpha)\delta^{1+\varepsilon}
\end{gather*}
  \item Оценим третье асимптотическое слагаемое. Так как
\begin{gather*}
R^2\sum_{d\mid\delta}d\log d\left(\log\frac{U_1}{\delta/d}+\alpha\log\frac{U_2}{\delta/d}\right)\ll
R^{2+\varepsilon}\delta^{1+\varepsilon}(\alpha^{1+\varepsilon}+\alpha^{-\varepsilon})
\end{gather*}
\end{enumerate}
Следовательно,
\begin{gather}\label{otvetlj}
\delta\sum\limits_{a\le R \atop \delta|a}\lambda(a,\alpha)=
\frac{4\pi}{15}R^{5/2}\sqrt{\alpha}\sum_{\tau|\delta}\frac{\varphi(\tau)}{\tau^2}+
O\left(R^{2+\varepsilon}\delta^{1+\varepsilon}(\alpha^{1+\varepsilon}+\alpha^{-\varepsilon})\right).
\end{gather}

Для дальнейших рассуждений нам потребуется следующая лемма.
\begin{Le}
\label{ljsigma}
\begin{gather*}
\sum\limits_{a\le R \atop \delta|a}\sigma_{-1}(a)a^{3/2}=
\frac{\pi^2}{15}\frac{R^{5/2}}{\delta}\sum_{\tau|\delta}\frac{\varphi(\tau)}{\tau^2}+
O\left(R^{3/2+\varepsilon}\delta^{\varepsilon}\right),
\end{gather*}
где
$$\sigma_{-1}(a)=\sum_{d\mid a}d^{-1}.$$
\begin{proof}
\begin{gather*}
\sum\limits_{a\le R \atop \delta|a}\sigma_{-1}(a)a^{3/2}=
\sum_{d\le R}\frac{1}{d}\sum_{\substack{a\le R/d\\\delta\mid ad}}a^{3/2}d^{3/2}=
\sum_{\tau|\delta}\sum_{\substack{d\le R\\(d,\delta)=\tau}}\sqrt{d}
\sum_{\substack{a\le R/d\\\frac{\delta}{\tau}\mid a}}a^{3/2}=\\=
\sum_{\tau|\delta}\sum_{\substack{d\le R\\(d,\delta)=\tau}}\sqrt{d}
\left(\frac{2R^{5/2}\tau}{5d^{5/2}\delta}+O\left(\frac{R^{3/2}}{d^{3/2}}\right)\right)=
\frac{2}{5}\frac{R^{5/2}}{\delta}\sum_{\tau|\delta}
\sum_{\substack{d\le R/\tau\\(d,\frac{\delta}{\tau})=1}}\frac{1}{\tau d^2}+
O\left(R^{3/2}\sum_{\tau|\delta}
\sum_{\substack{d\le R/\tau\\(d,\frac{\delta}{\tau})=1}}\frac{1}{\tau d}\right)
\end{gather*}
Совершая замену переменной $\delta_1=\frac{\delta}{\tau}$ и после этого обозначая $\tau=\delta_1$, получаем
\begin{gather*}
\sum\limits_{a\le R \atop \delta|a}\sigma_{-1}(a)a^{3/2}=
\frac{2}{5}\frac{R^{5/2}}{\delta^2}\sum_{\tau|\delta}\tau
\sum_{\substack{d\le R\tau/\delta\\(d,\tau)=1}}\frac{1}{d^2}+
O\left(\frac{R^{3/2}}{\delta}\sum_{\tau|\delta}\tau
\sum_{\substack{d\le R\tau/\delta\\(d,\tau)=1}}\frac{1}{d}\right)
\end{gather*}
Преобразуем обе внутренние суммы.
\begin{enumerate}
  \item Преобразуем первое слагаемое.
\begin{gather*}
\sum_{\substack{d\le R\tau/\delta\\(d,\tau)=1}}\frac{1}{d^2}=
\sum_{d\le R\tau/\delta}\frac{1}{d^2}\sum_{l\mid(d,\tau)}\mu(l)=
\sum_{l|\tau}\frac{\mu(l)}{l^2}\sum_{d\le\frac{R\tau}{l\delta}}\frac{1}{d^2}=\\=
\zeta(2)\sum_{l|\tau}\frac{\mu(l)}{l^2}+\sum_{l|\tau}O\left(\frac{\delta}{R\tau l}\right)=
\zeta(2)\sum_{l|\tau}\frac{\mu(l)}{l^2}+O\left(\frac{\delta\sigma_{-1}(\tau)}{R\tau}\right).
\end{gather*}
  \item Для оценки остаточного члена используем очевидное соотношение
\begin{gather*}
\sum_{\substack{d\le R\tau/\delta\\(d,\tau)=1}}\frac{1}{d}=
O\left(\log\frac{R}{\delta/\tau}\right).
\end{gather*}
Следовательно,
\begin{gather*}
O\left(\frac{R^{3/2}}{\delta}\sum_{\tau|\delta}\tau
\sum_{\substack{d\le R\tau/\delta\\(d,\tau)=1}}\frac{1}{d}\right)=
O\left(R^{3/2}\sum_{\tau|\delta}\frac{\log\frac{R}{\tau}}{\tau}\right)=
O\left(R^{3/2+\varepsilon}\delta^{\varepsilon}\right).
\end{gather*}
\end{enumerate}
Следовательно,
\begin{gather*}
\sum\limits_{a\le R \atop \delta|a}\sigma_{-1}(a)a^{3/2}=
\frac{2\zeta(2)}{5}\frac{R^{5/2}}{\delta^2}\sum_{\tau|\delta}\tau\sum_{l|\tau}\frac{\mu(l)}{l^2}+
O\left(\frac{R^{3/2}}{\delta}\sum_{\tau|\delta}\sigma_{-1}(\tau)\right)
+O\left(R^{3/2+\varepsilon}\delta^{\varepsilon}\right).
\end{gather*}
Используя следующее соотношение
\begin{gather*}
\sum_{\tau|\delta}\tau\sum_{l|\tau}\frac{\mu(l)}{l^2}=
\sum_{l|\delta}\frac{\mu(l)}{l^2}\sum_{\tau|\frac{\delta}{l}}l\tau=
\sum_{\tau|\delta}\tau\sum_{l|\frac{\delta}{\tau}}\frac{\mu(l)}{l}=
\sum_{\tau|\delta}\frac{\tau^2}{\delta}\varphi\left(\frac{\delta}{\tau}\right)=
\delta\sum_{\tau|\delta}\frac{\varphi(\tau)}{\tau^2}.
\end{gather*}
получаем
\begin{gather*}
\sum\limits_{a\le R \atop \delta|a}\sigma_{-1}(a)a^{3/2}=
\frac{\pi^2}{15}\frac{R^{5/2}}{\delta}\sum_{\tau|\delta}\frac{\varphi(\tau)}{\tau^2}+
O\left(R^{3/2+\varepsilon}\delta^{\varepsilon}\right).
\end{gather*}
Тем самым лемма доказана.
\end{proof}
\end{Le}
Подставляя результат леммы \ref{ljsigma} в формулу (\ref{otvetlj}), получаем
\begin{gather*}
\sum\limits_{a\le R \atop \delta|a}\lambda(a,\alpha)=
\sum\limits_{a\le R \atop \delta|a}\frac{4}{\pi}\sigma_{-1}(a)a^{3/2}\sqrt{\alpha}+
O\left(R^{2+\varepsilon}\delta^{\varepsilon}(\alpha^{1+\varepsilon}+\alpha^{-\varepsilon})\right).
\end{gather*}
Подставляя полученную формулу в \eqref{ljabel}, получаем
\begin{gather*}
\sum\limits_{a\le R \atop \delta|a}\frac{\lambda(a,\alpha)}{a}=\frac{1}{R}
\sum\limits_{a\le R \atop \delta|a}\frac{4}{\pi}\sigma_{-1}(a)a^{3/2}\sqrt{\alpha}+
O\left(R^{1+\varepsilon}\delta^{\varepsilon}(\alpha^{1+\varepsilon}+\alpha^{-\varepsilon})\right)+\\+
\int_{1}^{R}\Biggl(\frac{1}{t^2}\sum\limits_{a\le t \atop \delta|a}\frac{4}{\pi}\sigma_{-1}(a)a^{3/2}\sqrt{\alpha}+
O\left(t^{\varepsilon}\delta^{\varepsilon}(\alpha^{1+\varepsilon}+\alpha^{-\varepsilon})\right)\Biggl)dt
\end{gather*}
Применяя преобразование Абеля (лемма \ref{abel}), получаем
\begin{gather*}
\sum\limits_{a\le R \atop \delta|a}\frac{\lambda(a,\alpha)}{a}=
\sum\limits_{a\le R \atop \delta|a}\frac{4}{\pi}\sigma_{-1}(a)\sqrt{a}\sqrt{\alpha}+
O\left(R^{1+\varepsilon}\delta^{\varepsilon}(\alpha^{1+\varepsilon}+\alpha^{-\varepsilon})\right).
\end{gather*}
Следовательно, используя формулу \eqref{lj-lj}, получаем
\begin{gather*}
\sum\limits_{a\le R \atop \delta|a}\frac{\lambda^{*}(a,\alpha)}{a}=
\sum_{n|\delta}\sum\limits_{d_1d_2\le R \atop (d_1d_2,\delta)=n}\frac{\mu(d_1)\mu(d_2)}{d_1}
\sum\limits_{a\le \frac{R}{d_1d_2} \atop  \frac{\delta}{n}|a}
\frac{4}{\pi}\sigma_{-1}(a)\sqrt{a}\sqrt{\frac{d_1\alpha}{d_2}}+\\+
O\left(\sum_{n|\delta}\sum\limits_{d_1d_2\le R \atop (d_1d_2,\delta)=n}\frac{1}{d_1}
\frac{R^{1+\varepsilon}}{(d_1d_2)^{1+\varepsilon}}\frac{\delta^{\varepsilon}}{n^{\varepsilon}}
\left(\alpha^{1+\varepsilon}\frac{d_1^{1+\varepsilon}}{d_2^{1+\varepsilon}}+
\alpha^{-\varepsilon}\frac{d_1^{-\varepsilon}}{d_2^{-\varepsilon}}\right)
\right).
\end{gather*}
Совершая в обратной последовательности преобразования, поделанные в формуле \eqref{lj-lj}, получаем
\begin{gather*}
\sum_{n|\delta}\sum\limits_{d_1d_2\le R \atop (d_1d_2,\delta)=n}\frac{\mu(d_1)\mu(d_2)}{d_1}
\sum\limits_{a\le \frac{R}{d_1d_2} \atop  \frac{\delta}{n}|a}
\sigma_{-1}(a)\sqrt{a}\sqrt{\frac{d_1}{d_2}}=\\=
\sum_{n|\delta}\sum\limits_{d_1d_2\le R \atop (d_1d_2,\delta)=n}\frac{\mu(d_1)\mu(d_2)}{d_1}
\sum\limits_{a\le \frac{R}{d_1d_2} \atop  \frac{\delta}{n}d_1d_2|a}
\sigma_{-1}\left(\frac{a}{d_1d_2}\right)\sqrt{\frac{a}{d_1d_2}}\sqrt{\frac{d_1}{d_2}}=
\sum_{\substack{a\le R\\\delta\mid a}}\sqrt{a}\sum_{d_1d_2\mid a}
\frac{\mu(d_1)\mu(d_2)}{d_1d_2}\sigma_{-1}\left(\frac{a}{d_1d_2}\right).
\end{gather*}
Следовательно,
\begin{gather*}
\sum\limits_{a\le R \atop \delta|a}\frac{\lambda^{*}(a,\alpha)}{a}=
\frac{4}{\pi}\sqrt{\alpha}\sum_{\substack{a\le R\\\delta\mid a}}\sqrt{a}\sum_{d_1d_2\mid a}
\frac{\mu(d_1)\mu(d_2)}{d_1d_2}\sigma_{-1}\left(\frac{a}{d_1d_2}\right)+\\+
O(R^{1+\varepsilon}\delta^{\varepsilon}(\alpha^{1+\varepsilon}+\alpha^{-\varepsilon})).
\end{gather*}
Воспользовавшись формулой
\begin{gather*}
\sum_{d_1d_2\mid a}
\frac{\mu(d_1)\mu(d_2)}{d_1d_2}\sigma_{-1}\left(\frac{a}{d_1d_2}\right)=\frac{\varphi(a)}{a}
\end{gather*}
доказанной в статье H. Heilbronn~\cite[\S4]{Heilbronn}, получаем
\begin{gather*}
\sum\limits_{a\le R \atop \delta|a}\frac{\lambda^{*}(a,\alpha)}{a}=
\frac{4}{\pi}\sqrt{\alpha}\sum\limits_{a\le R \atop \delta|a}\frac{\varphi(a)}{\sqrt{a}}+
O(R^{1+\varepsilon}\delta^{\varepsilon}(\alpha^{1+\varepsilon}+\alpha^{-\varepsilon})).
\end{gather*}
Подставляя полученное выражение в лемму \ref{F}, получаем
\begin{gather*}
F=\sum_{a\le x_3N}\sum_{b\le x_1N}\sum\limits_{c\le x_2N \atop (a,b,c)=1}f(a,b,c)=
\sum\limits_{d_1d_2\le x_3N \atop (d_1,d_2)=1}d_1d_2\sum_{\delta_1\le \frac{x_3N}{d_1}}\sum_{\delta_2\le \frac{x_3N}{d_2}}
\frac{\mu(\delta_1)}{\delta_1}\frac{\mu(\delta_2)}{\delta_2}(\delta_1,\delta_2)\\
\int_0^{\frac{x_1N}{d_1}}\int_0^{\frac{x_2N}{d_2}}\Biggl(
\sum\limits_{a\le \frac{x_3N}{d_1d_2} \atop \delta|a}\frac{8}{\pi}\frac{\varphi(a)}{\sqrt{a}}\sqrt{t_1t_2}+
\left(\frac{x_3N}{d_1d_2}\right)^{1+\varepsilon}\delta^{\varepsilon}
O\left(\frac{t_2^{1+\varepsilon}}{t_1^{\varepsilon}}+\frac{t_1^{1+\varepsilon}}{t_2^{\varepsilon}}\right)\Biggl)dt_1dt_2+
O\left(x_1x_2x_3^{2+\varepsilon}N^{4+\varepsilon}\right).
\end{gather*}
Учитывая замечание \ref{G}, получаем
\begin{gather*}
F-G=\sum_{a\le x_3N}\sum_{b\le x_1N}\sum\limits_{c\le x_2N \atop (a,b,c)=1}\left(f(a,b,c)-\frac{8}{\pi}\sqrt{abc}\right)
=\\=
\sum\limits_{d_1d_2\le x_3N \atop (d_1,d_2)=1}d_1d_2\sum_{\delta_1\le \frac{x_3N}{d_1}}\sum_{\delta_2\le \frac{x_3N}{d_2}}
\frac{\mu(\delta_1)}{\delta_1}\frac{\mu(\delta_2)}{\delta_2}(\delta_1,\delta_2)
\int\limits_0^{\frac{x_1N}{d_1}}\int\limits_0^{\frac{x_2N}{d_2}}\left(\frac{x_3N}{d_1d_2}\right)^{1+\varepsilon}\delta^{\varepsilon}
O\left(\frac{t_2^{1+\varepsilon}}{t_1^{\varepsilon}}+\frac{t_1^{1+\varepsilon}}{t_2^{\varepsilon}}\right)dt_1dt_2+\\+
O\left(x_1x_2x_3^{2+\varepsilon}N^{4+\varepsilon}\right)=\\=
x_3^{1+\varepsilon}N^{4+\varepsilon}\sum\limits_{d_1d_2\le x_3N \atop (d_1,d_2)=1}
\left(\frac{x_1^{1-\varepsilon}x_2^{2+\varepsilon}}{d_1^{1+\varepsilon}d_2^{2+\varepsilon}}+
\frac{x_2^{1-\varepsilon}x_1^{2+\varepsilon}}{d_2^{1+\varepsilon}d_1^{2+\varepsilon}}\right)
\sum_{\delta_1\le \frac{x_3N}{d_1}}\sum_{\delta_2\le \frac{x_3N}{d_2}}
\frac{(\delta_1,\delta_2)}{\delta_1\delta_2}\delta^{\varepsilon}+
O\left(x_1x_2x_3^{2+\varepsilon}N^{4+\varepsilon}\right).
\end{gather*}
Так как
$$\delta=HOK\left(\frac{\delta_1}{(\delta_1,d_2)},\frac{\delta_2}{(\delta_2,d_1)}\right)
\le HOK(\delta_1,\delta_2)=\frac{(\delta_1,\delta_2)}{\delta_1\delta_2},$$
то
\begin{gather*}
\sum_{\delta_1\le \frac{x_3N}{d_1}}\sum_{\delta_2\le \frac{x_3N}{d_2}}
\frac{(\delta_1,\delta_2)}{\delta_1\delta_2}\delta^{\varepsilon}\ll
\sum_{d\le\frac{x_3N}{\max(d_1,d_2)}}\sum_{\delta_1\le \frac{x_3N}{dd_1}}\sum_{\delta_2\le \frac{x_3N}{dd_2}}
\frac{1}{(d\delta_1\delta_2)^{1-\varepsilon}}\ll
\frac{x_3^{\varepsilon}N^{\varepsilon}}{(d_1d_2)^{\varepsilon}}.
\end{gather*}
Следовательно,
\begin{gather*}
F-G=O\left(N^{4+\varepsilon}x_1x_2x_3\left(\frac{x_2^{1+\varepsilon}x_3^{\varepsilon}}{x_1^{\varepsilon}}+
\frac{x_1^{1+\varepsilon}x_3^{\varepsilon}}{x_2^{\varepsilon}}+x_3^{1+\varepsilon}
\right)\right).
\end{gather*}
Тем самым теорема \ref{Th2} доказана. Из теоремы \ref{Th2} сразу следует утверждение теоремы \ref{Th1}.
\end{proof}
\end{Th}
\section{Приложение}
Следующая лемма так же общеизвестна (формула суммирования Эйлера-Маклорена)
\begin{Le}
\label{euler}
Пусть  функции $\rho(x)$, $\sigma(x)$ определяются равенствами
$$\rho(x)=\frac{1}{2}-\{x\}, \quad \sigma(x)=\int\limits_0^x\rho(t)dt.$$
Тогда
\begin{enumerate}
  \item[(a)]
  Если $f(x)$ дважды непрерывно дифференцируема на отрезке $[a,b]$, то
$$
\sum_{a<x\le b}f(x)=
\int\limits_a^bf(x)dx+\rho(b)f(b)-\rho(a)f(a)+\sigma(a)f’(a)-\sigma(b)f’(b)+\int\limits_a^b\sigma(x)f"(x)dx.
$$
  \item[(b)] Если $f(x)$ непрерывно дифференцируема на отрезке $[a,b]$, то
$$
\sum_{a<x\le b}f(x)=
\int\limits_a^bf(x)dx+\rho(b)f(b)-\rho(a)f(a)-\int\limits_a^b\rho(x)f’(x)dx.
$$
\end{enumerate}
\begin{proof} см в книге А.А.Карацубы ~\cite{Karazuba}. \end{proof}
\end{Le}
Нам понадобятся следствия леммы \ref{euler}.
\begin{Sl}
\label{ostatok}
Если $f(x)$ непрерывно дифференцируема на отрезке $[a,b]$ и монотонна, то
$$
\sum_{a<x\le b}f(x)=
\int\limits_a^bf(x)dx+O(f(b))+O(f(a)).
$$
\begin{proof}
Действительно, если функция монотонна то
\begin{gather*}
\left|\int\limits_a^b\rho(x)f’(x)dx\right|\ll\int\limits_a^b|f’(x)|dx=O(f(b))+O(f(a)).
\end{gather*}
Теперь утверждение непосредственно следует из пункта (b) леммы \ref{euler}.
\end{proof}
\end{Sl}
\begin{Sl}
\label{sl}
Справедливы следующие асимптотические формулы
\begin{enumerate}
 \item[(a)]
 Для любого натурального $p$ выполнено
  $$\sum_{n\le R}n^p=\frac{R^{p+1}}{p+1}+O\left(R^p\right),$$
 \item[(b)]
 Для любого натурального $p$ выполнено
 $$\sum_{n\le R}n^p\log\frac{R}{n}=\frac{R^{p+1}}{(p+1)^2}+O\left(R^p\log R\right),$$
 \item[(c)]
  $$\sum_{n\le R}n\frac{R-n}{R+n}=\left(\frac{3}{2}-2\log2\right)R^2+O\left(R\right),$$
 \item[(d)]
  $$\sum_{n\le R}n^2\frac{(R-n)^2}{(R+n)^2}=\left(\frac{25}{3}-12\log2\right)R^3+O\left(R^2\right),$$
\item[(e)]
  $$\sum_{n\le R}\frac{(R-n)}{n^2+R^2}=\frac{\pi}{4}-\frac{1}{2}\log2+O\left(\frac{1}{R}\right),$$
\item[(f)]
  $$\sum_{n\le R}n^2\log\left(1+\frac{R}{n}\frac{R-n}{R+n}\right)=
  \left(\frac{5}{6}-\frac{\log 2}{3}-\frac{\pi}{6}\right)R^3+O\left(R^2\log R\right),$$
 \item[(g)]
  $$\sum_{n\le R}n^2\frac{R-n}{R^2+n^2}=
  \left(\frac{1}{2}+\frac{\log 2}{2}-\frac{\pi}{4}\right)R^2+O\left(R\right),$$
\item[(h)]
  $$\sum_{n\le R}n^3\frac{R-n}{R+n}=
  \left(\frac{17}{12}-2\log 2\right)R^4+O\left(R^3\right),$$
\item[(i)]
  $$\sum_{n\le R}n^4\log\left(1+\frac{R}{n}\frac{R-n}{R+n}\right)=
  \left(\frac{\pi}{10}-\frac{3}{20}-\frac{\log 2}{5}\right)R^5+O\left(R^4\log R\right),$$
\item[(j)]
  $$\sum_{n\le R}n^4\frac{R-n}{R^2+n^2}=
  \left(\frac{\pi}{4}-\frac{5}{12}-\frac{\log 2}{2}\right)R^4+O\left(R^3\right),$$
  \item[(k)]
  $$\sum_{n\le R}n^2\frac{R-n}{R+n}=
  \left(2\log 2-\frac{4}{3}\right)R^3+O\left(R^2\right),$$
  \item[(l)]
  $$\sum_{n\le R}n\frac{(R-n)^2}{(R+n)^2}=
  \left(8\log 2-\frac{11}{2}\right)R^2+O\left(R\right),$$
   \item[(m)]
  $$\sum_{n\le R}n\log\left(1+\frac{n}{R}\right)=
  \frac{R^2}{4}+O\left(R\right),$$
  \item[(n)]
  $$\sum_{n\le R}n^2\log\left(1+\frac{n}{R}\right)=
  \left(\frac{2}{3}\log2-\frac{5}{18}\right)R^3+O\left(R^2\right),$$
  \item[(0)]
  $$\sum_{n\le R}\frac{1}{n^2}\log\left(1+\frac{n^2}{R^2}\right)=
  \left(\frac{\pi}{2}-\log2\right)\frac{1}{R}+O\left(\frac{1}{R^2}\right),$$
  \item[(p)]
  $$\sum_{n\le R}\frac{(R-n)^2}{(R+n)^2}=
  \left(3-4\log 2\right)R+O\left(1\right),$$
   \item[(r)]
  $$\sum_{n\le R}n\log\frac{R^2+n^2}{n(R+n)}=
  O\left(R^2\right),$$
  \item[(s)]
  $$\sum_{n\le R}\frac{(R+n)}{n^2+R^2}=O\left(1\right).$$
  \end{enumerate}
\begin{proof}
Все пункты  доказываются непосредственным применением пункта (b) \mbox{леммы \ref{euler}.}
\begin{enumerate}
\item[(b)]
Достаточно применить пункт (b) леммы \ref{euler} и формулу
\begin{gather}\label{nlog}
\int x^p\log x=\frac{x^{p+1}}{p+1}\log x-\frac{x^{p+1}}{(p+1)^2}.
\end{gather}
\item[(e)]
Достаточно применить пункт (b) леммы \ref{euler} и формулу
\begin{gather*}
\int\frac{1}{x^2+a^2}dx=\frac{1}{a}\arctan\frac{x}{a}, \quad
\int\frac{x}{x^2+a^2}dx=\frac{1}{2}\log\left(x^2+a^2\right).
\end{gather*}
  \item[(f)]
Преобразуем суммируемую функцию
\begin{gather*}
n^2\log\left(1+\frac{R}{n}\frac{R-n}{R+n}\right)=
n^2\log\left(R^2+n^2\right)-
n^2\log n-
n^2\log\left(R+n\right).
\end{gather*}
К каждой из трех сумм применим пункт (b) леммы \ref{euler}, используя формулу \eqref{nlog} и формулу
$$
\int x^2\log(x^2+a^2)dx=\frac{1}{3}\left(x^3\log(x^2+a^2)-\frac{2}{3}x^3+2xa^2-2a^3\arctan\frac{x}{a}\right).
$$
\item[(g)]
Преобразуем суммируемую функцию
\begin{gather*}
n^2\frac{R-n}{R^2+n^2}=\left(R-n\right)+R^2\frac{n}{R^2+n^2}-R^3\frac{1}{R^2+n^2}
\end{gather*}
К каждой из трех сумм применим пункт (b) леммы \ref{euler}.
  \item[(i)]
Преобразуем суммируемую функцию
\begin{gather*}
n^4\log\left(1+\frac{R}{n}\frac{R-n}{R+n}\right)=
n^4\log\left(R^2+n^2\right)-
n^4\log n-
n^4\log\left(R+n\right).
\end{gather*}
К каждой из трех сумм применим пункт (b) леммы \ref{euler}, используя формулу \eqref{nlog} и формулу
$$
\int x^4\log(x^2+a^2)dx=\frac{1}{5}\left(x^5\log(x^2+a^2)-\frac{2}{5}x^5+
\frac{2}{3}x^3a^2-2xa^4+2a^5\arctan\frac{x}{a}\right).
$$
\item[(j)]
Преобразуем суммируемую функцию
\begin{gather*}
n^4\frac{R-n}{R^2+n^2}=\left(-n^3+Rn^2+R^2n-R^3\right)+R^4\frac{R-n}{R^2+n^2}.
\end{gather*}
К каждой сумме применим пункт (b) леммы \ref{euler}.
\item[(0)]
Преобразуем суммируемую функцию
\begin{gather*}
\frac{1}{n^2}\log\left(1+\frac{n^2}{R^2}\right)=
\frac{\log(R^2+n^2)}{n^2}-\frac{\log(R^2)}{n^2}.
\end{gather*}
К каждой из двух сумм применим пункт (b) леммы \ref{euler}, используя  формулу
$$
\int\frac{\log(a^2+x^2)}{x^2}dx=-\frac{\log(a^2+x^2)}{x}+2\arctan\frac{x}{a}.
$$

\end{enumerate}

\end{proof}
\end{Sl}

\end{document}